\documentclass[preprint,3p,times]{elsarticle}

\usepackage[section]{placeins}
\usepackage{multirow}
\usepackage{makecell}
\usepackage{longtable}
\usepackage{booktabs}
\usepackage{supertabular}
\usepackage[table]{xcolor}
\usepackage{colortbl}
\usepackage{array}
\usepackage{color}

\usepackage{bm}
\usepackage{amssymb}   
\usepackage{amsmath} 

\usepackage{subeqnarray}
\usepackage{cases}  

\usepackage{arydshln}

\usepackage{tikz}
\usepackage{etoolbox}
\usepackage{enumitem}
\newcommand*{\circled}[1]{\lower.7ex\hbox{\tikz\draw (0pt, 0pt)%
    circle (.5em) node {\makebox[1em][c]{\small #1}};}}
\robustify{\circled}
\newcommand*{\blackcircled}[1]{{\lower.7ex\hbox\tikz\draw[fill=black] (0pt, 0pt)%
    circle (.5em) node {\makebox[1em][c]{\small\color{white} #1}};}}
\robustify{\blackcircled}

\usepackage{graphicx}
\usepackage{epstopdf}

\usepackage[bf]{caption2}

\usepackage[ruled,linesnumbered,lined,commentsnumbered]{algorithm2e}

\usepackage[figuresright]{rotating}

\newtheorem{theorem}{Theorem}
\newtheorem{lemma}{Lemma}

\newtheorem{definition}{Definition}

\newtheorem{remark}{Remark}
\newtheorem{example}{Example}

\journal{Elsevier}

\begin{document}

\begin{frontmatter}

\title{Locally Order-Preserving Mapping for WENO Methods}


\author[a]{Ruo Li}
\ead{rli@math.pku.edu.cn}

\author[b,c]{Wei Zhong\corref{cor1}}
\ead{zhongwei2016@pku.edu.cn}

\cortext[cor1]{Corresponding author}

\address[a]{CAPT, LMAM and School of Mathematical Sciences, Peking
University, Beijing 100871, China}

\address[b]{School of Mathematical Sciences, Peking University,
Beijing 100871, China}

\address[c]{Northwest Institute of Nuclear Technology, Xi'an 
710024, China}

\begin{abstract}

 In our previous studies \cite{MOP-WENO-ACMk,MOP-WENO-X}, the commonly reported issue that most of the existing mapped WENO schemes suffer from either losing high resolutions or generating spurious oscillations in long-run simulations of hyperbolic problems has been successfully addressed, by devising the improved mapped WENO schemes, namely MOP-WENO-X, where ``X'' stands for the version of the existing mapped WENO scheme. However, all the MOP-WENO-X schemes bring about the serious deficiency that their resolutions in the region with high-frequency but smooth waves are dramatically decreased compared to their associated WENO-X schemes. The purpose of this paper is to overcome this drawback. We firstly present the definition of the \textit{locally order-preserving (LOP)} mapping. Then, by using a new proposed posteriori adaptive technique, we apply this \textit{LOP} property to obtain the new mappings from those of the WENO-X schemes. The essential idea of the posteriori adaptive technique is to identify the global stencil in which the existing mappings fail to preserve the \textit{LOP} property, and then replace the mapped weights with the weights of the classic WENO-JS scheme to recover the \textit{LOP} property. We build the resultant mapped WENO schemes and denote them as LOP-WENO-X. The numerical experiments demonstrate that the resolutions in the region with high-frequency but smooth waves of the LOP-WENO-X schemes are similar or even better than those of their associated WENO-X schemes and naturally much higher than the MOP-WENO-X schemes. Furthermore, the LOP-WENO-X schemes gain all the great advantages of the MOP-WENO-X schemes, such as attaining high resolutions and in the meantime preventing spurious oscillations near discontinuities when solving the one-dimensional linear advection problems with long output times, and significantly reducing the post-shock oscillations in the simulations of the two-dimensional problems with shock waves.

\end{abstract}



\end{frontmatter}

\section{Introduction}
\label{secIntroduction}
In the numerical calculation of hyperbolic conservation laws in the 
form
\begin{equation}
\dfrac{\partial \mathbf{u}}{\partial t} +
\nabla \cdot \mathbf{F}(\mathbf{u}) = 0,
\label{eq:governingEquation}
\end{equation}
with proper initial conditions and boundary conditions, the class of 
weighted essentially non-oscillatory (WENO) schemes 
\cite{ENO-Shu1988,ENO-Shu1989,WENO-LiuXD,WENO-JS,WENOoverview}, 
which is an evolution of the essentially non-oscillatory (ENO) 
schemes \cite{ENO1987JCP71, ENO1987V24, ENO1986, ENO1987JCP83}, is 
widely used due to its success in using a nonlinear convex 
combination of all the candidate stencils to automatically achieve 
high-order accuracy in smooth regions and without destroying the 
non-oscillatory property near shocks or other discontinuities. The 
classic WENO-JS proposed by Jiang and Shu \cite{WENO-JS} is the most 
popular one of the WENO schemes. By using the sum of the normalized 
squares of the scaled $L_{2}-$norms of all the derivatives of $r$ 
local interpolating polynomials, a new measurement of the smoothness 
of the numerical solutions on substencils, named smoothness 
indicators, is devised to help to obtain $(2r-1)$th-order of 
accuracy from the $r$th-order ENO schemes.

It is well-known that the WENO-JS scheme is a quite successful 
methodology for solving problems modeled by the hyperbolic 
conservation laws in the form of Eq.\eqref{eq:governingEquation}. 
However, despite its great advantages of efficient 
implementation and high-order accuracy, its convergence orders 
dropped for many cases such as at or near critical points of order 
$n_{\mathrm{cp}}=1$ in the smooth regions. Here, we refer to
$n_{\mathrm{cp}}$ as the order of the critical point; e.g.,
$n_{\mathrm{cp}} = 1$ corresponds to $f'=0, f'' \neq 0$ and
$n_{\mathrm{cp}} = 2$ corresponds to $f'=0, f'' = 0, f''' \neq 0$, 
etc. In \cite{WENO-M}, Henrick et al. pointed out that the 
fifth-order WENO-JS scheme fails to yield the 
optimal convergence order at or near critical points where the first 
derivative vanishes but the third derivative does not 
simultaneously. In the same article, they derived the necessary and 
sufficient conditions on the nonlinear weights for optimality of the 
convergence rate of the fifth-order WENO schemes and these 
conditions were reduced to a simpler sufficient condition by Borges 
et al. \cite{WENO-Z}, which could be easily extended to the 
$(2r-1)$th-order WENO schemes \cite{WENO-IM}. Then, in order to 
address the drawback of the WENO-JS scheme discussed above, Henrick 
et al. \cite{WENO-M} innovatively designed a mapping function 
satisfying the sufficient condition to achieve the optimal order of 
accuracy, leading to the original mapped WENO scheme, named WENO-M 
hereafter. Since then, obeying the similar criteria proposed by 
Henrick et al. \cite{WENO-M}, many different kinds of mapped WENO 
schemes have been successfully proposed \cite{WENO-PM,WENO-IM,
WENO-PPM5,WENO-RM260,WENO-MAIMi,WENO-ACM}. 

In \cite{WENO-IM}, by rewriting the mapping function of the WENO-M 
scheme in a simpler and more meaningful form and then extending it 
to a general class of improved mapping functions, Feng et al. 
proposed the group of WENO-IM($k, A$) schemes and WENO-IM(2, 0.1) 
was recommended. WENO-IM(2, 0.1) can significantly decrease the 
dissipations of the WENO-M scheme leading to higher resolutions.
However, in \cite{WENO-RM260}, Wang et al. indicated that the 
seventh- and ninth- order WENO-IM(2, 0.1) schemes generated evident 
spurious oscillations near discontinuities with a long output time. 
Furthermore, the present authors have reported that 
\cite{WENO-MAIMi,MOP-WENO-ACMk}, even for the fifth-order 
WENO-IM(2,0.1) scheme, the spurious oscillations are also produced 
when the grid number increases. In \cite{WENO-PM}, Feng et al. found 
that, when the WENO-M scheme was used for solving the problems with 
discontinuities, its mapping function may amplify the effect from 
the non-smooth stencils leading to a potential loss of accuracy near 
discontinuities and this loss of accuracy could be accumulated when the simulating time is large. To fix this issue, two additional 
requirements, that is, $g'(0) = 0$ and $g'(1)=0$ ($g(x)$ denotes the 
mapping function), to the original criteria in \cite{WENO-M} was 
proposed and then a piecewise polynomial mapping function satisfying 
these additional requirements was devised. The resultant scheme was 
denoted as WENO-PM$k$ and $k = 6$ was recommended. The WENO-PM6 
scheme \cite{WENO-PM} obtained significantly 
higher resolution than the WENO-M scheme when computing the 
one-dimensional linear advection problem with long output times. 
However, it may generate the non-physical oscillations near the 
discontinuities as shown in Fig. 8 of \cite{WENO-IM} and Figs. 3-8 
of \cite{WENO-RM260}. 

Besides WENO-IM($2, 0.1$) and WENO-PM6, many other modified mapped 
WENO schemes have been successfully proposed to enhance the 
conventional WENO-JS scheme's performance, e.g., WENO-PPM$n$
\cite{WENO-PPM5}, WENO-RM($mn0$) \cite{WENO-RM260}, WENO-MAIM$i$
\cite{WENO-MAIMi}, WENO-ACM \cite{WENO-ACM}, MIP-WENO-ACM$k$ 
\cite{MOP-WENO-ACMk} and et al. Despite that, as reported in 
literatures \cite{WENO-IM,WENO-RM260}, most of these existing 
modified mapped WENO schemes can hardly avoid the spurious 
oscillations near discontinuities, especially for long output time 
simulations. In addition, when computing the 2D problems with shock 
waves, the post-shock oscillations become very serious for most of 
the existing modified mapped WENO schemes \cite{WENO-ACM}.

It was reported \cite{MOP-WENO-ACMk} that, for many existing mapped 
WENO schemes, e.g., WENO-PM6 \cite{WENO-PM}, WENO-IM(2, 0.1) 
\cite{WENO-IM}, WENO-MAIM1 \cite{WENO-MAIMi}, MIP-WENO-ACM$k$ 
\cite{MOP-WENO-ACMk} and et al, the order of the nonlinear 
weights for the substencils of the same global stencil has been 
changed at many points in the mapping process. This is caused by 
weights increasing of non-smooth substencils and weights decreasing 
of smooth substencils. As far as it is known, this phenomenon occurs 
in all existing mapped WENO schemes. After a systematic theoretical 
analysis and a further verification with extensive numerical 
experiments, the authors claimed that the order-change of the 
mapped nonlinear weights may essentially cause the resolution loss 
or generate the non-physical numerical oscillations by existing 
mapped WENO schemes, when making long output time 
simulations. Then, the concept of \textit{order-preserving} mapping 
has been defined and the \textit{order-preserving} property was 
suggested as an additional criterion in the design of the mapping 
function. Following this criterion, the new mapped WENO scheme, say, 
MOP-WENO-ACM$k$, was proposed. It was examined by numerical tests 
that the MOP-WENO-ACM$k$ scheme can obtain the optimal convergence 
rates in smooth regions even in the presence of critical points. And 
also, it is able to remove spurious oscillations around 
discontinuities and to reduce the numerical dissipations so that the 
resolutions are very high for long output times. Furthermore, the 
MOP-WENO-ACM$k$ scheme has a significant advantage in decreasing the 
post-shock oscillations when solving the 2D tests with shock waves. 
Lately, the idea of \textit{order-preserving} mapping was 
successfuly introduced into other existing mapped WENO schemes and 
the resultant improved mapped schemes \cite{MOP-WENO-X}, say, 
MOP-WENO-X, gained all the benefits of the MOP-WENO-ACM$k$ scheme. 
However, disappointedly, all the MOP-WENO-X schemes including 
MOP-WENO-ACM$k$ fail to achieve high resolutions in the region with 
high-frequency but smooth waves, such as the Shu-Osher problem 
\cite{ENO-Shu1989} and the Titarev-Toro problem \cite{Titarev-Toro-1,
Titarev-Toro-2,Titarev-Toro-3}. Indeed, their resolutions are even 
much lower than the associated WENO-X schemes (see subsection 
\ref{subsec:1DEuler} below).

Our major purpose in this study is to address the aforementioned 
shortcoming of the MOP-WENO-X schemes while maintaining their 
benefits. We propose the definition of the 
\textit{locally order-preserving (LOP)} mapping, which is a 
development of the \textit{order-preserving (OP)} mapping given in 
\cite{MOP-WENO-ACMk}. By using a posteriori adaptive technique, we 
apply the \textit{LOP} property to various existing mapped WENO 
schemes leading to a new class of mapped WENO schemes, denoted as 
LOP-WENO-X. Firstly, a new function named \textbf{postINDEX} used 
to implement the posteriori adaptive technique is defined (see 
Definition \ref{def:postINDEX} in subsection \ref{subsec:newMapping} 
below). Then, a general algorithm to construct \textit{LOP} 
mappings based on the existing mappings by using the 
posteriori adaptive technique is proposed. We present the 
properties and the necessary proofs or analyses of the mappings of 
the LOP-WENO-X schemes. The convergence rates of accuracy of the 
LOP-WENO-X schemes have also been given. Solutions for 1D linear 
advection problems with initial conditions including high-order 
critical points and discontinuities at large output times have been 
discussed in detail. We demonstrate the great advantages of 
the LOP-WENO-X schemes in the region with high-frequency but smooth 
waves by solving the Shu-Osher and Titarev-Toro problems. At last, 
for 2D Euler equations, numerical experiments of accuracy tests and 
a benchmark problem with shock waves, are run to show the good
performances of the LOP-WENO-X schemes.

We organize the remainder of this paper as follows. In Section 
\ref{secMappedWENO}, we briefly review the preliminaries to 
understand the procedures of the WENO-JS \cite{WENO-JS}, WENO-M 
\cite{WENO-M} and some other versions of mapped WENO schemes. The 
main contribution of this paper will be presented in Section 
\ref{LOP-WENO-X}, where we provide the posteriori adaptive 
technique to build a general method to introduce the \textit{locally 
order-preserving} mapping and hence derive the LOP-WENO-X schemes 
for improving the existing mapped WENO-X schemes. Some numerical 
results of 2D Euler equations are provided in Section 
\ref{NumericalExperiments} to illustrate the performance and 
advantages of the proposed WENO schemes. Finally, we close this 
paper with concluding remarks in Section \ref{secConclusions}.


\section{Preliminaries}
\label{secMappedWENO}
\subsection{The fifth-order WENO-JS scheme}
For the hyperbolic conservation laws in 
Eq.\eqref{eq:governingEquation}, without loss of generality, we 
discuss its simplest form of the one-dimensional scalar equation 
\begin{equation}
u_{t} + f(u)_{x} = 0.
\label{eq:1D-hyperbolicLaw}
\end{equation}
Let $\{I_{j}\}$ be a control volume of the given computational 
domain $[x_{l}, x_{r}]$ with the $j$th cell $I_{j}:=[x_{j-1/2}, 
x_{j+1/2}]$. The center and boundaries of $I_{j}$ are denoted by 
$x_{j}=x_{l} + (j - 1/2)\Delta x$ and $x_{j \pm 1/2} = x_{j} \pm 
\Delta x/2$ with the cell size $\Delta x = \frac{x_{r} - x_{l}}{N}$ 
leading to the uniform meshes. Let $\bar{u}_{j}(t)$ be the numerical 
approximation to the cell average $\bar{u}(x_{j}, t)=
\dfrac{1}{\Delta x}\int_{x_{j-1/2}}^{x_{j+1/2}}u(\xi,t)\mathrm{d}\xi$
, then the semi-discretization form of Eq.\eqref{eq:1D-hyperbolicLaw}
can be written as
\begin{equation}
\dfrac{\mathrm{d}\bar{u}_{j}(t)}{\mathrm{d}t}
\approx - \dfrac{1}{\Delta x}\Big( \hat{f}_{j+1/2} - \hat{f}_{j-1/2} 
\Big),
\label{eq:discretizedFunction}
\end{equation}
where $\hat{f}_{j \pm 1/2} = \hat{f}(u_{j \pm 1/2}^{-}, 
u_{j \pm 1/2}^{+})$ is the numerical flux used to approximate the 
physical flux function $f(u)$ at the cell boundaries $x_{j \pm 1/2}$
. In this paper, the values of $u_{j \pm 1/2}^{\pm}$ are calculated 
by the WENO reconstructions narrated later, and hereafter, we will 
only describe how $u_{j + 1/2}^{-}$ is approximated as the formulas 
for $u_{j + 1/2}^{+}$ are symmetric to $u_{j + 1/2}^{-}$ with 
respect to $x_{j+1/2}$. Also, for brevity, we will drop the ``-'' 
sign in the superscript. 

In the fifth-order WENO-JS scheme, a 5-point global stencil 
$S^{5}=\{I_{j-2}, I_{j-1}, I_{j}, I_{j+1},I_{j+2}\}$ is used to 
construct the values of $u_{j + 1/2}$ from known cell average values 
$\overline{u}_{j}$. The global stencil is subdivided into three 
3-point substencils $S_{s} = \{I_{j+s-2}, I_{j+s-1}, I_{j+s}\}$ with 
$s=0, 1, 2$. It is known that the third-order approximations of 
$u(x_{j+1/2}, t)$ associated with these substencils are explicitly 
given by
\begin{equation}
\begin{array}{l}
\begin{aligned}
&u_{j+1/2}^{0} = \dfrac{1}{3}\bar{u}_{j-2} - 
\dfrac{7}{6}\bar{u}_{j-1}
+ \dfrac{11}{6}\bar{u}_{j}, \\
&u_{j+1/2}^{1} = - \dfrac{1}{6}\bar{u}_{j-1} + 
\dfrac{5}{6}\bar{u}_{j}
+ \dfrac{1}{3}\bar{u}_{j+1}, \\
&u_{j+1/2}^{2} = \dfrac{1}{3}\bar{u}_{j} + 
\dfrac{5}{6}\bar{u}_{j+1}
- \dfrac{1}{6}\bar{u}_{j+2}.
\end{aligned}
\end{array}
\label{eq:approx_ENO}
\end{equation}
Through a convex combination of those third-order approximations of 
substencils, the $u_{j + 1/2}$ of global stencil $S^{5}$ is computed 
as follows
\begin{equation}
u_{j + 1/2} = \sum\limits_{s = 0}^{2}\omega_{s}u_{j + 1/2}^{s}.
\label{eq:approx_WENO}
\end{equation}
The nonlinear weights of WENO-JS is calculated by
\begin{equation} 
\omega_{s}^{\mathrm{JS}} = \dfrac{\alpha_{s}^{\mathrm{JS}}}{\sum_{l =
 0}^{2} \alpha_{l}^{\mathrm{JS}}}, \alpha_{s}^{\mathrm{JS}} = \dfrac{
 d_{s}}{(\epsilon + \beta_{s})^{2}}, \quad s = 0,1,2,
\label{eq:weights:WENO-JS}
\end{equation} 
where the ideal linear weights $d_{0}=0.1, d_{1}=0.6, d_{2} = 0.3$, 
$\epsilon > 0$ is a very small number so that the denominator will 
not be zero, and the smoothness indicators $\beta_{s}$ are given as 
\cite{WENO-JS}
\begin{equation*}
\begin{array}{l}
\begin{aligned}
\beta_{0} &= \dfrac{13}{12}\big(\bar{u}_{j - 2} - 2\bar{u}_{j - 1} + 
\bar{u}_{j} \big)^{2} + \dfrac{1}{4}\big( \bar{u}_{j - 2} - 4\bar{u}_
{j - 1} + 3\bar{u}_{j} \big)^{2}, \\
\beta_{1} &= \dfrac{13}{12}\big(\bar{u}_{j - 1} - 2\bar{u}_{j} + \bar
{u}_{j + 1} \big)^{2} + \dfrac{1}{4}\big( \bar{u}_{j - 1} - \bar{u}_{
j + 1} \big)^{2}, \\
\beta_{2} &= \dfrac{13}{12}\big(\bar{u}_{j} - 2\bar{u}_{j + 1} + \bar
{u}_{j + 2} \big)^{2} + \dfrac{1}{4}\big( 3\bar{u}_{j} - 4\bar{u}_{j 
+ 1} + \bar{u}_{j + 2} \big)^{2}.
\end{aligned}
\end{array}
\end{equation*}

The WENO-JS scheme fails to obtain the designed convergence rates of 
accuracy at or near critical points. More details can be 
found in \cite{WENO-M}.

\subsection{The mapped WENO approach}
\label{subsecWENO-M}
To recover the designed convergence order at or near critical points,
Henrich et al. \cite{WENO-M} designed a mapping function taking the 
form
\begin{equation}
\big( g^{\mathrm{M}} \big)_{s}(\omega) = \dfrac{ \omega \big( d_{s} +
d_{s}^2 - 3d_{s}\omega + \omega^{2} \big) }{ d_{s}^{2} + (1 - 2d_
{s})\omega }, \quad \quad s = 0, 1, 2.
\label{mapFuncWENO-M}
\end{equation}
One can easily verify that $\big( g^{\mathrm{M}}\big)_{s}(\omega)$ 
is a non-decreasing monotone function on $[0, 1]$ with finite slopes 
and satisfies the following properties.
\begin{lemma} 
The mapping function $\big( g^{\mathrm{M}} \big)_{s}(\omega)$
defined by Eq.(\ref{mapFuncWENO-M}) satisfies: \\

C1. $0 \leq \big(g^{\mathrm{M}}\big)_{s}(\omega)\leq 1, \big(
g^{\mathrm{M}}\big)_{s}(0)=0, \big(g^{\mathrm{M}}\big)_{s}(1) = 1$;

C2. $\big( g^{\mathrm{M}} \big)_{s}(d_{s}) = d_{s}$;

C3. $\big( g^{\mathrm{M}} \big)_{s}'(d_{s}) = \big(g^{\mathrm{M}}
\big)_{s}''(d_{s}) = 0$. 
\label{lemmaWENO-Mproperties}
\end{lemma}

After Henrick's innovative work, many other mapping functions have 
been successfully designed \cite{WENO-IM,WENO-PM,WENO-PPM5,
WENO-RM260,WENO-MAIMi,WENO-ACM,MOP-WENO-ACMk}. Here, we directly 
express some mapping functions of these schemes succinctly as shown 
in Table \ref{tab:map:WENO-X}. For more details and other versions 
of the mapping function, we refer to the references.

\begin{sidewaystable}[!ht]
\footnotesize
\centering
\caption{Mapping functions for different mapped WENO schemes.}
\label{tab:map:WENO-X}
\begin{tabular*}{\hsize}
{@{}@{\extracolsep{\fill}}llllll@{}}
\toprule
Scheme, WENO-X  & $\big(g^{\mathrm{X}}\big)_{s}(\omega)$ 
& Coefficients  & Parameter settings  & $n_{\mathrm{X}}$ \\
\hline
{} & {} & {} & {} & {} & {} \\
WENO-JS, \cite{WENO-JS} & $\big(g^{\mathrm{JS}}\big)_{s}(\omega) = 
\omega$ & None         & None         & -    \\
{} & {} & {} & {} & {} & {} \\
WENO-M, \cite{WENO-M}   & $\big( g^{\mathrm{M}} \big)_{s}(\omega) = 
\dfrac{\omega \big( d_{s} + d_{s}^2 - 3d_{s}\omega + \omega^{2} \big) }{ d_{s}^{2} + (1 - 2d_{s})\omega }$ & None    & None       & 2 \\
{} & {} & {} & {} & {} & {} \\
WENO-PM$k$, \cite{WENO-PM}  & $\big( g^{\mathrm{PM}} \big)_{s}(\omega) = C_{1}(\omega - d_{s})^{k+1}(\omega + C_{2}) + d_{s}$ 
& $\big(C_{1}, C_{2}\big) = \left\{
\begin{array}{ll}
\begin{aligned}
&\Bigg((-1)^{k}\dfrac{k+1}{d_{s}^{k+1}}, \dfrac{d_{s}}{k+1}\Bigg), & 0 \leq \omega \leq d_{s}, \\
&\Bigg(-\dfrac{k+1}{(1-d_{s})^{k+1}}, \dfrac{d_{s}-(k+2)}{k+1}\Bigg),    & d_{s} < \omega \leq 1,
\end{aligned}
\end{array}\right.$
& $k = 6$ & $k$ \\
{} & {} & {} & {} & {} & {} \\
WENO-IM($k,A$), \cite{WENO-IM} & $\big( g^{\mathrm{IM}} \big)_{s}(\omega) = d_{s} + \dfrac{\big(\omega - d_{s} \big)^{k + 1}A}{\big( \omega - d_{s} \big)^{k}A + \omega(1 - \omega)}, \quad A > 0, k = 2n, n \in \mathbb{N}^{+}$ & None & $\left\{\begin{array}{l} k=2,\\
 A = 0.1. \end{array}\right.$  & $k$  \\ 
{} & {} & {} & {} & {} & {} \\
WENO-PPM$n$, \cite{WENO-PPM5}  & $\big( g_{s}^{\mathrm{PPM}5} \big)_{s}(\omega) = \left\{
\begin{array}{ll}
\begin{aligned}
&\big( g_{s,\mathrm{L}}^{\mathrm{PPM}5} \big)_{s}(\omega) = d_{s}
\big( 1 + (a - 1)^{5} \big), & \omega 
\in [0, d_{s}]\\
&\big( g_{s,\mathrm{R}}^{\mathrm{PPM}5} \big)_{s}(\omega) = d_{s} + 
b^{4}\big( \omega - d_{s} \big)^{5}, & \omega 
\in (d_{s}, 1],
\end{aligned}
\end{array}
\right. $  & $\left\{\begin{array}{l} a = \dfrac{\omega}{d_{s}},\\
b = \dfrac{1}{d_{s} - 1} \end{array}\right.$ & $n = 5$  
& $4$  \\
{} & {} & {} & {} & {} & {} \\
WENO-RM($mn0$), \cite{WENO-RM260} & $\big( g^{\mathrm{RM}} \big)_{s}(\omega) = d_{s} + \dfrac{(\omega - 
d_{s})^{7}}{a_{0} + a_{1}\omega + a_{2}\omega^{2} + a_{3}\omega^{3}}$ & $\left\{\begin{array}{l} a_{0} = d_{s}^{6}, \\ 
a_{1} = -7d_{s}^{5}, \\ 
a_{2} = 21d_{s}^{4}, \\ 
a_{3} = (1 - d_{s})^{6} - \sum\limits_{i = 0}^{2}a_{i}. \end{array}\right.$ & $\left\{\begin{array}{l} m = 2,\\ n= 6. \end{array}\right.$
& $3, 4$  \\
{} & {} & {} & {} & {} & {} \\
WENO-ACM, \cite{WENO-ACM}   & $\big( g^{\mathrm{ACM}} \big)_{s}(\omega) = \left\{
\begin{array}{ll}
\begin{aligned}
&\dfrac{d_{s}}{2}\mathrm{sgm}(\omega - \mathrm{CFS}_{s}, \delta_{s}, 
A, k) + \dfrac{d_{s}}{2}, & \omega \leq d_{s}, \\
&\dfrac{1-d_{s}}{2}\mathrm{sgm}(\omega -\overline{\mathrm{CFS}}_{s}, 
\delta_{s}, A, k) + \dfrac{1 + d_{s}}{2}, & \omega > d_{s},
\end{aligned}
\end{array}
\right.$ & $\mathrm{sgm}\big( x, \delta, A, k \big) = \left\{ 
\begin{array}{ll}
\begin{aligned}
&\dfrac{x}{|x|}, & |x| \geq \delta, \\ 
&\dfrac{x}{\Big(A\big( \delta ^2 - x^2 \big)\Big)^{k + 3} + \ |x|}, 
& |x| < \delta.
\end{aligned}
\end{array} \right.$ & $\left\{\begin{array}{l}B = 20,\\
 k = 2, \\
 \delta_{s} = 1\mathrm{e-}6, \\
\mathrm{CFS}_{s} = d_{s}/10. \end{array}\right.$ 
& $\infty$ \\
\bottomrule
\end{tabular*}
\end{sidewaystable}

\clearpage
\subsection{Time discretization}
In order to advance the ODEs (see Eq.\eqref{eq:discretizedFunction}) 
resulting from the semi-discretized PDEs in time, we use the 
following explicit, third-order, strong stability preserving (SSP) 
Runge-Kutta method \cite{ENO-Shu1988,SSPRK1998,SSPRK2001}
\begin{equation*}
\begin{array}{l}
\begin{aligned}
&\overrightarrow{U}^{(1)} = \overrightarrow{U}^{n} + \Delta t 
\mathcal{L}(\overrightarrow{U}^{n}), \\
&\overrightarrow{U}^{(2)} = \dfrac{3}{4} \overrightarrow{U}^{n} + 
\dfrac{1}{4} \overrightarrow{U}^{(1)} + \dfrac{1}{4}\Delta t 
\mathcal{L}(\overrightarrow{U}^{(1)}), \\
&\overrightarrow{U}^{n + 1} = \dfrac{1}{3} \overrightarrow{U}^{n} + 
\dfrac{2}{3}\overrightarrow{U}^{(2)} + \dfrac{2}{3} \Delta t 
\mathcal{L}(\overrightarrow{U}^{(2)}),
\end{aligned}
\end{array}
\end{equation*}
where $\mathcal{L}(\cdot) := -\dfrac{1}{\Delta x}\Big(\hat{f}_{j+1/2}
- \hat{f}_{j-1/2} \Big)$, $\overrightarrow{U}^{(1)}$, 
$\overrightarrow{U}^{(2)}$ are the intermediate stages, 
$\overrightarrow{U}^{n}$ is the value of $\overrightarrow{U}$ at 
time level $t^{n} = n\Delta t$, and $\Delta t$ is the time step 
satisfying some proper CFL condition. The WENO reconstructions will 
be applied to compute $\mathcal{L}(\cdot)$. The well-known global 
Lax-Friedrichs flux, that is, $\hat{f}(a,b) = \frac{1}{2}\big[f(a) + 
f(b) - \alpha(b - a)\big]$, will be employed.


\section{The locally order-preserving (LOP) mapped WENO schemes}
\label{LOP-WENO-X}

\subsection{Definition of the locally order-preserving (LOP) mapping}
\label{subsec:LOP}
In \cite{MOP-WENO-ACMk}, the authors innovatively proposed the 
definition of the \textit{order-preserving (OP)} and 
\textit{non-order-preserving (non-OP)} mapping and claimed that the 
\textit{OP} property plays an essential role in obtaining high 
resolution and avoiding spurious oscillations meanwhile for long 
output time simulations. However, the requirement, that is to make 
sure the mapping functions to be \textit{OP} in the whole range of 
$\omega \in (0, 1)$, is a sufficient, but not a necessary, condition 
for the low dissipation and robustness. Actually, this requirement 
is too strict in some sense. Therefore, we develop the 
\textit{locally order-preserving (LOP)} mapping.

\begin{definition}(locally order-preserving mapping) 
For $\forall x_{j}$, let $S^{2r-1}$ denote the $(2r-1)$-point global 
stencil centered around $x_{j}$. Assume that $\{\omega_{0},\cdots,
\omega_{r-1}\}$ are the nonlinear weights associated with the 
$r$-point substencils $\{S_{0},\cdots,S_{r-1}\}$, and
$\big(g^{\mathrm{X}}\big)_{s}(\omega), s=0,\cdots,r-1$ is the 
mapping function of the mapped WENO-X scheme. If for $\forall 
m, n \in \{0, \cdots, r - 1\}$, when $\omega_{m} > \omega_{n}$, we 
have $\big( g^{\mathrm{X}}\big)_{m}(\omega_{m}) \geq 
\big(g^{\mathrm{X}} \big)_{n}(\omega_{n})$, and when $\omega_{m}=
\omega_{n}$, we have $\big( g^{\mathrm{X}}\big)_{m}(\omega_{m}) = 
\big(g^{\mathrm{X}} \big)_{n}(\omega_{n})$, then we say the set of 
mapping functions \Big\{$\big( g^{\mathrm{X}}\big)_{s}(\omega), s=0,
\cdots,r-1$\Big\} is \textbf{locally order-preserving (LOP)}.
\label{def:LOPM}
\end{definition}

To maintain coherence and for the convenience of the readers, we 
state the definition of \textit{OP/non-OP} point proposed in 
\cite{MOP-WENO-ACMk}.

\begin{definition}(OP/non-OP point) We say that a \textbf{non-OP} 
mapping process occurs at $x_{j}$, if $\exists m, n \in \{0,\cdots,
r-1\}$, s.t.
\begin{equation}\left\{
\begin{array}{ll}
\begin{aligned}
&\big(\omega_{m} - \omega_{n}\big)\bigg(\big(g^{\mathrm{X}}\big)_{m}
(\omega_{m}) - \big(g^{\mathrm{X}}\big)_{n}(\omega_{n})\bigg) < 0, 
&\mathrm{if} \quad \omega_{m} \neq \omega_{n},\\
&\big(g^{\mathrm{X}}\big)_{m}(\omega_{m}) \neq \big(g^{\mathrm{X}}
\big)_{n}(\omega_{n}), &\mathrm{if} \quad \omega_{m}=\omega_{n}.
\end{aligned}
\end{array}\right.
\end{equation}
And we say $x_{j}$ is a \textbf{non-OP point}. Otherwise, we say 
$x_{j}$ is an \textbf{OP point}.
\label{def:OP-point}
\end{definition}

\begin{remark}
Naturally, if the set of mapping functions 
\Big\{$\big( g^{\mathrm{X}}\big)_{s}(\omega), s = 0,\cdots,r-1$\Big\}
is not \textbf{LOP}, it must be \textbf{non-OP}.
\end{remark}

\subsection{Design of the LOP mapped WENO schemes}
\label{subsec:newMapping}
For illustrative purposes in the present study we mainly consider a 
limited number of existing WENO schemes as shown in Table 
\ref{tab:map:WENO-X}, where we present their setting parameters. The 
notation $n_{\mathrm{X}}$ denotes the order of the specified 
critical point, namely $\omega = d_{s}$, of the mapping function of 
the WENO-X scheme, that is, $\big( g^{\mathrm{X}} \big)'_{s}(d_{s}) =
\cdots = \big( g^{\mathrm{X}} \big)^{(n_{\mathrm{X}})}_{s}(d_{s})=0, 
\big( g^{\mathrm{X}} \big)^{(n_{\mathrm{X}} + 1)}_{s}(d_{s})\neq 0$. 
To simplify the presentation below, we have already presented 
$n_{\mathrm{X}}$ of the WENO-X scheme in the last column of Table 
\ref{tab:map:WENO-X}. 

\begin{lemma}
For the WENO-X scheme shown in Table \ref{tab:map:WENO-X}, the 
mapping function $\big( g^{\mathrm{X}} \big)_{s}(\omega), s=0,1,
\cdots,r-1$ is monotonically increasing over $[0, 1]$.
\label{lemma:general_mapping}
\end{lemma}
\textbf{Proof.}
See the references given in the last column of Table 
\ref{tab:map:WENO-X}.
$\hfill\square$ \\

Before proposing Algorithm \ref{alg:posteriori} to 
devise the posteriori adaptive \textit{OP} mapping, we firstly give 
the \textbf{postINDEX} function and a set of function 
$\mathbb{S}^{\mathrm{X}}$ by the following definitions.

\begin{definition}({\rm{\textbf{postINDEX}}} function)
The {\rm{\textbf{postINDEX}}} function is defined as follows
\begin{equation}
\mathbf{postINDEX}(a, b, \mathrm{X}) = 
\Big( \omega_{a}^{\mathrm{JS}} - \omega_{b}^{\mathrm{JS}} 
\Big)\Bigg( g_{a}^{\mathrm{X}}\Big(\omega_{a}^{\mathrm{JS}}\Big) - 
g_{b}^{\mathrm{X}}\Big(\omega_{b}^{\mathrm{JS}}\Big) \Bigg),
\label{eq:postINDEX}
\end{equation}
where $a, b = 0, \cdots, r-1$, $\omega_{a}^{\mathrm{JS}},
\omega_{b}^{\mathrm{JS}}$ are the nonlinear weights of the WENO-JS 
scheme, and $g_{s}^{\mathrm{X}}(\omega)$ is the mapping function of 
the existing mapped WENO-X scheme as shown in Table \ref{tab:map:WENO-X}.
\label{def:postINDEX}
\end{definition}

\begin{definition} Define a set of function $\mathbb{S}^{\mathrm{X}}$
as follows
\begin{equation}
\mathbb{S}^{\mathrm{X}}=\bigg\{\mathbf{postINDEX}(a,b,\mathrm{X}):   
\mathbf{postINDEX}(a, b, \mathrm{X}) > 0\bigg\} \bigcup \bigg\{ 
\mathbf{postINDEX}(a, b, \mathrm{X}): \omega_{a}^{\mathrm{JS}} - 
\omega_{b}^{\mathrm{JS}} = g_{a}^{\mathrm{X}}\Big(\omega_{a}^{
\mathrm{JS}}\Big) - g_{b}^{\mathrm{X}}\Big(\omega_{b}^{\mathrm{JS}}
\Big) = 0 \bigg\}.
\label{eq:postINDEX:SET}
\end{equation}
\label{def:postINDEX:SET}
\end{definition}

For any existing $(2r - 1)$th-order mapped WENO schemes, e.g., the 
mapped WENO-X schemes in Table \ref{tab:map:WENO-X}, we have the 
following property.

\begin{lemma}
At $x_{j}$, for $\forall a, b = 0, 1, \cdots, r - 1$ and $a \neq b$, 
if $\textbf{postINDEX}(a, b, \mathrm{X}) \in \mathbb{S}^{\mathrm{X}}$
, then $x_{j}$ is an \textit{OP} point to the WENO-X scheme. 
Otherwise, if $\exists a,b = 0,1,\cdots,r-1$ and $a \neq b$, s.t. 
$\textbf{postINDEX}(a, b, \mathrm{X}) \notin \mathbb{S}^{\mathrm{X}}$
, then $x_{j}$ is a \textit{non-OP} point to the WENO-X scheme.
\label{lem:postINDEX}
\end{lemma}
\textbf{Proof.} 
According to Definitions \ref{def:OP-point}, \ref{def:postINDEX} and 
\ref{def:postINDEX:SET}, we can trivially finish the proof. 
$\hfill\square$ \\

By using the \textbf{postINDEX} function, we build a general method 
to introduce \textit{LOP} mappings into the existing \textit{non-OP} 
mapped WENO schemes, as given in Algorithm \ref{alg:posteriori}.

\begin{algorithm}[htb]
\caption{A general method to construct \textit{LOP} mappings.}
\label{alg:posteriori}
\SetKwInOut{Input}{input}\SetKwInOut{Output}{output}
\Input{$s$, index indicating the substencil $S_{s}$ and $s = 0,1,
\cdots, r-1$ \\
$d_{s}$, optimal weights \\
$\alpha^{\mathrm{JS}}_{s}$, nonnormalized nonlinear weights computed 
by the WENO-JS scheme \\
$\omega^{\mathrm{JS}}_{s}$, nonlinear weights computed by the 
WENO-JS scheme \\
$\big( g^{\mathrm{X}} \big)_{s}(\omega)$, nonnormalized nonlinear 
weights computed by the mapped WENO-X scheme}
\Output{$\Big\{\big(g^{\mathrm{LOP-X}}\big)_{s}(\omega^{
\mathrm{JS}}_{s}), s=0,1,\cdots,r-1\Big\}$, the new set of mapping 
functions that is \textit{LOP}}
\BlankLine
\emph{$\big( g^{\mathrm{X}} \big)_{s}(\omega), s = 0, 1,\cdots,r-1$ 
is a monotonically increasing mapping function over $[0, 1]$, and 
the set of mapping functions $\Big\{\big(g^{\mathrm{X}} \big)_{s}(
\omega), s = 0, 1,\cdots,r-1 \Big\}$ is \textit{non-OP}}\;
\tcp{implementation of the ``postINDEX'' function in Definition 
\ref{def:postINDEX}}
\For{$s_{1}=0; s_{1} \leq r - 2; s_{1} ++$}{
	\For{$s_{2} = s_{1} + 1; s_{2} \leq r - 1; s_{2} ++$}{
	   $\kappa=\mathrm{\textbf{postINDEX}}(s_{1},s_{2},\mathrm{X})$\;
	   \eIf{$\kappa \in \mathbb{S}^{\mathrm{X}}$}{
		   $\lambda = 1$\;
	   }{
		   $\lambda = 0$\;
		   \textbf{Break}\;
	   }	
	}
	\If{$\lambda = 0$}{
		\textbf{Break}\;
	}
}
\tcp{get $\big(g^{\mathrm{LOP-X}}\big)_{s}(\omega^{\mathrm{JS}}_{s})$}
\For{$s=0; s \leq r - 1; s ++$}{
	\eIf{$\lambda = 1$}{
		$\big(g^{\mathrm{LOP-X}}\big)_{s}(\omega^{\mathrm{JS}}_{s})
		= \big(g^{\mathrm{X}}\big)_{s}(\omega^{\mathrm{JS}}_{s})$\;
	}{
		$\big(g^{\mathrm{LOP-X}}\big)_{s}(\omega^{\mathrm{JS}}_{s})
		= \alpha_{s}^{\mathrm{JS}}$. 
		\tcp{$\alpha_{s}^{\mathrm{JS}}$ is computed by 
		Eq.\eqref{eq:weights:WENO-JS}}
	}
}
\end{algorithm}

\begin{theorem}
The set of mapping functions 
$\Big\{\big(g^{\mathrm{LOP-X}}\big)_{s}(\omega^{\mathrm{JS}}_{s}), 
s=0,1,\cdots,r-1\Big\}$ obtained through Algorithm 
\ref{alg:posteriori} is \textit{LOP}.
\label{theorem:g_MOP}
\end{theorem}
\textbf{Proof.} Naturally, the WENO-JS scheme could be treated as a 
mapped WENO scheme whose mapping function is defined as 
$\big( g^{\mathrm{JS}} \big)_{s}(\omega) = \omega,s= 0,\cdots,r-1$, 
and it is easy to verify that the set of mapping functions $\Big\{ 
\big( g^{\mathrm{JS}} \big)_{s}(\omega), s=0,\cdots,r-1\Big\}$ is 
\textit{LOP} while the widths of its optimal weight intervals 
(standing for the intervals about $\omega = d_{s}$ over which the 
mapping process attempts to use the optimal weights, see 
\cite{WENO-MAIMi,WENO-ACM}) are zero. Thus, for the case of 
$\lambda = 0$ in Algorithm \ref{alg:posteriori} (see line 22), 
the set of $\Big\{\big(g^{\mathrm{LOP-X}}\big)_{s}(\omega^{\mathrm{
JS}}_{s}), s = 0,1,\cdots, r-1 \Big\}$ is \textit{LOP} as 
$\alpha^{\mathrm{JS}}_{s}$ is the unnormalized weights associated 
with $\omega^{\mathrm{JS}}_{s}$. For the other case of $\lambda = 1$ 
in Algorithm \ref{alg:posteriori}, according to Lemma 
\ref{lem:postINDEX}, we can directly get that the set of
$\Big\{\big(g^{\mathrm{LOP-X}}\big)_{s}(\omega^{\mathrm{JS}}_{s}), 
s = 0,1,\cdots, r-1 \Big\}$ is \textit{LOP}. Now, we have finished 
the proof. 
$\hfill\square$ \\

We now define the modified weights which are \textit{LOP} as follows
\begin{equation}
\omega_{s}^{\mathrm{LOP-X}}=\dfrac{\alpha_{s}^{\mathrm{LOP-X}}}{
\sum_{l=0}^{r-1}\alpha_{l}^{\mathrm{LOP-X}}}, \quad \alpha_{s}^{
\mathrm{LOP-X}} = \big(g^{\mathrm{LOP-X}}\big)_{s}(\omega^{
\mathrm{JS}}_{s}),
\label{eq:LOP-mapping}
\end{equation}
where $s = 0,\cdots,r-1$ and $\big(g^{\mathrm{LOP-X}}\big)_{s}(\omega^{\mathrm{JS}}_{s})$ is computed through Algorithm \ref{alg:posteriori}. We denote by LOP-WENO-X the associated schemes.

\subsection{Convergence properties}
\begin{figure}[ht]
\centering
  \includegraphics[height=0.41\textwidth]
  {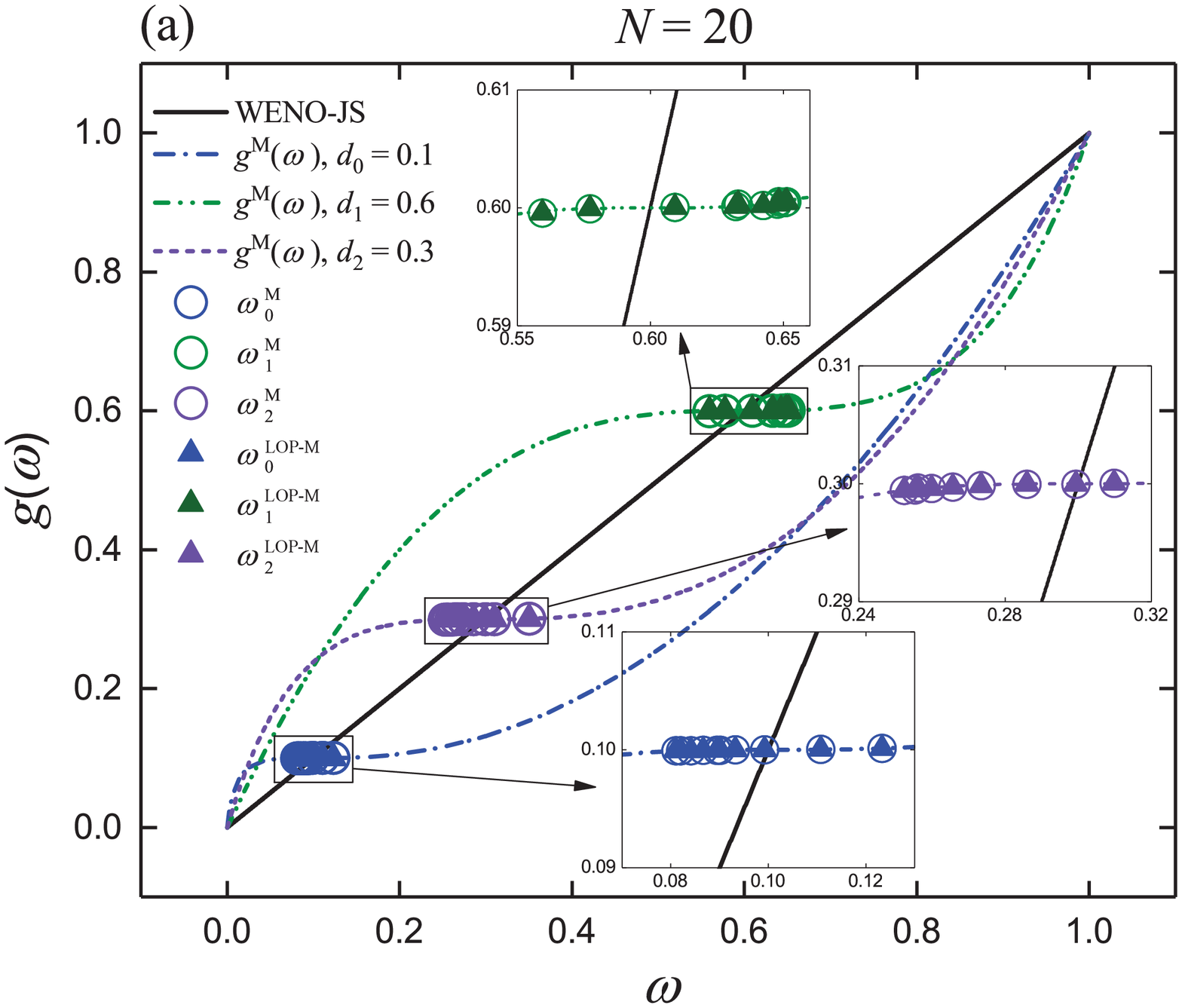}
  \includegraphics[height=0.41\textwidth]
  {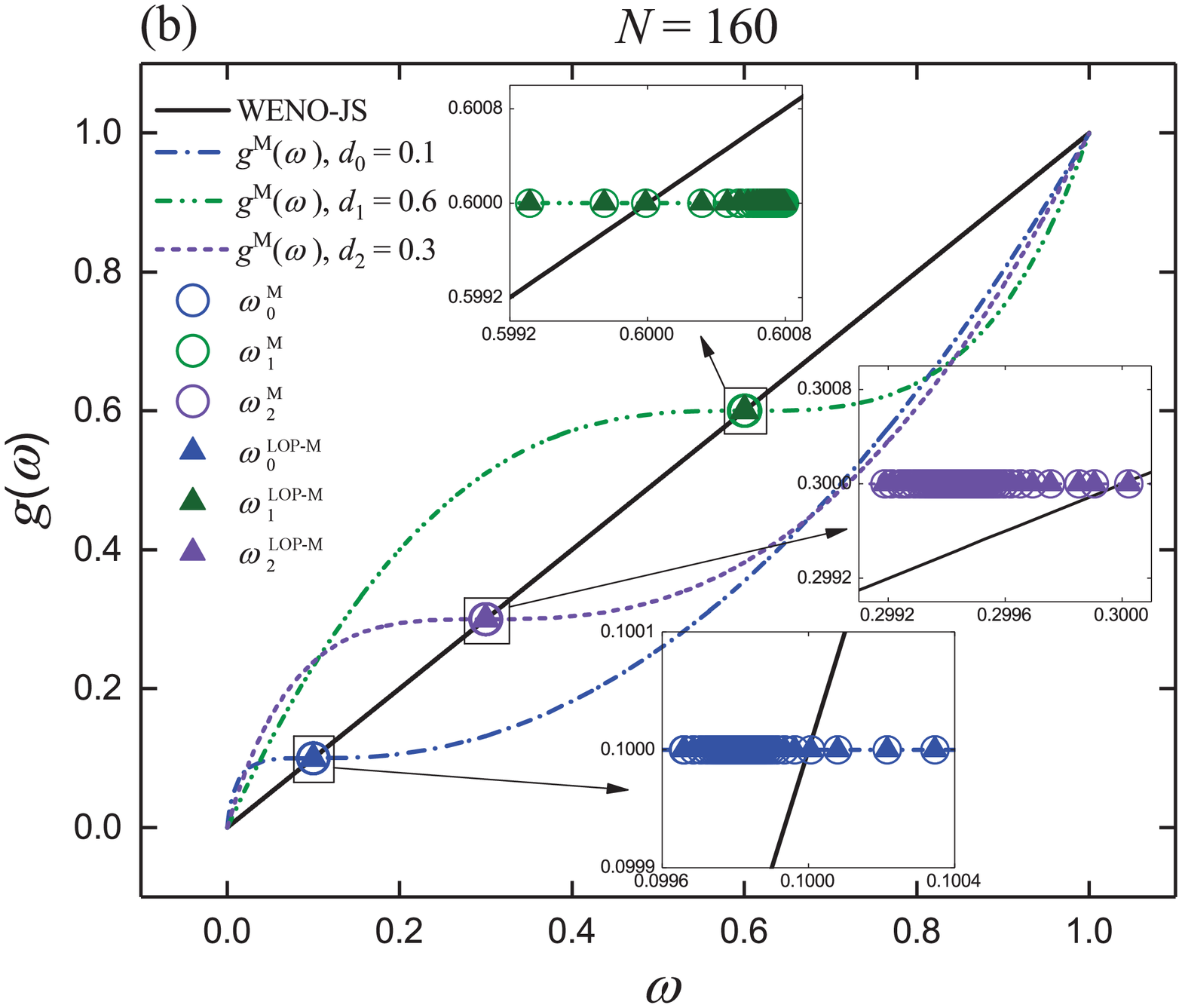} 
  \caption{Comparison for the real-time mappings on solving 
  Eq.\eqref{eq:LAE} with $u(x,0) = \sin(\pi x)$ using WENO-M and 
  LOP-WENO-M.}
\label{fig:realgOmega:LAE1:M}
\end{figure}

\begin{figure}[ht]
\centering
  \includegraphics[height=0.41\textwidth]
  {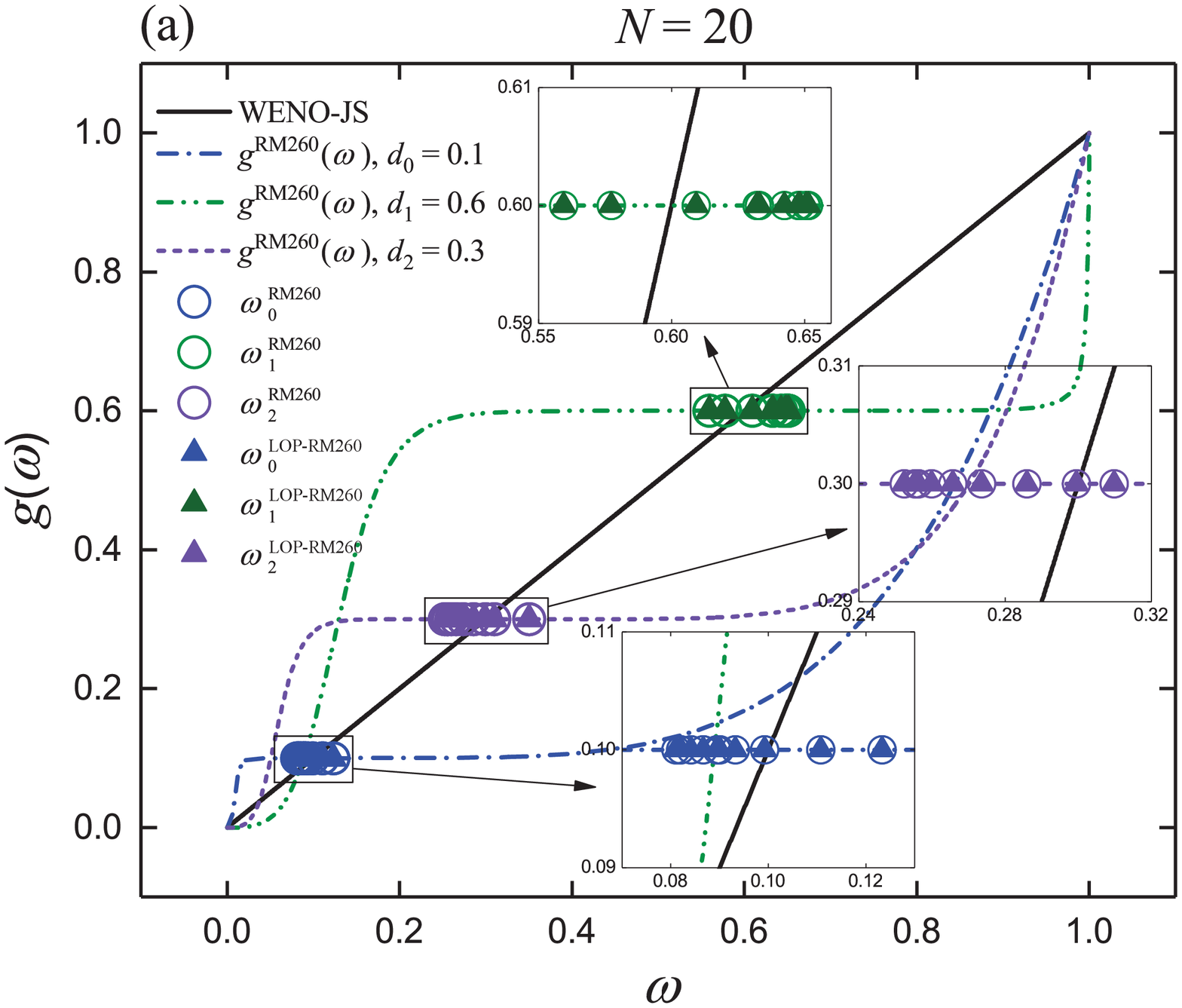}
  \includegraphics[height=0.41\textwidth]
  {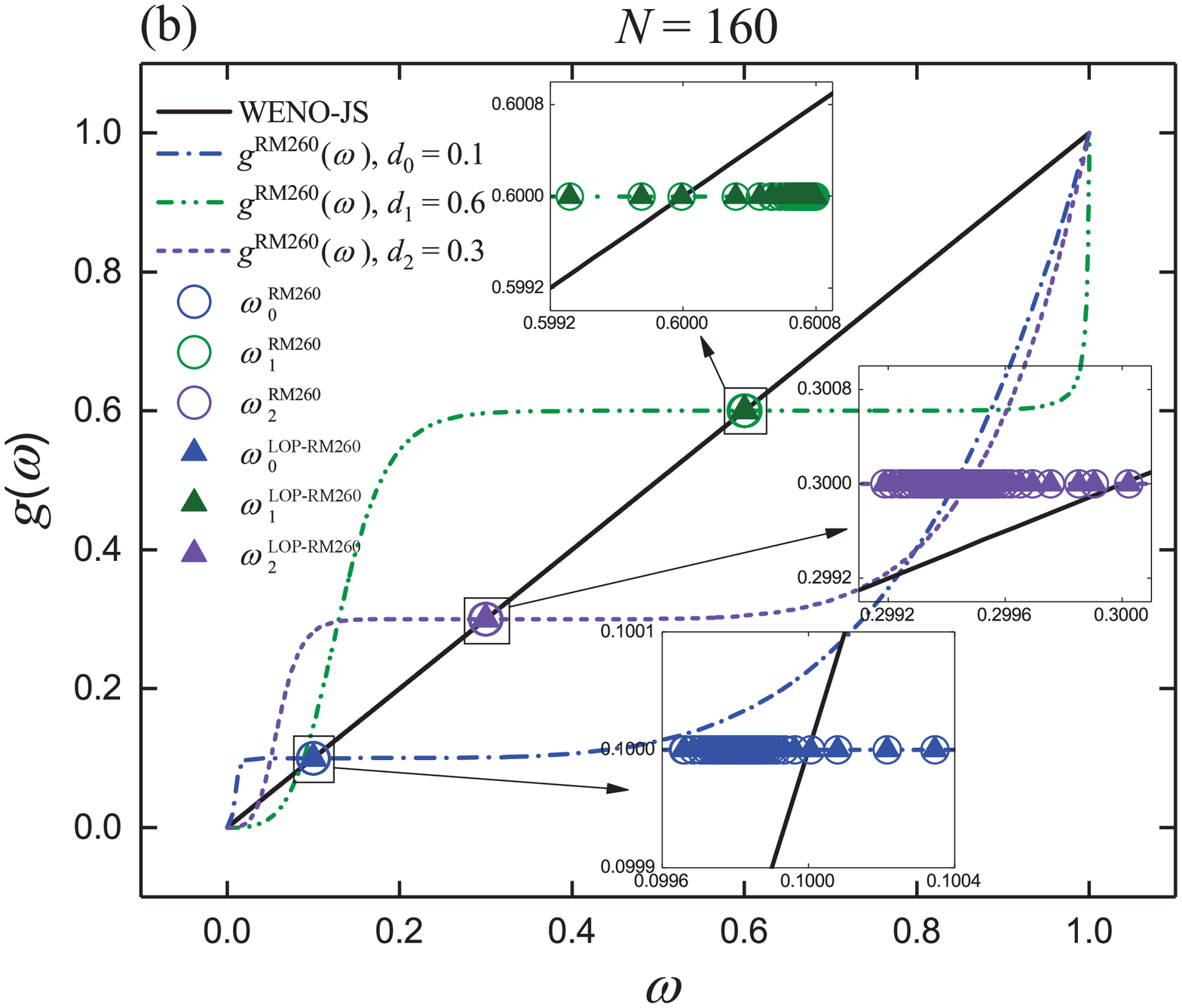} 
  \caption{Comparison for the real-time mappings on solving 
  Eq.\eqref{eq:LAE} with $u(x,0) = \sin(\pi x)$ using WENO-RM(260)
  and LOP-WENO-RM(260).}
\label{fig:realgOmega:LAE1:RM260}
\end{figure}

To study the properties of the mapping functions of the LOP-WENO-X 
schemes, we make a detailed analysis of the real-time mapping 
relationship. In constrast to commonly used mapping relationships 
that are directly computed by the designed mapping functions, the 
real-time mapping relationship here is obtained from the calculation 
of some specific problem at specified output time. Without loss of 
generality, we consider the following one-dimensional linear 
advection equation
\begin{equation}
\dfrac{\partial u}{\partial t} + \dfrac{\partial u}{\partial x} = 0,
\quad -1 \leq x \leq 1,
\label{eq:LAE}
\end{equation}
with the initial condition of $u(x,0) = \sin(\pi x)$. In 
Fig. \ref{fig:realgOmega:LAE1:M} and Fig. 
\ref{fig:realgOmega:LAE1:RM260}, we plot the real-time mapping 
relationships of the LOP-WENO-M and WENO-M schemes, the 
LOP-WENO-RM(260) and WENO-RM(260) schemes, as well as the 
designed mappings of the WENO-M and WENO-RM(260) schemes with 
$t = 2.0$. We find that the real-time mapping relationships of the 
LOP-WENO-M and LOP-WENO-RM(260) schemes are identical to those 
of the WENO-M and WENO-RM(260) schemes respectively. Actually, after 
extensive tests, the same results are observed for all other 
considered LOP-WENO-X and WENO-X schemes, and we do not present 
them here just for simplicity. Thus, we summarize this property as 
follows.
\begin{lemma}
The real-time mapping relationship of the LOP-WENO-X scheme is 
identical to that of the corresponding existing mapped WENO-X scheme 
presented in Table \ref{tab:map:WENO-X} in smooth regions where 
$n_{\mathrm{cp}} < r - 1$.
\label{lemma:real-timeMapping}
\end{lemma}

Then, we can trivially get the following Theorem.
\begin{theorem}
If $n_{\mathrm{cp}} < r - 1$, the mapping function 
$\big(g^{\mathrm{LOP-X}}\big)_{s}(\omega)$ obtained from Algorithm 
\ref{alg:posteriori} satisfies the 
following properties:\\

C1. for $\forall \omega \in (0,1), 
\big(g^{\mathrm{LOP-X}}\big)'_{s}(\omega) \geq 0$;

C2. for $\forall \omega \in \Omega$, $0\leq \big(g^{\mathrm{LOP-X}}
\big)_{s}(\omega) \leq 1$;

C3. $\big(g^{\mathrm{LOP-X}}\big)_{s}(d_{s}) = d_{s}, 
\big(g^{\mathrm{LOP-X}}\big)'_{s}(d_{s}) = \cdots = 
\big(g^{\mathrm{LOP-X}}\big)^{(n_{\mathrm{X}})}_{s}
(d_{s}) = 0$ where $n_{\mathrm{X}}$ is given in Table 
\ref{tab:map:WENO-X};

C4. $\big(g^{\mathrm{LOP-X}}\big)_{s}(0) = 0, 
\big(g^{\mathrm{LOP-X}}\big)_{s}(1) = 1, 
\big(g^{\mathrm{LOP-X}} \big)'_{s}(0) = 
\big(g^{\mathrm{X}} \big)'_{s}(0), 
\big(g^{\mathrm{LOP-X}} \big)'_{s}(1) = 
\big(g^{\mathrm{X}} \big)'_{s}(1)$.
\label{theorem:LOP-mapping}
\end{theorem}

\begin{figure}[!ht]
\centering
  \includegraphics[height=0.41\textwidth]
  {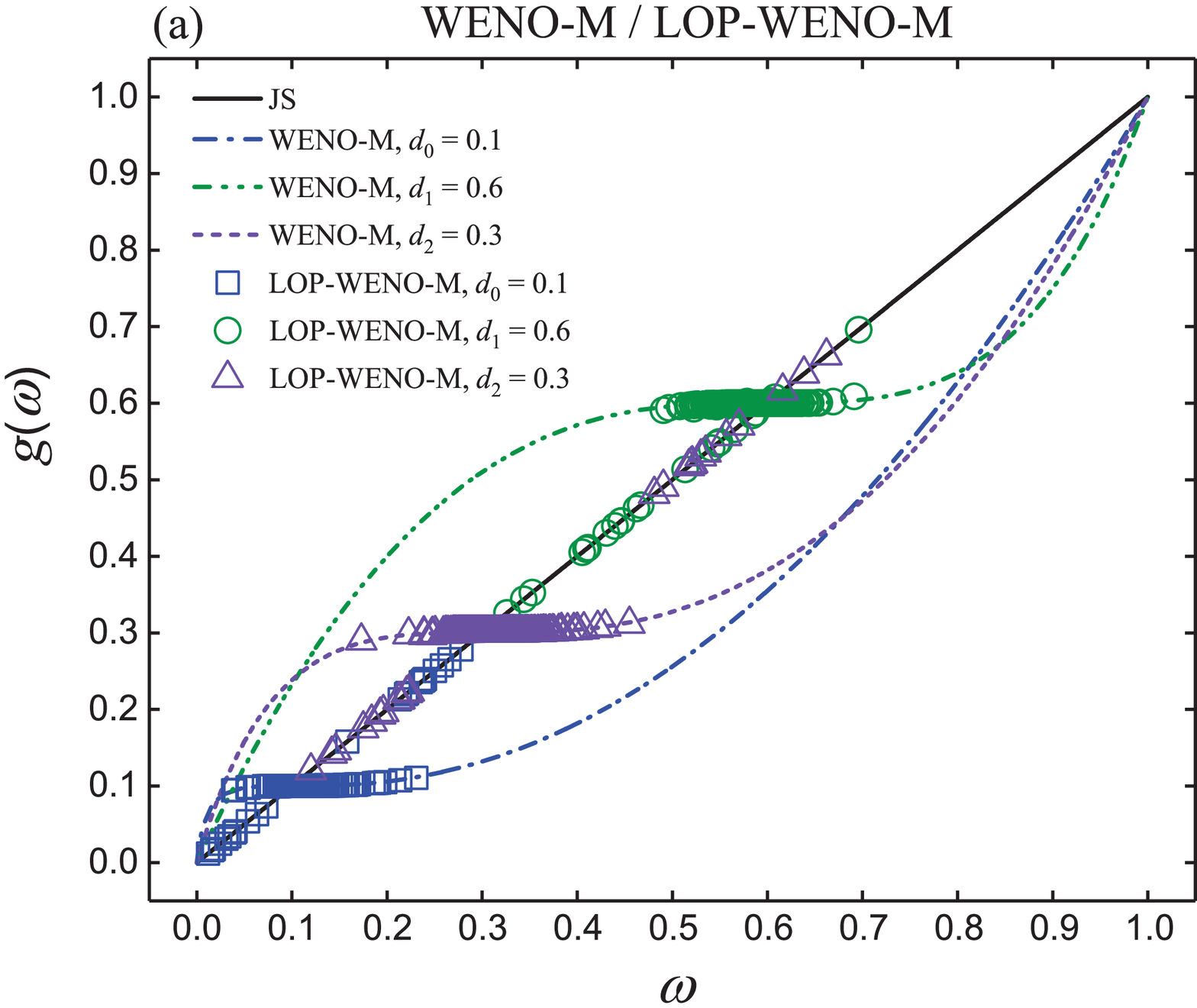}
  \includegraphics[height=0.41\textwidth]
  {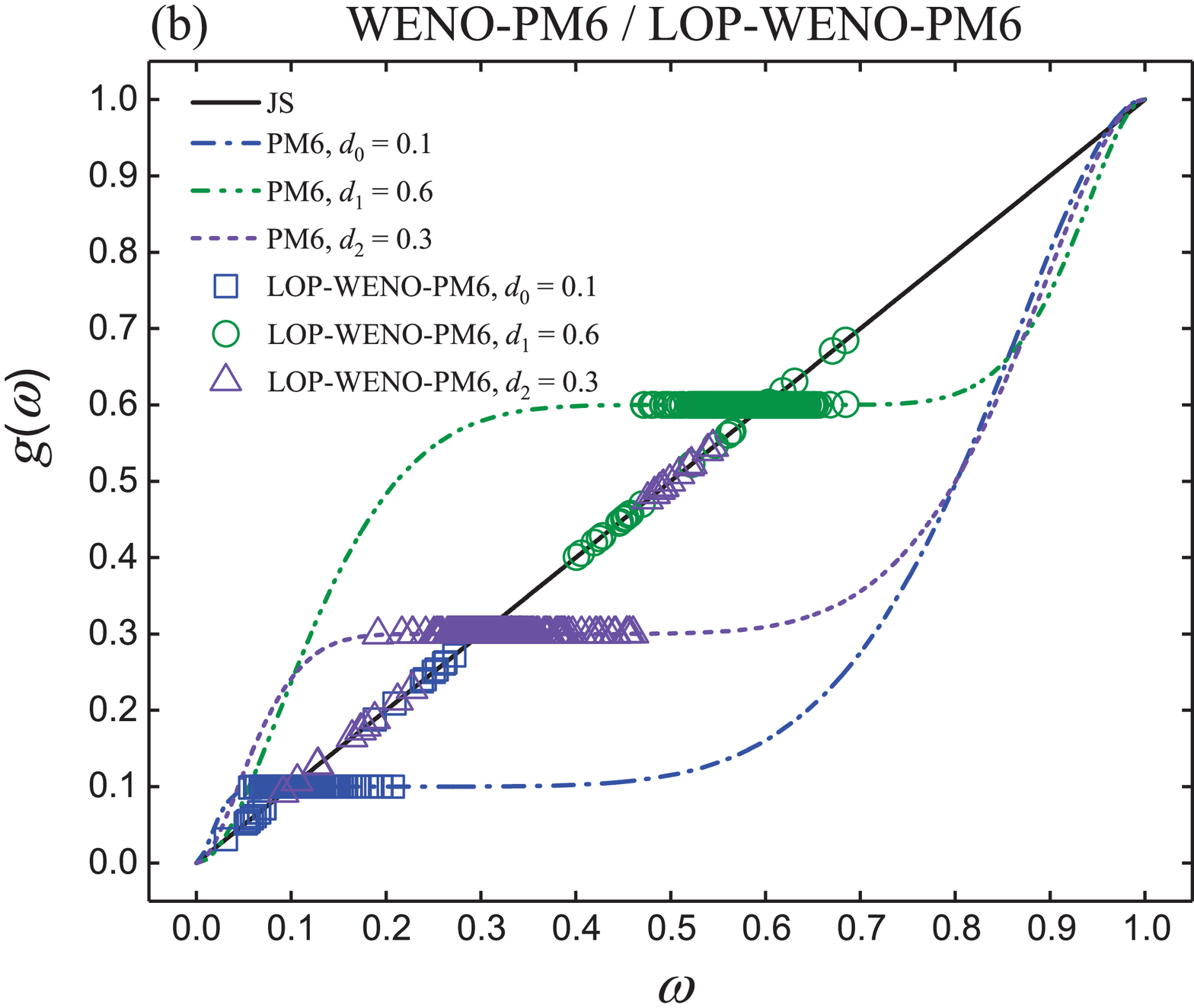} \\
  \includegraphics[height=0.41\textwidth]
  {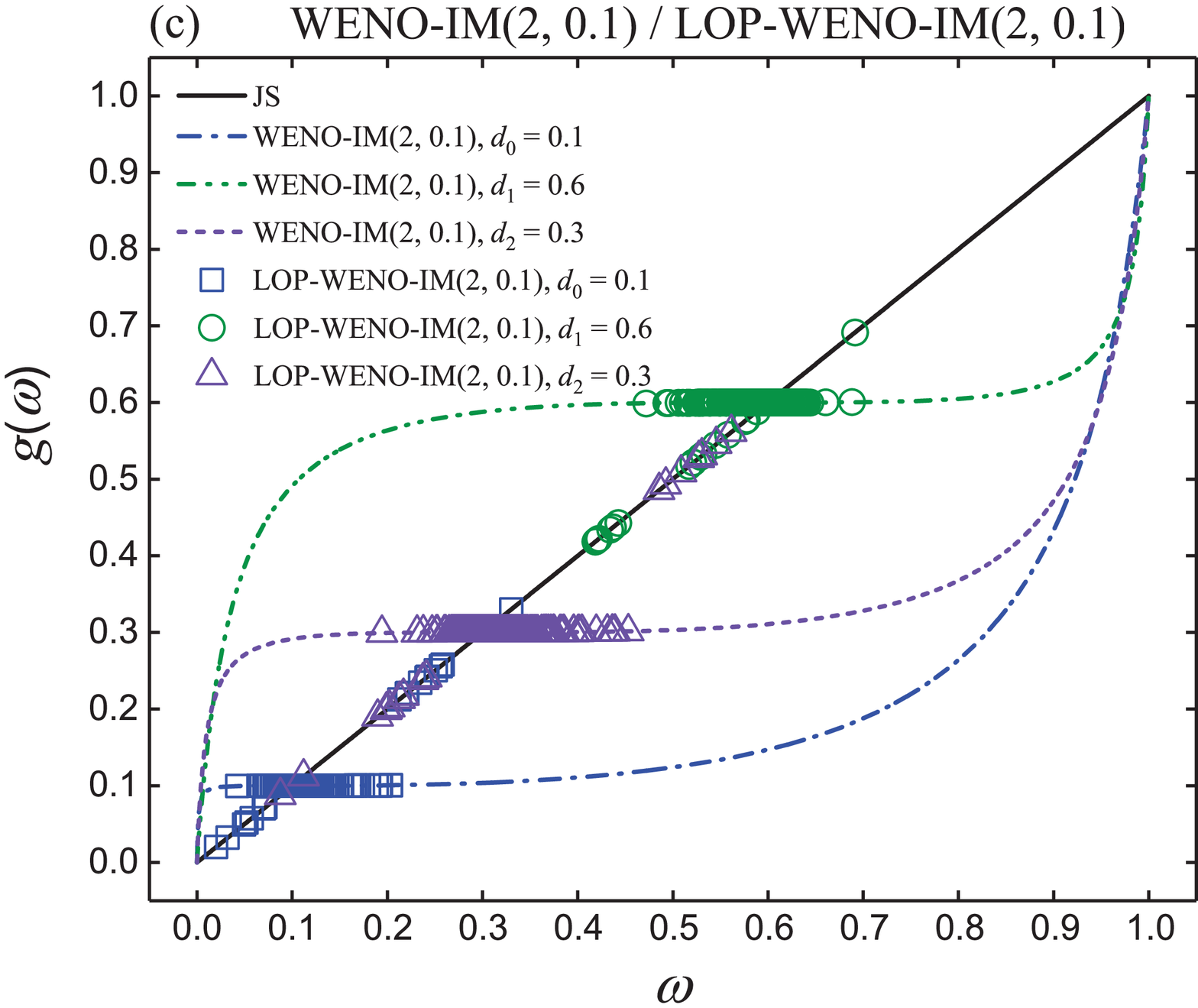}  
  \includegraphics[height=0.41\textwidth]
  {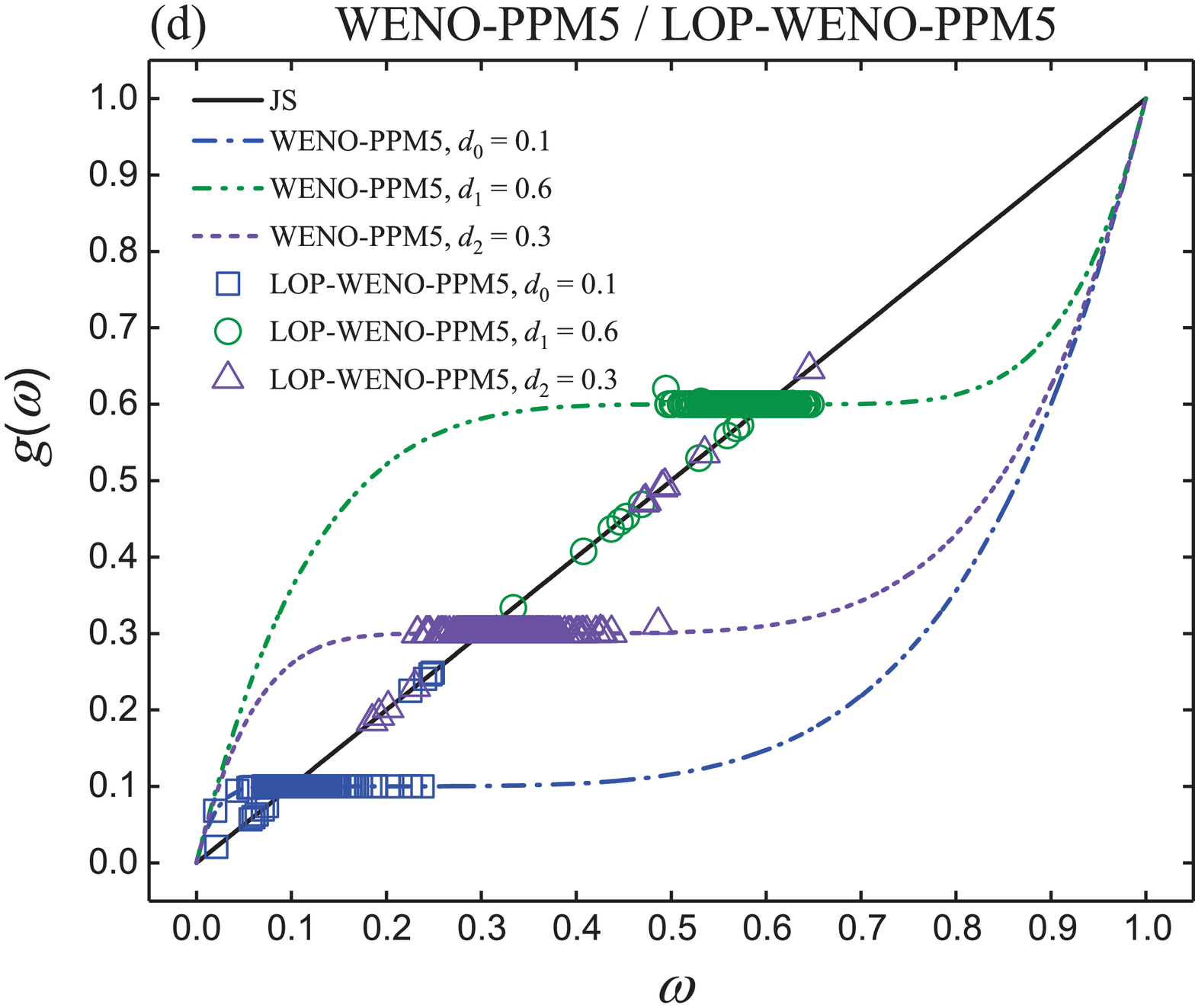}\\
  \includegraphics[height=0.41\textwidth]
  {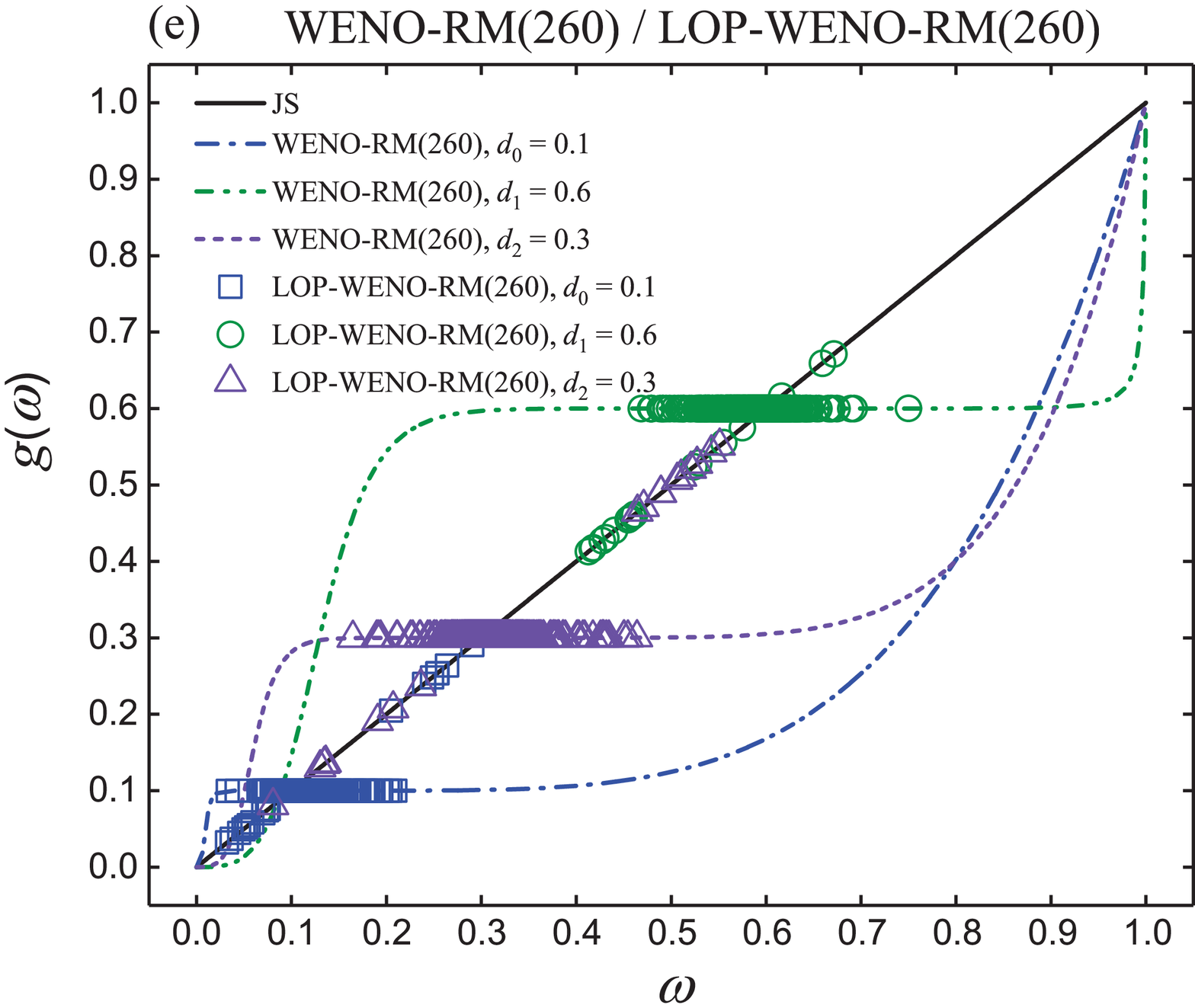} 
  \includegraphics[height=0.41\textwidth]
  {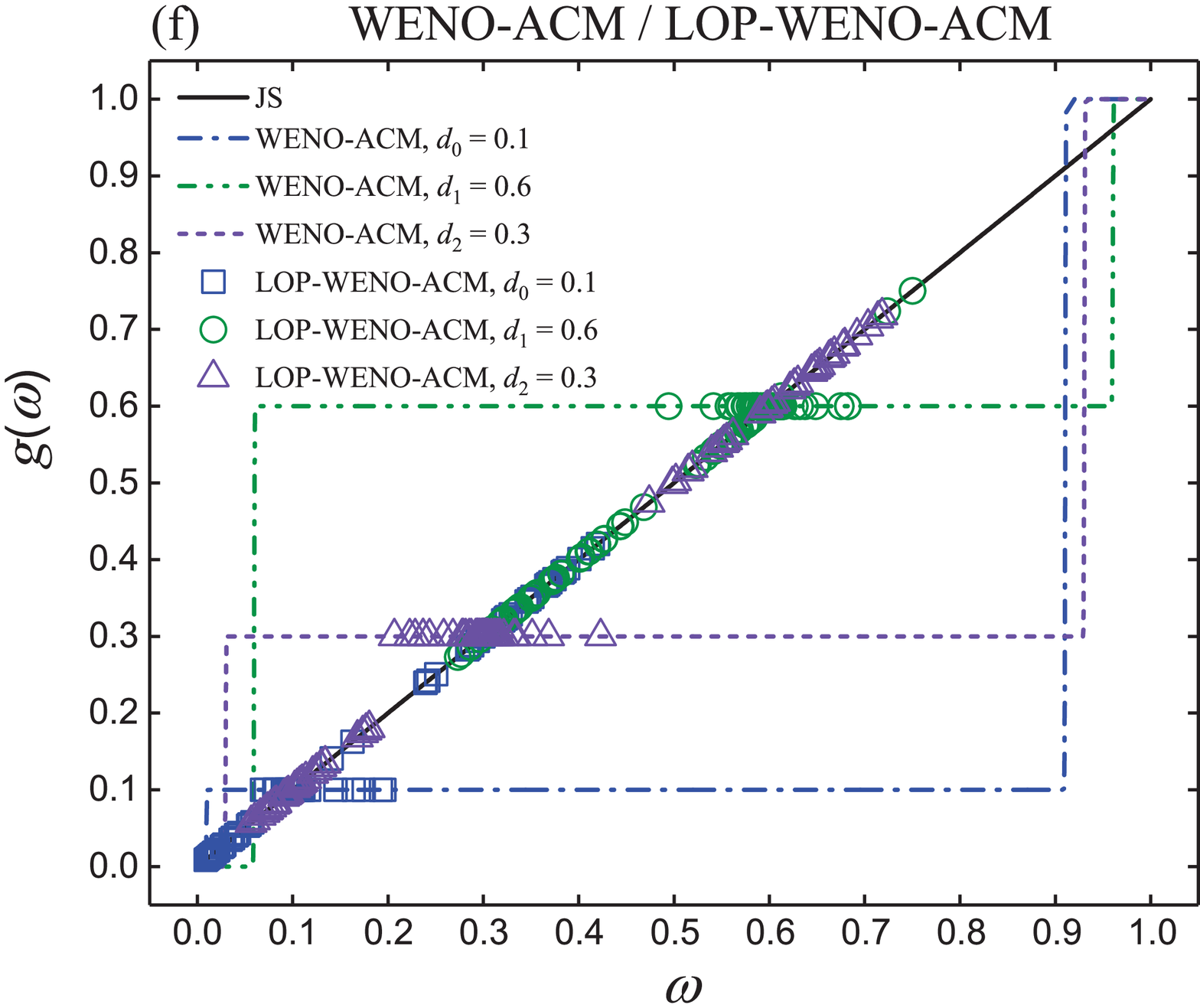}    
  \caption{Comparison for the real-time mappings on solving SLP
  using WENO-X and LOP-WENO-X.}
\label{fig:realgOmega:SLP}
\end{figure}

It is worthy to indicate that Lemma \ref{lemma:real-timeMapping} is 
not always true when the initial condition includes 
discontinuities. To show this, we solve Eq.\eqref{eq:LAE} with an 
initial condition as follows
\begin{equation}
\begin{array}{l}
u(x, 0) = \left\{
\begin{array}{ll}
\dfrac{1}{6}\big[ G(x, \beta, z - \hat{\delta}) + 4G(x, \beta, z) + G
(x, \beta, z + \hat{\delta}) \big], & x \in [-0.8, -0.6], \\
1, & x \in [-0.4, -0.2], \\
1 - \big\lvert 10(x - 0.1) \big\rvert, & x \in [0.0, 0.2], \\
\dfrac{1}{6}\big[ F(x, \alpha, a - \hat{\delta}) + 4F(x, \alpha, a) +
 F(x, \alpha, a + \hat{\delta}) \big], & x \in [0.4, 0.6], \\
0, & \mathrm{otherwise},
\end{array}\right. 
\end{array}
\label{eq:LAE:SLP}
\end{equation}
where $G(x, \beta, z) = \mathrm{e}^{-\beta (x - z)^{2}}, F(x, \alpha
, a) = \sqrt{\max \big(1 - \alpha ^{2}(x - a)^{2}, 0 \big)}$, and 
the constants are $z = -0.7, \hat{\delta} = 0.005, \beta = \dfrac{
\log 2}{36\hat{\delta} ^{2}}, a = 0.5$ and $\alpha = 10$. For 
brevity in the presentation, we call this \textit{\textbf{L}inear 
\textbf{P}roblem} SLP as it is presented by \textit{\textbf{S}hu} et 
al. in \cite{WENO-JS}. It is known that this problem consists of a 
Gaussian, a square wave, a sharp triangle and a semi-ellipse.  
In the calculations here, the periodic boundary condition is used 
and the CFL number is taken to be $0.1$. We take a uniform cell 
number of $N = 800$ and a short output time of $t = 2$. In Fig. 
\ref{fig:realgOmega:SLP}, we give the real-time mapping 
relationships of the LOP-WENO-X schemes, as well as the WENO-X 
schemes. It can be seen that the real-time mapping relationships 
of the LOP-WENO-X schemes are very different from those of the 
WENO-X schemes.

In Theorem \ref{theorem:LOP-convergence}, we give the convergence 
properties of the $(2r - 1)$th-order LOP-WENO-X schemes. As 
Theorem \ref{theorem:LOP-mapping} is true, the proof of Theorem 
\ref{theorem:LOP-convergence} is almost identical to that of the 
associated WENO-X schemes in the references presented in Table 
\ref{tab:map:WENO-X}.

\begin{theorem}
The requirements for the $(2r-1)$th-order LOP-WENO-X schemes to 
achieve the optimal order of accuracy are identical to those of their
associated $(2r-1)$th-order WENO-X scheme.
\label{theorem:LOP-convergence}
\end{theorem}

\subsection{Long-run simulations of linear advection equation for 
comparison}
\label{longRunLAE}
\subsubsection{With high-order critical points}
In order to demonstrate the ability of the LOP-WENO-X schemes that 
they can preserve high resolutions for the problem including 
high-order critical points with long output times, we conduct the 
following test.
\begin{example}
\rm{We solve Eq.\eqref{eq:LAE} with the periodic boundary condition 
by considering an initial condition that has high-order critical 
points, taking the form} 
\label{ex:long-run:case1}
\end{example} 
\begin{equation}
u(x, 0) = \exp\Big(-\big(x - 9.0\big)^{5}\cos^{9}\big(\pi(x - 9.0)
\big) \Big), 
\label{eq:long-run:case1}
\end{equation}
In this test, the computational domain is $x \in (7.5, 10.5)$ and 
the CFL number is $(\Delta x)^{2/3}$. 

The following $L_{1}$ and $L_{\infty}$ errors are computed to test 
the dissipations of the schemes
\begin{equation}
\displaystyle
\begin{aligned}
L_{1} = h \cdot \displaystyle\sum\nolimits_{i=1}^{N}\Big\lvert 
u_{i}^{\mathrm{exact}} - (u_{h})_{i} \Big\rvert, \quad 
L_{\infty} = \displaystyle\max_{1\leq i\leq N} \Big\lvert 
u_{i}^{\mathrm{exact}} -(u_{h})_{i}
\Big\rvert,
\end{aligned}
\label{normsDefinition}
\end{equation}
where $N$ is the number of the cells and $h$ is the associated 
uniform spatial step size. $(u_{h})_{i}$ is the numerical solution 
and $u_{i}^{\mathrm{exact}}$ is the exact solution. It is trivial to 
verify that the exact solution is $u(x,t) = \exp\Big(-\big((x - t) - 
9.0\big)^{5}\cos^{9}\big(\pi((x - t) - 9.0)\big) \Big)$. 

In addition, the following increased errors are considered
\begin{equation*}
\begin{aligned}
\chi_{1}=\frac{L_{1}^{\mathrm{Y}}(t)-L_{1}^{\mathrm{ILW}}(t)}{L_{1}^{\mathrm{ILW}}(t)}\times100\%, \quad 
\chi_{\infty}=\frac{L_{\infty}^{\mathrm{Y}}(t)-L_{\infty}^{\mathrm{ILW}}(t)}{L_{\infty}^{\mathrm{ILW}}(t)}\times100\%,
\end{aligned}
\end{equation*}
where $L_{1}^{\mathrm{ILW}}(t)$ and $L_{\infty}^{\mathrm{ILW}}(t)$ 
are the $L_{1}$ and $L_{\infty}$ errors of the WENO5-ILW scheme 
(say, the WENO5 scheme using ideal linear weights), and similarly, 
$L_{1}^{\mathrm{Y}}(t)$ and $L_{\infty}^{\mathrm{Y}}(t)$ are those 
of the scheme ``Y''. 

\begin{table}[!ht]
\begin{myFontSize}
\centering
\caption{$L_{1}$, $L_{\infty}$ errors and the increased errors (in 
percentage) of various considered schemes on solving $u_{t}+u_{x}=0$
with initial condition Eq.\eqref{eq:long-run:case1} and $N = 300$.}
\label{tab::long-run:N300}
\begin{tabular*}{\hsize}
{@{}@{\extracolsep{\fill}}rlrlrlrlr@{}}
\hline
\space    &\multicolumn{4}{l}{\cellcolor{gray!35}{WENO5-ILW}}  
          &\multicolumn{4}{l}{\cellcolor{gray!35}{WENO-JS}}  \\
\cline{2-5}  \cline{6-9}
Time, $t$             & $L_{1}$ error      & $\chi_{1}$
                      & $L_{\infty}$ error & $\chi_{\infty}$
		  			  & $L_{1}$ error      & $\chi_{1}$
      				  & $L_{\infty}$ error & $\chi_{\infty}$ \\
\Xhline{0.65pt}
15                    & 9.39243E-04        & -
                      & 1.43469E-03        & -
                      & 1.53437E-03        & 63\%
                      & 2.70581E-03        & 89\% \\
60                    & 1.18694E-03        & -
                      & 2.20682E-03        & -
                      & 4.94939E-03        & 317\%
                      & 8.38132E-03        & 280\% \\
150                   & 2.85184E-03        & -
                      & 4.96667E-03        & -
                      & 2.15858E-02        & 657\%
                      & 5.25441E-02        & 958\% \\ 
300                   & 5.39974E-03        & -
                      & 8.81363E-03        & -
                      & 7.93589E-02        & 1370\%
                      & 1.34321E-01        & 1424\% \\
600                   & 9.94133E-03        & -
                      & 1.50917E-02        & -
                      & 2.10016E-01        & 2013\%
                      & 3.04860E-01        & 1920\% \\
900                   & 1.38061E-02        & -
                      & 1.96281E-02        & -
                      & 2.84632E-01        & 1962\%
                      & 4.20080E-01        & 2040\% \\
1200                  & 1.74067E-02        & -
                      & 2.39652E-02        & -
                      & 3.26687E-01        & 1777\%
                      & 5.14072E-01        & 2045\% \\
\hline
\space    &\multicolumn{4}{l}{\cellcolor{gray!35}{WENO-M}}
          &\multicolumn{4}{l}{\cellcolor{gray!35}{LOP-WENO-M}}\\
\cline{2-5}  \cline{6-9}
Time, $t$             & $L_{1}$ error      & $\chi_{1}$         
                      & $L_{\infty}$ error & $\chi_{\infty}$     
		  			  & $L_{1}$ error      & $\chi_{1}$         
      				  & $L_{\infty}$ error & $\chi_{\infty}$      \\
\Xhline{0.65pt}
15                    & 9.63640E-04        & 3\%
                      & 1.43486E-03        & 0\%
                      & 1.12245E-03        & 20\%
                      & 1.45965E-03        & 2\% \\
60                    & 1.59176E-03        & 34\%
                      & 4.93572E-03        & 124\%
                      & 3.13037E-03        & 164\%
                      & 9.68122E-03        & 339\% \\
150                   & 6.91167E-03        & 142\%
                      & 3.05991E-02        & 516\%
                      & 9.04418E-03        & 217\%
                      & 1.62967E-02        & 228\% \\ 
300                   & 3.02283E-02        & 460\%
                      & 1.13618E-01        & 1189\%
                      & 1.14248E-02        & 112\%
                      & 1.86573E-02        & 112\% \\
600                   & 9.06252E-02        & 812\%
                      & 1.76325E-01        & 1068\%
                      & 2.28157E-02        & 130\%
                      & 4.01795E-02        & 166\% \\
900                   & 1.52637E-01        & 1006\%
                      & 3.25300E-01        & 1557\%
                      & 2.84909E-02        & 106\%
                      & 3.99876E-02        & 104\% \\
1200                  & 1.95044E-01        & 1021\%
                      & 3.32396E-01        & 1287\%
                      & 2.68748E-02        & 54\%
                      & 3.08547E-02        & 29\% \\
\hline
\space    &\multicolumn{4}{l}{\cellcolor{gray!35}{WENO-PM6}} 
          &\multicolumn{4}{l}{\cellcolor{gray!35}{LOP-WENO-PM6}}\\
\cline{2-5}  \cline{6-9}
Time, $t$             & $L_{1}$ error      & $\chi_{1}$         
                      & $L_{\infty}$ error & $\chi_{\infty}$     
		  			  & $L_{1}$ error      & $\chi_{1}$         
      				  & $L_{\infty}$ error & $\chi_{\infty}$      \\
\Xhline{0.65pt}
15                    & 9.33637E-04        & -1\%
                      & 1.43467E-03        & 0\%
                      & 1.11346E-03        & 19\%
                      & 1.43572E-03        & 0\% \\
60                    & 1.13744E-03        & -4\%
                      & 2.20467E-03        & 0\%
                      & 3.03330E-03        & 156\%
                      & 9.06783E-03        & 311\% \\
150                   & 2.68760E-03        & -6\%
                      & 4.88656E-03        & -2\%
                      & 9.07249E-03        & 218\%
                      & 1.65412E-02        & 233\% \\ 
300                   & 4.97186E-03        & -8\%
                      & 8.79888E-03        & 0\%
                      & 1.09817E-02        & 103\%
                      & 1.96942E-02        & 123\% \\
600                   & 8.51389E-03        & -14\%
                      & 1.49420E-02        & -1\%
                      & 1.89507E-02        & 91\%
                      & 2.99502E-02        & 98\% \\
900                   & 1.12016E-02        & -19\%
                      & 1.90706E-02        & -3\%
                      & 2.55009E-02        & 85\%
                      & 4.82978E-02        & 146\% \\
1200                  & 1.36304E-02        & -22\%
                      & 2.30084E-02        & -4\%
                      & 2.47046E-02        & 42\%
                      & 3.57705E-02        & 49\% \\
\hline
\space    &\multicolumn{4}{l}{\cellcolor{gray!35}{WENO-IM(2, 0.1)}} 
    &\multicolumn{4}{l}{\cellcolor{gray!35}{LOP-WENO-IM(2, 0.1)}}\\
\cline{2-5}  \cline{6-9}
Time, $t$             & $L_{1}$ error      & $\chi_{1}$         
                      & $L_{\infty}$ error & $\chi_{\infty}$     
		  			  & $L_{1}$ error      & $\chi_{1}$         
      				  & $L_{\infty}$ error & $\chi_{\infty}$      \\
\Xhline{0.65pt}
15                    & 9.41045E-04        & 0\%
                      & 1.43471E-03        & 0\%
                      & 1.11781E-03        & 19\%
                      & 1.76330E-03        & 23\% \\
60                    & 1.19465E-03        & 1\%
                      & 2.21080E-03        & 0\%
                      & 3.07175E-03        & 159\%
                      & 9.21010E-03        & 317\% \\
150                   & 2.79176E-03        & -2\%
                      & 4.99982E-03        & 1\%
                      & 9.02303E-03        & 216\%
                      & 1.65281E-02        & 233\% \\ 
300                   & 5.09842E-03        & -6\%
                      & 8.83560E-03        & 0\%
                      & 1.09585E-02        & 103\%
                      & 1.88913E-02        & 114\% \\
600                   & 8.84945E-03        & -11\%
                      & 1.48704E-02        & -1\%
                      & 2.07426E-02        & 109\%
                      & 3.28878E-02        & 118\% \\
900                   & 1.17416E-02        & -15\%
                      & 1.91100E-02        & -3\%
                      & 2.64361E-02        & 91\%
                      & 4.61971E-02        & 135\% \\
1200                  & 1.43988E-02        & -17\%
                      & 2.31394E-02        & -3\%
                      & 2.54519E-02        & 46\%
                      & 2.93893E-02        & 23\% \\
\hline
\space    &\multicolumn{4}{l}{\cellcolor{gray!35}{WENO-PPM5}} 
          &\multicolumn{4}{l}{\cellcolor{gray!35}{LOP-WENO-PPM5}}\\
\cline{2-5}  \cline{6-9}
Time, $t$             & $L_{1}$ error      & $\chi_{1}$         
                      & $L_{\infty}$ error & $\chi_{\infty}$     
		  			  & $L_{1}$ error      & $\chi_{1}$         
      				  & $L_{\infty}$ error & $\chi_{\infty}$      \\
\Xhline{0.65pt}
15                    & 9.34095E-04        & -1\%
                      & 1.43467E-03        & 0\%
                      & 1.11267E-03        & 18\%
                      & 1.43465E-03        & 0\% \\
60                    & 1.14062E-03        & -4\%
                      & 2.20440E-03        & 0\%
                      & 3.02337E-03        & 155\%
                      & 9.03860E-03        & 310\% \\
150                   & 2.69185E-03        & -6\%
                      & 4.89513E-03        & -1\%
                      & 9.07221E-03        & 218\%
                      & 1.65292E-02        & 233\% \\ 
300                   & 4.98381E-03        & -8\%
                      & 8.79434E-03        & 0\%
                      & 1.10228E-02        & 104\%
                      & 2.01451E-02        & 129\% \\
600                   & 8.54488E-03        & -14\%
                      & 1.49710E-02        & -1\%
                      & 1.89501E-02        & 91\%
                      & 2.96365E-02        & 96\% \\
900                   & 1.12525E-02        & -18\%
                      & 1.91124E-02        & -3\%
                      & 2.55066E-02        & 85\%
                      & 4.90556E-02        & 150\% \\
1200                  & 1.37002E-02        & -21\%
                      & 2.30576E-02        & -4\%
                      & 2.48630E-02        & 43\%
                      & 3.66520E-02        & 53\% \\
\hline
\space    &\multicolumn{4}{l}{\cellcolor{gray!35}{WENO-RM(260)}} 
      &\multicolumn{4}{l}{\cellcolor{gray!35}{LOP-WENO-RM(260)}}\\
\cline{2-5}  \cline{6-9}
Time, $t$             & $L_{1}$ error      & $\chi_{1}$         
                      & $L_{\infty}$ error & $\chi_{\infty}$     
		  			  & $L_{1}$ error      & $\chi_{1}$         
      				  & $L_{\infty}$ error & $\chi_{\infty}$      \\
\Xhline{0.65pt}
15                    & 9.39107E-04        & 0\%
                      & 1.43469E-03        & 0\%
                      & 1.11731E-03        & 19\%
                      & 1.44384E-03        & 1\% \\
60                    & 1.16183E-03        & -2\%
                      & 2.20830E-03        & 0\%
                      & 3.06524E-03        & 158\%
                      & 9.17780E-03        & 316\% \\
150                   & 2.70494E-03        & -5\%
                      & 4.89842E-03        & -1\%
                      & 9.02171E-03        & 216\%
                      & 1.65533E-02        & 233\% \\ 
300                   & 4.99464E-03        & -8\%
                      & 8.80531E-03        & 0\%
                      & 1.09337E-02        & 102\%
                      & 1.88283E-02        & 114\% \\
600                   & 8.56131E-03        & -14\%
                      & 1.49164E-02        & -1\%
                      & 2.03795E-02        & 105\%
                      & 3.22925E-02        & 114\% \\
900                   & 1.12528E-02        & -18\%
                      & 1.90455E-02        & -3\%
                      & 2.62193E-02        & 90\%
                      & 4.60197E-02        & 134\% \\
1200                  & 1.36724E-02        & -21\%
                      & 2.29797E-02        & -4\%
                      & 2.53001E-02        & 45\%
                      & 2.94014E-02        & 23\% \\
\hline
\space    &\multicolumn{4}{l}{\cellcolor{gray!35}{WENO-ACM}} 
          &\multicolumn{4}{l}{\cellcolor{gray!35}{LOP-WENO-ACM}}\\
\cline{2-5}  \cline{6-9}
Time, $t$             & $L_{1}$ error      & $\chi_{1}$         
                      & $L_{\infty}$ error & $\chi_{\infty}$     
		  			  & $L_{1}$ error      & $\chi_{1}$         
      				  & $L_{\infty}$ error & $\chi_{\infty}$      \\
\Xhline{0.65pt}
15                    & 9.39110E-04        & 0\%
                      & 1.43469E-03        & 0\%
                      & 1.11724E-03        & 19\%
                      & 1.44365E-03        & 1\% \\
60                    & 1.17055E-03        & -1\%
                      & 2.20771E-03        & 0\%
                      & 3.06492E-03        & 158\%
                      & 9.17648E-03        & 316\% \\
150                   & 2.70928E-03        & -5\%
                      & 4.92156E-03        & -1\%
                      & 9.02312E-03        & 216\%
                      & 1.65537E-02        & 233\% \\ 
300                   & 5.02511E-03        & -7\%
                      & 8.79766E-03        & 0\%
                      & 1.09333E-02        & 102\%
                      & 1.88278E-02        & 114\% \\
600                   & 8.63079E-03        & -13\%
                      & 1.49585E-02        & -1\%
                      & 2.03612E-02        & 105\%
                      & 3.22898E-02        & 114\% \\
900                   & 1.13152E-02        & -18\%
                      & 1.91098E-02        & -3\%
                      & 2.62113E-02        & 90\%
                      & 4.59477E-02        & 134\% \\
1200                  & 1.37302E-02        & -21\%
                      & 2.30445E-02        & -4\%
                      & 2.53002E-02        & 45\%
                      & 2.94013E-02        & 23\% \\
\hline
\end{tabular*}
\end{myFontSize}
\end{table}

In Table \ref{tab::long-run:N300}, we present the $L_{1}$ and
$L_{\infty}$ errors and the increased errors with $N=300$ at various 
final times of $t = 15, 60, 150, 300, 600, 900, 1200$. We can find 
that: (1) the WENO-JS scheme produces the largest $L_{1}$ and
$L_{\infty}$ errors,  leading to the largest increased errors, among 
all considered schemes for each output time; (2) when the output 
time is small, like $t \leq 60$, the WENO-M scheme provides more 
accurate results than the LOP-WENO-M scheme, leading to smaller 
increased errors; (3) however, when the output time gets larger, 
like $t \geq 300$, the increased errors of the LOP-WENO-M 
scheme evidently decrease and get closer to those of the WENO5-ILW 
scheme, while the errors of the WENO-M scheme increase significantly 
leading to extremely larger increased errors; (4) in spite that the 
errors of the LOP-WENO-X schemes (except the case of ``X = M'') 
are larger than those of their associated WENO-X schemes, these 
errors can maintain a acceptable level leading to tolerable 
increased errors that are far lower than those of the WENO-JS and 
WENO-M schemes, and moreover, the LOP-WENO-X schemes are able to 
avoid the spurious oscillations on solving problems with 
discontinuities for long output times while their associated WENO-X 
schemes fail to prevent the spurious oscillations (for example, see 
Fig. \ref{fig:Z:N1600} and Fig. \ref{fig:Z:N3200}). 

\begin{figure}[!ht]
\centering
  \includegraphics[height=0.33\textwidth]
  {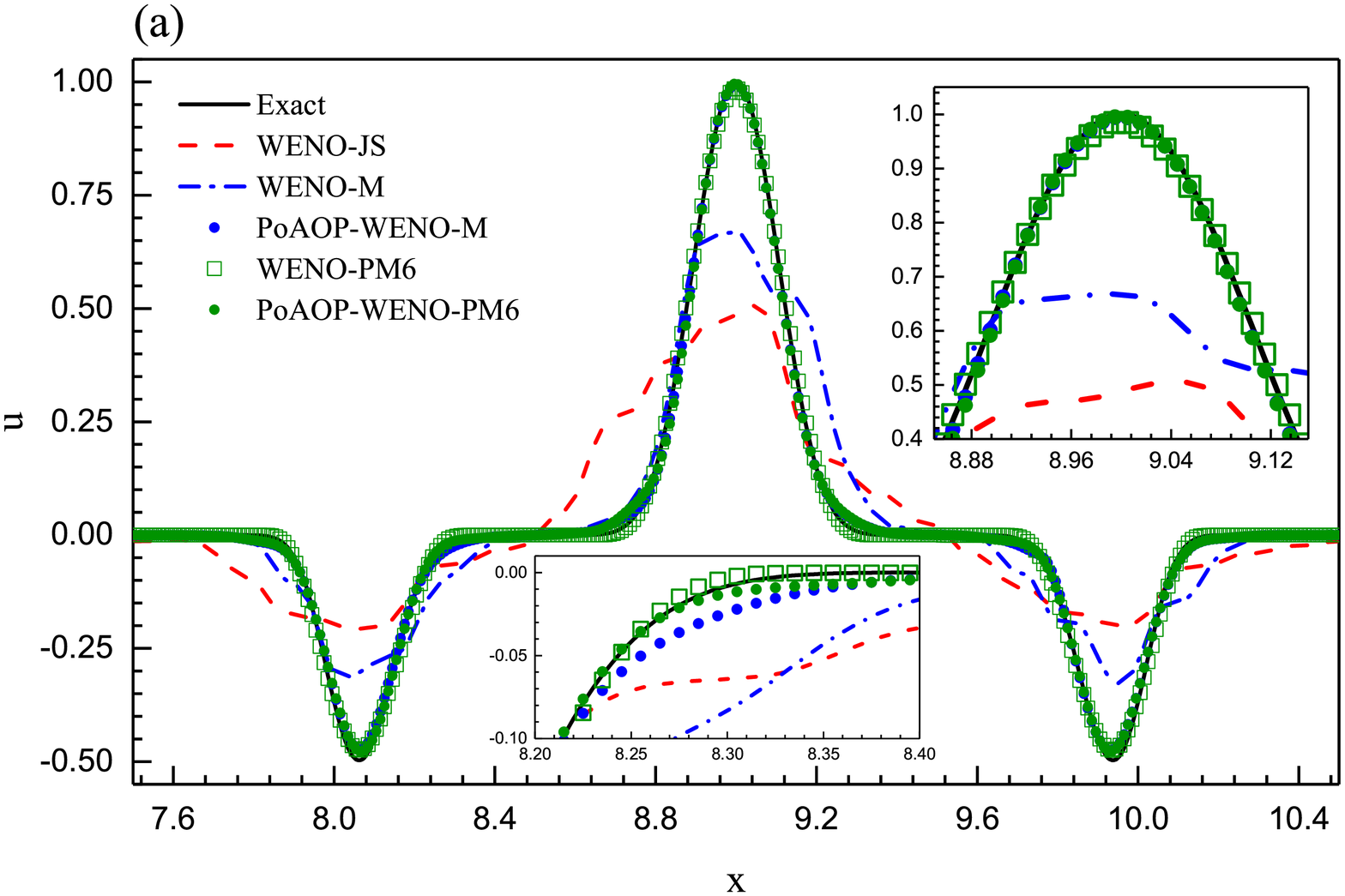}
  \includegraphics[height=0.33\textwidth]
  {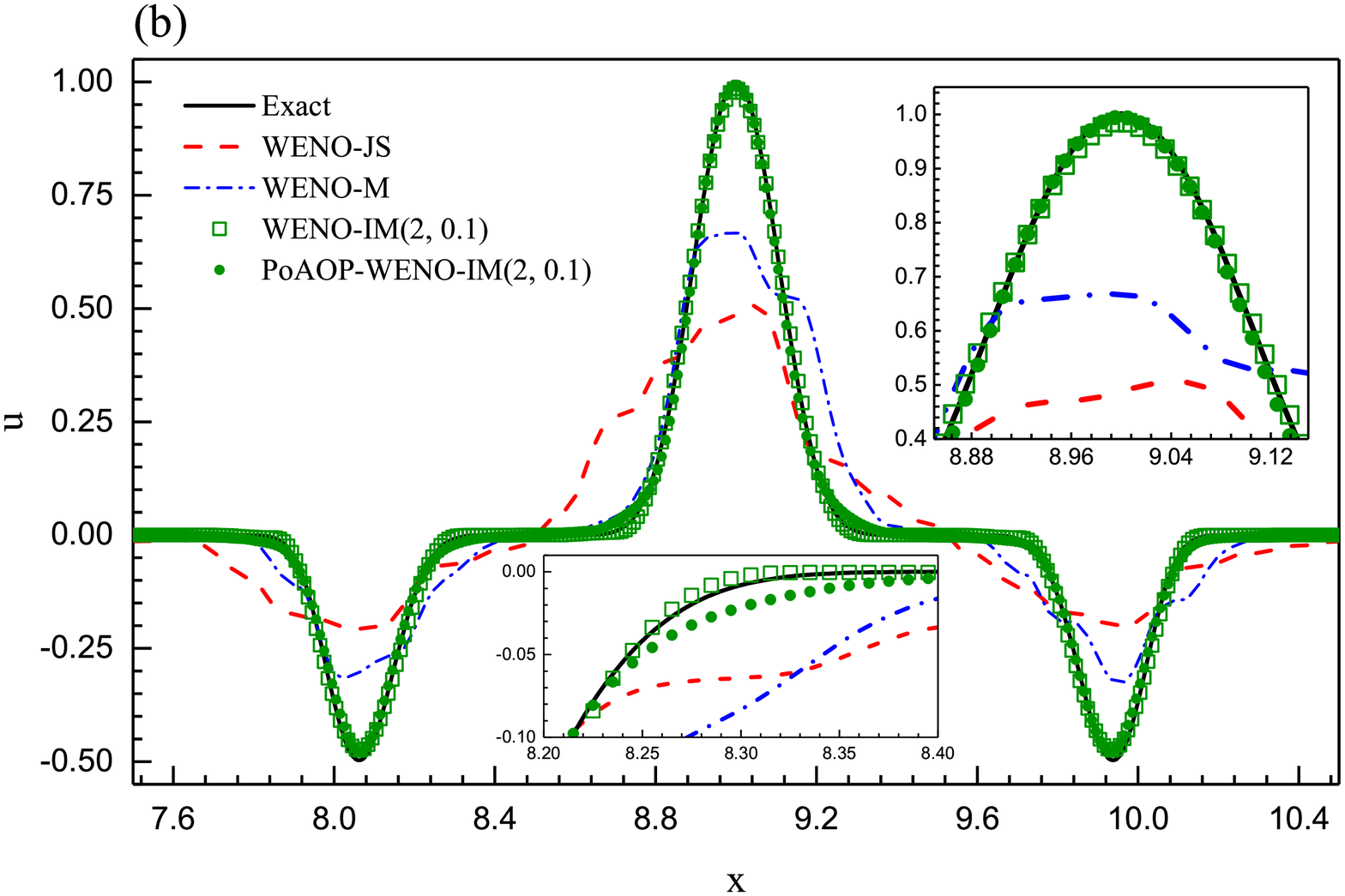}\\
  \includegraphics[height=0.33\textwidth]
  {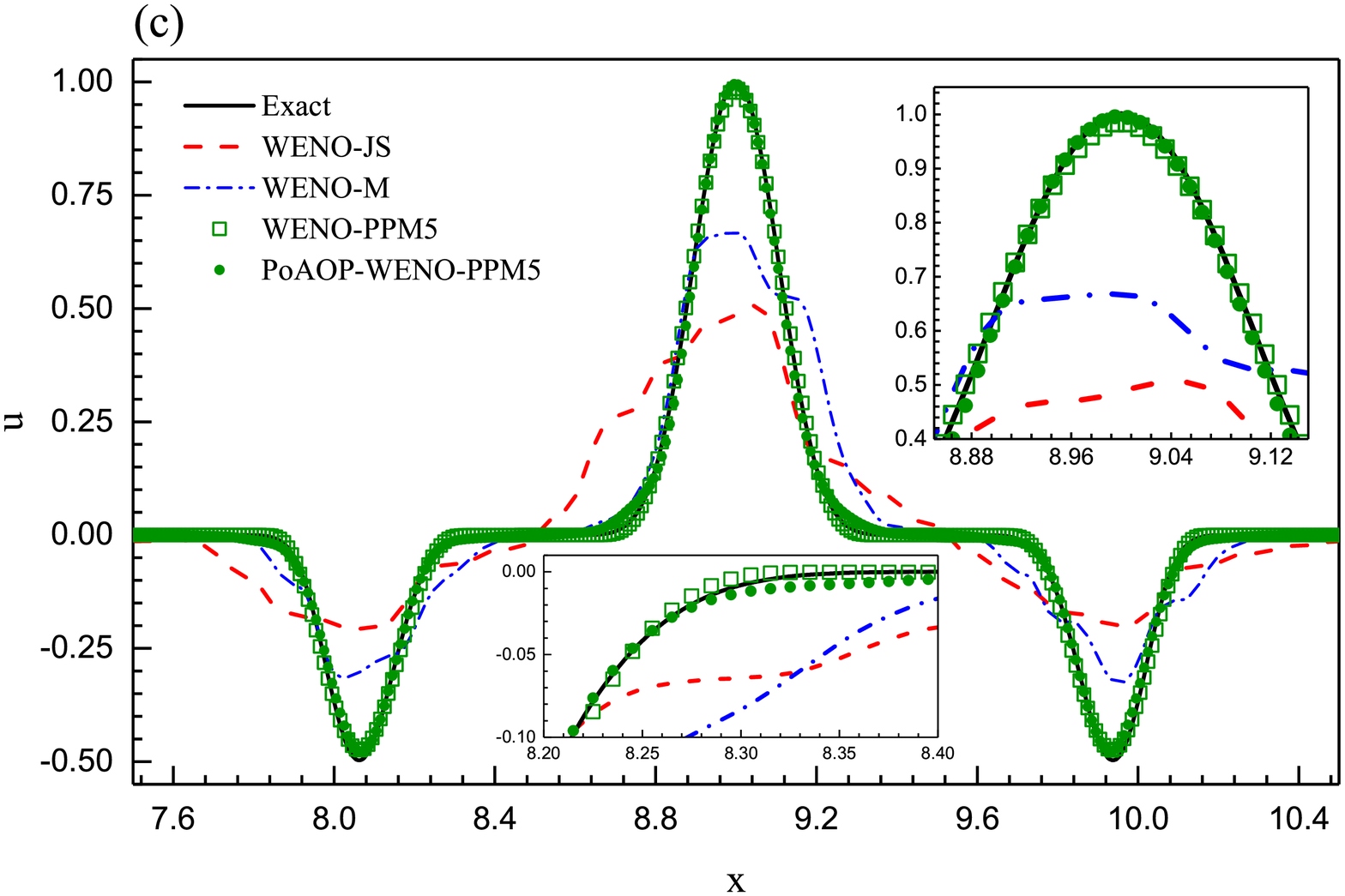}  
  \includegraphics[height=0.33\textwidth]
  {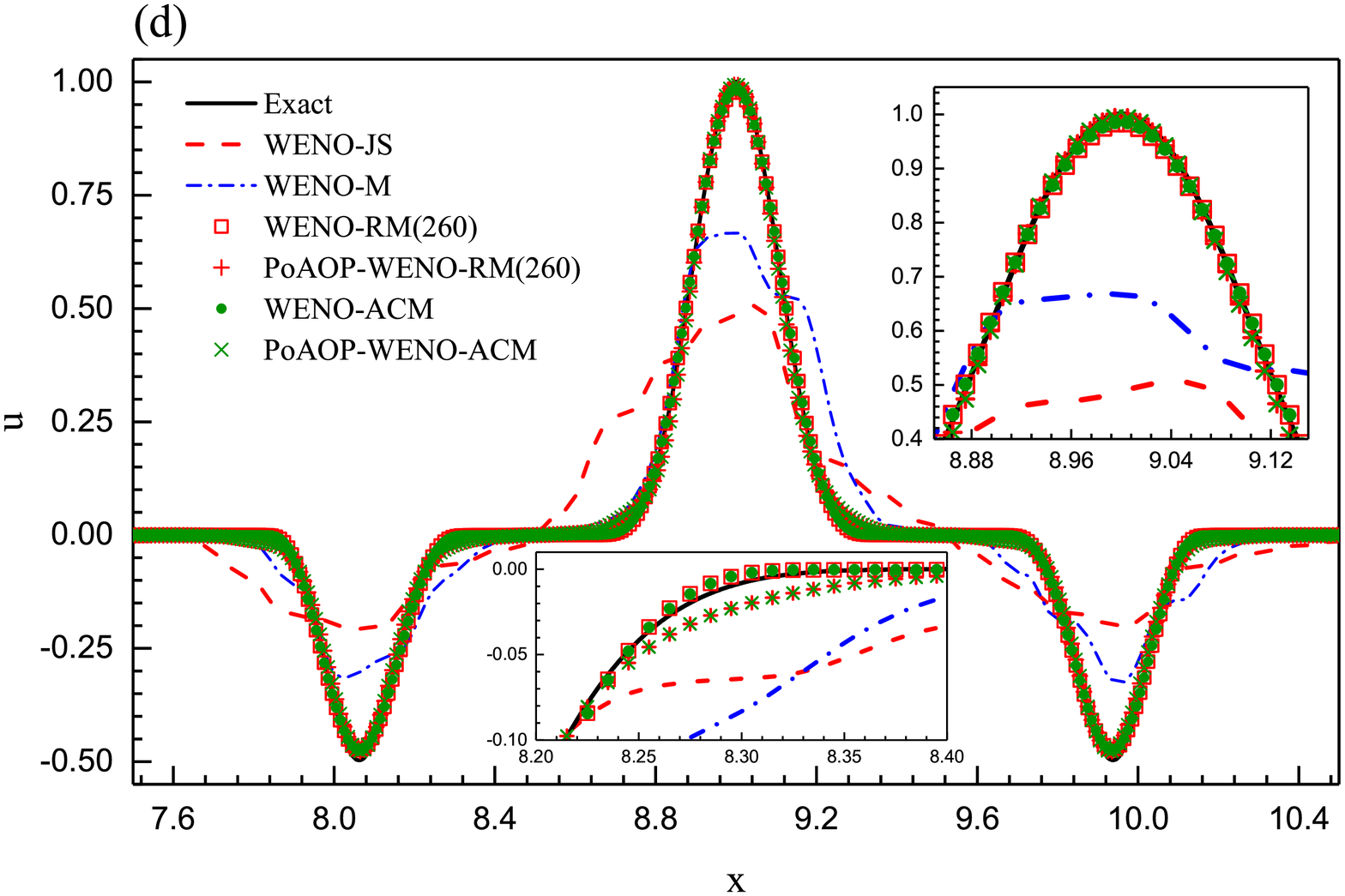}   
  \caption{Performance of various WENO schemes for Example 
  \ref{ex:long-run:case1} at output time $t = 1200$ with $N = 300$.}
\label{fig:ex:long-run:N300}
\end{figure}

Fig. \ref{fig:ex:long-run:N300} shows the performances of various 
considered schemes at output time $t = 1200$ with the grid number 
of $N = 300$. From Fig. \ref{fig:ex:long-run:N300}, we can find 
that: (1) the WENO-JS scheme shows the lowest resolution, followed 
by the WENO-M scheme whose resolution is far lower than its 
associated LOP-WENO-M scheme; (2) the other LOP-WENO-X schemes 
give results with slightly lower resolutions than their associated 
WENO-X schemes but they still show far better resolutions than the 
WENO-M and WENO-JS schemes.

\subsubsection{With discontinuities}
To show the advantage of the LOP-WENO-X schemes that they can not 
only preserve high resolutions but also prevent spurious 
oscillations especially for long output time computations, we solve 
Eq.\eqref{eq:LAE} with the periodic boundary condition by 
considering the following initial condition.

\begin{example}
The initial condition is given by
\begin{equation}
\begin{array}{l}
u(x, 0) = \left\{
\begin{array}{ll}
1.0,   & x \in [-1.0, 0.0], \\
0.0,   & x \in (0.0, 1.0].
\end{array}\right. 
\end{array}
\label{eq:LAE:Z}
\end{equation}
\label{ex:AccuracyTest:Z}
\end{example}

It simply consists of two constant states separated by sharp 
discontinuities at $x = 0.0, \pm 1.0$. 

Firstly, we examine the convergence properties of the considered 
schemes with the output time $t =2000$. Here, the CFL number is taken to be $0.1$. For the purpose of comparison, we also present the results computed by the WENO5-ILW scheme.

We give the $L_{1}$, $L_{\infty}$ errors and 
the corresponding convergence orders in Table 
\ref{table:AccuracyTest:Z:t2000}. We can see that: (1) the WENO-JS scheme produces significantly larger 
numerical errors than all other schemes and this indicates that it 
has the highest dissipation among all schemes; (2) the numerical 
errors generated by the LOP-WENO-M scheme are much smaller than 
its associated WENO-M scheme, especially for the $L_{1}$ errors 
for the computing cases of $N=400, 800$, and this demonstrates the 
advantage of the LOP-WENO-M scheme of decreasing the dissipation; 
(3) the $L_{1}$ orders of the other mapped WENO-X schemes are 
clearly lower than those of their associated LOP-WENO-X schemes 
although their corresponding numerical errors are slightly smaller; 
(4) the $L_{\infty}$ errors of the LOP-WENO-X schemes 
are very close to, or even smaller for many cases than, their 
associated mapped WENO-X schemes. Moreover, if we take a view of the 
$x - u$ profiles, we can find that the resolution of the result 
computed by the WENO-M scheme is significantly lower than that of 
the LOP-WENO-M scheme, and the other mapped WENO-X schemes 
generate spurious oscillations but their associated LOP-WENO-X 
schemes do not. To manifest this, detailed tests will be conducted 
and the solutions will be presented carefully in the following pages.

\begin{table}[!ht]
\begin{myFontSize}
\centering
\caption{Numerical errors and convergence orders of accuracy on 
Example \ref{ex:AccuracyTest:Z} at $t = 2000.0$.}
\label{table:AccuracyTest:Z:t2000}
\begin{tabular*}{\hsize}
{@{}@{\extracolsep{\fill}}cllllllll@{}}
\hline
\space    &\multicolumn{4}{l}{\cellcolor{gray!35}{WENO5-ILW}}  
          &\multicolumn{4}{l}{\cellcolor{gray!35}{WENO-JS}}  \\
\cline{2-5}  \cline{6-9}
$N$                   & $L_{1}$ error      & $L_{1}$ order 
                      & $L_{\infty}$ error & $L_{\infty}$ order
		  			  & $L_{1}$ error      & $L_{1}$ order 
      				  & $L_{\infty}$ error & $L_{\infty}$ order \\
\Xhline{0.65pt}
200                   & 1.03240E-01        & -
                      & 4.67252E-01        & -
                      & 4.48148E-01        & -
                      & 5.55748E-01        & - \\
400                   & 5.79848E-02        & 0.8323
                      & 4.70837E-01        & -0.0110
                      & 3.37220E-01        & 0.4103
                      & 5.77105E-01        & -0.0544 \\
800                   & 3.25843E-02        & 0.8315
                      & 4.74042E-01        & -0.0098
                      & 2.93752E-01        & 0.1991
                      & 5.17829E-01        & 0.1564 \\
\hline
\space    &\multicolumn{4}{l}{\cellcolor{gray!35}{WENO-M}}
          &\multicolumn{4}{l}{\cellcolor{gray!35}{LOP-WENO-M}}\\
\cline{2-5}  \cline{6-9}
$N$                   & $L_{1}$ error      & $L_{1}$ order 
                      & $L_{\infty}$ error & $L_{\infty}$ order
		  			  & $L_{1}$ error      & $L_{1}$ order 
      				  & $L_{\infty}$ error & $L_{\infty}$ order \\
\Xhline{0.65pt}
200                   & 1.76398E-01        & -
                      & 5.27583E-01        & -
                      & 1.22201E-01        & -
                      & 5.04793E-01        & - \\
400                   & 1.67082E-01        & 0.0783
                      & 5.73328E-01        & -0.1200
                      & 6.77592E-02        & 0.8508
                      & 4.88315E-01        & 0.0479 \\
800                   & 2.00760E-01        & -0.2649
                      & 5.47150E-01        & 0.0674
                      & 3.67281E-02        & 0.8835
                      & 4.90550E-01        & -0.0066 \\
\hline
\space    &\multicolumn{4}{l}{\cellcolor{gray!35}{WENO-PM6}} 
          &\multicolumn{4}{l}{\cellcolor{gray!35}{LOP-WENO-PM6}}\\
\cline{2-5}  \cline{6-9}
$N$                   & $L_{1}$ error      & $L_{1}$ order 
                      & $L_{\infty}$ error & $L_{\infty}$ order
		  			  & $L_{1}$ error      & $L_{1}$ order 
      				  & $L_{\infty}$ error & $L_{\infty}$ order \\
\Xhline{0.65pt}
200                   & 8.67541E-02        & -
                      & 5.02070E-01        & -
                      & 1.19011E-01        & -
                      & 4.75985E-01        & - \\
400                   & 5.29105E-02        & 0.7134
                      & 5.09366E-01        & -0.0208
                      & 6.45626E-02        & 0.8823
                      & 4.95054E-01        & -0.0567 \\
800                   & 2.97704E-02        & 0.8297
                      & 5.15102E-01        & -0.0162
                      & 3.52222E-02        & 0.8742
                      & 4.75635E-01        & 0.0577 \\
\hline
\space    &\multicolumn{4}{l}{\cellcolor{gray!35}{WENO-IM(2, 0.1)}} 
    &\multicolumn{4}{l}{\cellcolor{gray!35}{LOP-WENO-IM(2, 0.1)}}\\
\cline{2-5}  \cline{6-9}
$N$                   & $L_{1}$ error      & $L_{1}$ order 
                      & $L_{\infty}$ error & $L_{\infty}$ order
		  			  & $L_{1}$ error      & $L_{1}$ order 
      				  & $L_{\infty}$ error & $L_{\infty}$ order \\
\Xhline{0.65pt}
200                   & 7.94092E-02        & -
                      & 4.64949E-01        & -
                      & 1.22302E-01        & -
                      & 5.08308E-01        & - \\
400                   & 4.61209E-02        & 0.7839
                      & 4.76074E-01        & -0.0341
                      & 6.64627E-02        & 0.8798
                      & 5.02003E-01        & 0.0180 \\
800                   & 2.65533E-02        & 0.7965
                      & 4.91316E-01        & -0.0455
                      & 3.61408E-02        & 0.8789
                      & 4.79591E-01        & 0.0659 \\
\hline
\space    &\multicolumn{4}{l}{\cellcolor{gray!35}{WENO-PPM5}} 
          &\multicolumn{4}{l}{\cellcolor{gray!35}{LOP-WENO-PPM5}}\\
\cline{2-5}  \cline{6-9}
$N$                   & $L_{1}$ error      & $L_{1}$ order 
                      & $L_{\infty}$ error & $L_{\infty}$ order
		  			  & $L_{1}$ error      & $L_{1}$ order 
      				  & $L_{\infty}$ error & $L_{\infty}$ order \\
\Xhline{0.65pt}
200                   & 9.20390E-02        & -
                      & 4.99999E-01        & -
                      & 1.17886E-01        & -
                      & 4.84251E-01        & - \\
400                   & 5.27679E-02        & 0.8026
                      & 5.07952E-01        & -0.0228
                      & 6.58012E-02        & 0.8412
                      & 5.04572E-01        & -0.0593 \\
800                   & 2.96879E-02        & 0.8298
                      & 5.14059E-01        & -0.0172
                      & 3.58152E-02        & 0.8775
                      & 4.99765E-01        & 0.0138 \\
\hline
\space    &\multicolumn{4}{l}{\cellcolor{gray!35}{WENO-RM(260)}} 
      &\multicolumn{4}{l}{\cellcolor{gray!35}{LOP-WENO-RM(260)}}\\
\cline{2-5}  \cline{6-9}
$N$                   & $L_{1}$ error      & $L_{1}$ order 
                      & $L_{\infty}$ error & $L_{\infty}$ order
		  			  & $L_{1}$ error      & $L_{1}$ order 
      				  & $L_{\infty}$ error & $L_{\infty}$ order \\
\Xhline{0.65pt}
200                   & 8.64542E-02        & -
                      & 5.02486E-01        & -
                      & 1.19069E-01        & -
                      & 5.09991E-01        & - \\
400                   & 5.17965E-02        & 0.7391
                      & 5.08770E-01        & -0.0179
                      & 6.58446E-02        & 0.8547
                      & 5.02010E-01        & 0.0228 \\
800                   & 2.91482E-02        & 0.8294
                      & 5.14009E-01        & -0.0148
                      & 3.63654E-02        & 0.8565
                      & 4.78897E-01        & 0.0680 \\
\hline
\space    &\multicolumn{4}{l}{\cellcolor{gray!35}{WENO-ACM}} 
          &\multicolumn{4}{l}{\cellcolor{gray!35}{LOP-WENO-ACM}}\\
\cline{2-5}  \cline{6-9}
$N$                   & $L_{1}$ error      & $L_{1}$ order 
                      & $L_{\infty}$ error & $L_{\infty}$ order
		  			  & $L_{1}$ error      & $L_{1}$ order 
      				  & $L_{\infty}$ error & $L_{\infty}$ order \\
\Xhline{0.65pt}
200                   & 8.87640E-02        & -
                      & 5.06230E-01        & -
                      & 1.21982E-01        & -
                      & 5.14204E-01        & - \\
400                   & 5.16217E-02        & 0.7820
                      & 5.11512E-01        & -0.0150
                      & 6.55457E-02        & 0.8961
                      & 4.98088E-01        & 0.0459 \\
800                   & 2.94211E-02        & 0.8111
                      & 5.15990E-01        & -0.0126
                      & 3.61428E-02        & 0.8588
                      & 4.79224E-01        & 0.0557 \\
\hline
\end{tabular*}
\end{myFontSize}
\end{table}

To provide a better illustration, we re-calculate Example \ref{ex:AccuracyTest:Z} by considered WENO schemes with the output time $t=200$ using 
the uniform meshes of $N = 1600$ and $N = 3200$, respectively. 

\begin{figure}[!ht]
\centering
\includegraphics[height=0.32\textwidth]
{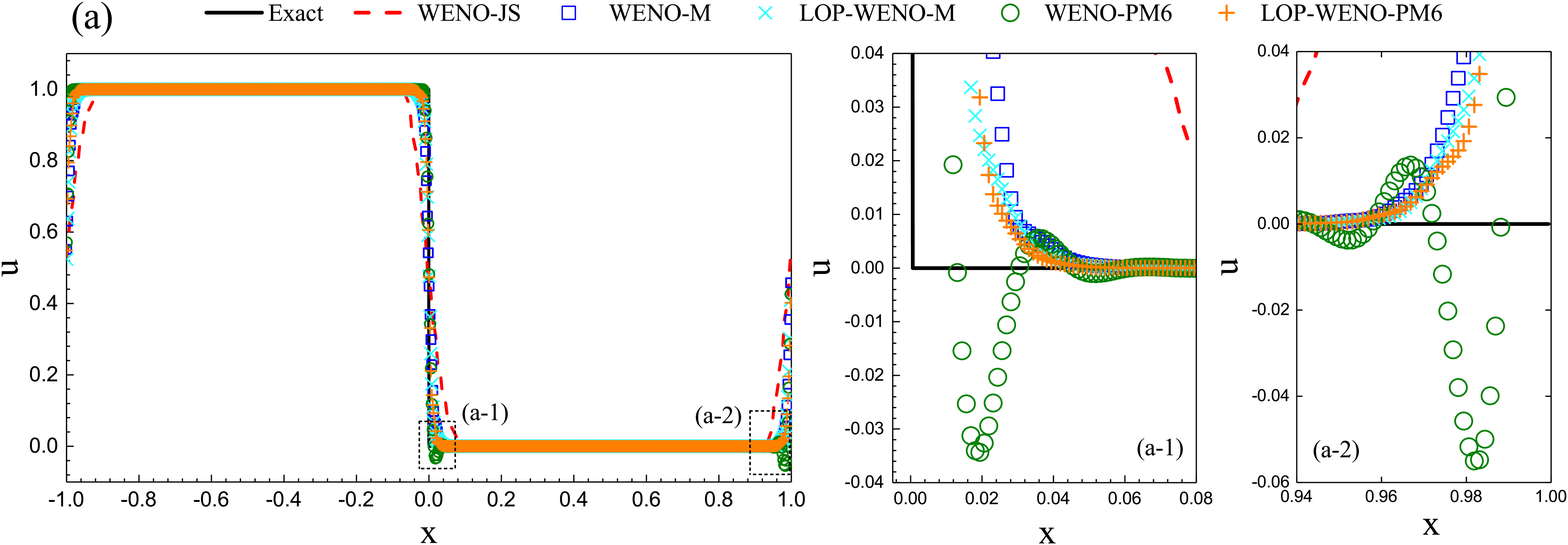}\\ \vspace{0.8cm}
\includegraphics[height=0.32\textwidth]
{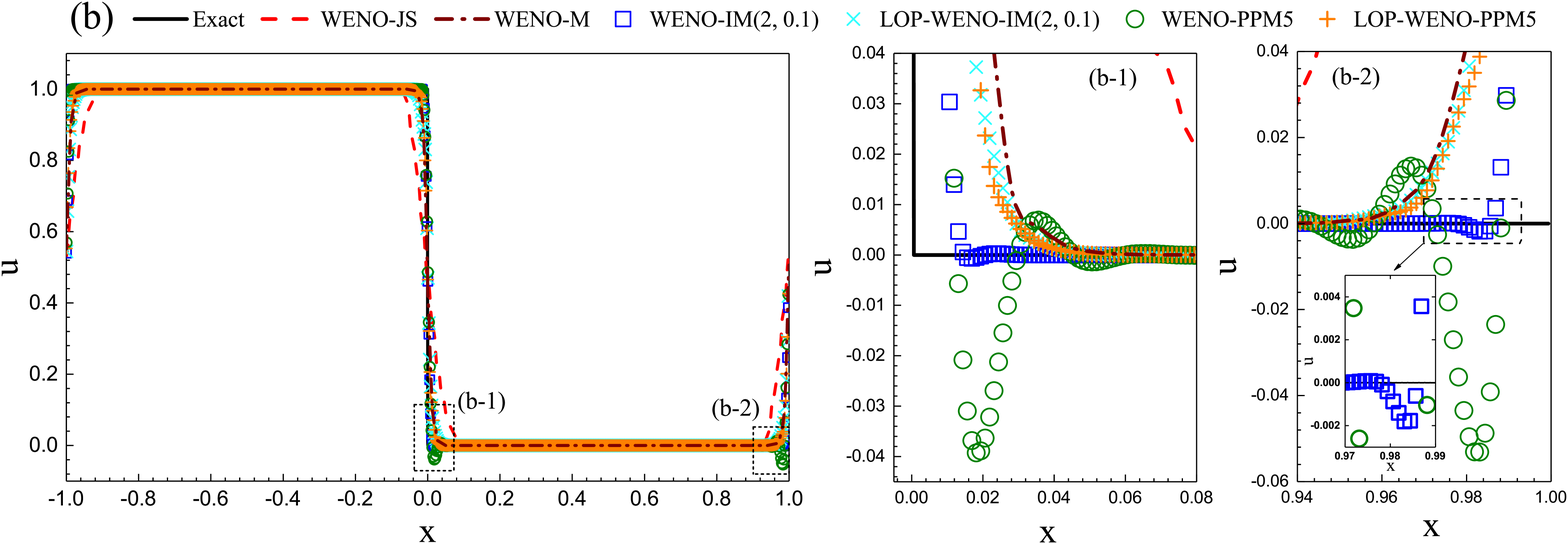}\\ \vspace{0.8cm}
\includegraphics[height=0.32\textwidth]
{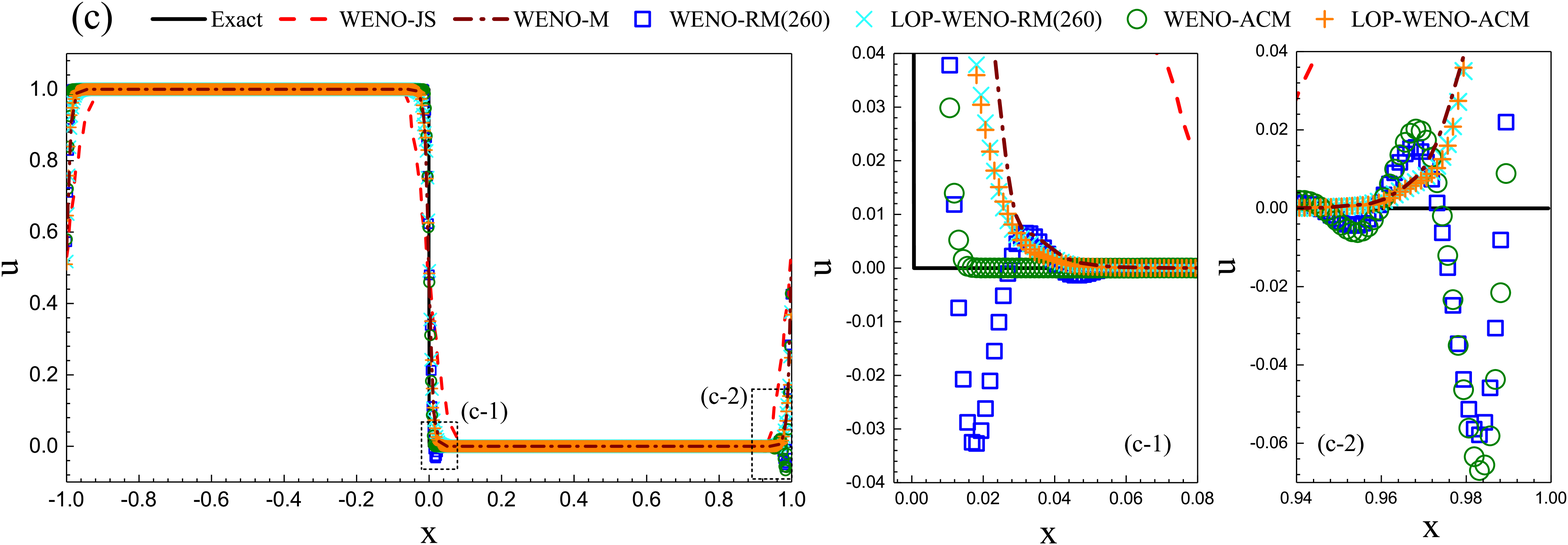}
\caption{Performance of the considered WENO schemes for Example 
\ref{ex:AccuracyTest:Z} at output time $t = 200$ with $N = 1600$.}
\label{fig:Z:N1600}
\end{figure}

Fig. \ref{fig:Z:N1600} shows the comparison of considered schemes with $t = 200$ and $N = 1600$. We can observe that: (1) all the 
LOP-WENO-X schemes provide the numerical results with 
significantly higher resolutions than those of the WENO-JS and 
WENO-M schemes, and moreover, they are all able to avoid the 
spurious oscillations that will be inevitably generated by most of 
their associated mapped WENO-X schemes; (2) it seems that 
the WENO-IM(2, 0.1) scheme almost does not generate spurious 
oscillations and it gains better resolutions than the 
LOP-WENO-IM(2, 0.1) scheme in most of the region; (3) however, if 
we take a closer look, we can see that the WENO-IM(2, 0.1) scheme 
also generates very slight spurious oscillations as shown in Fig. 
\ref{fig:Z:N1600}(b-2). 

\begin{figure}[!ht]
\centering
\includegraphics[height=0.32\textwidth]
{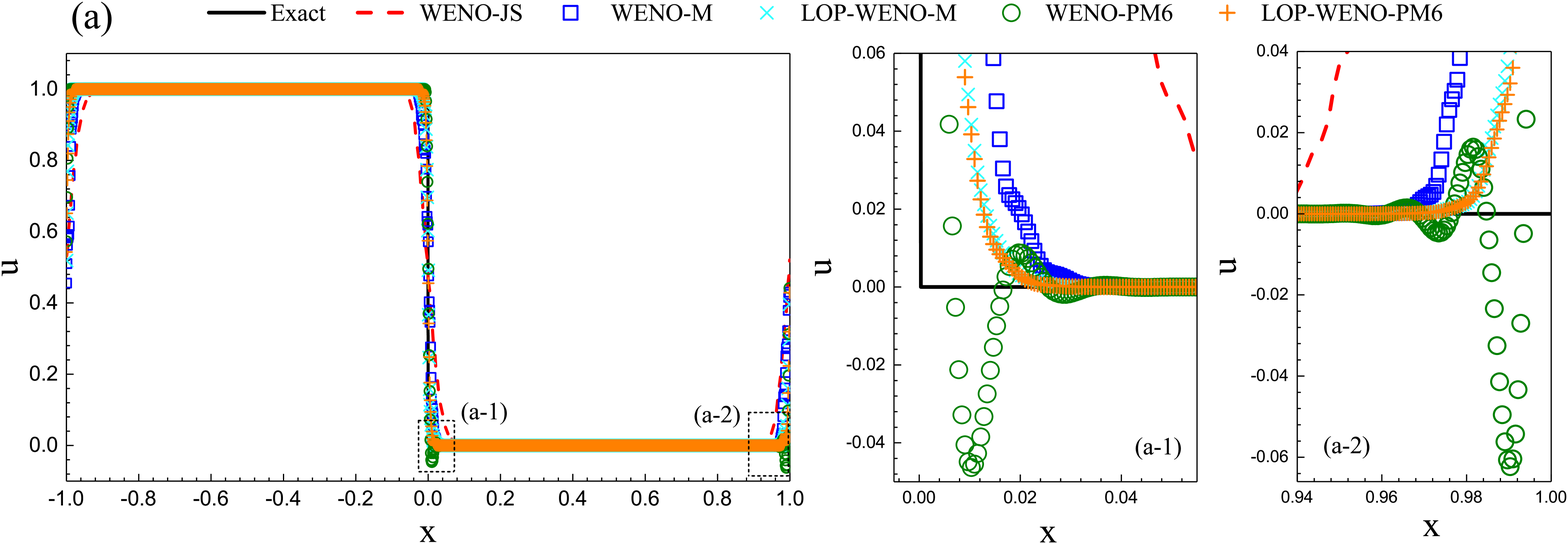}\\ \vspace{0.8cm}
\includegraphics[height=0.32\textwidth]
{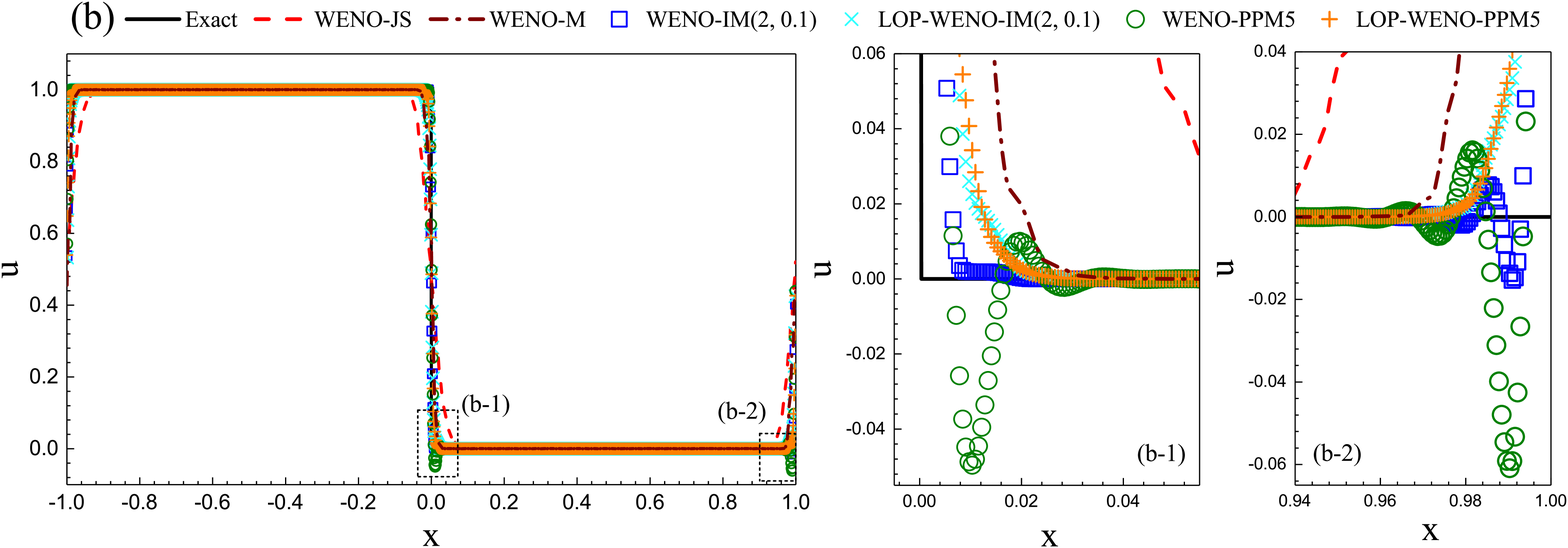}\\ \vspace{0.8cm}
\includegraphics[height=0.32\textwidth]
{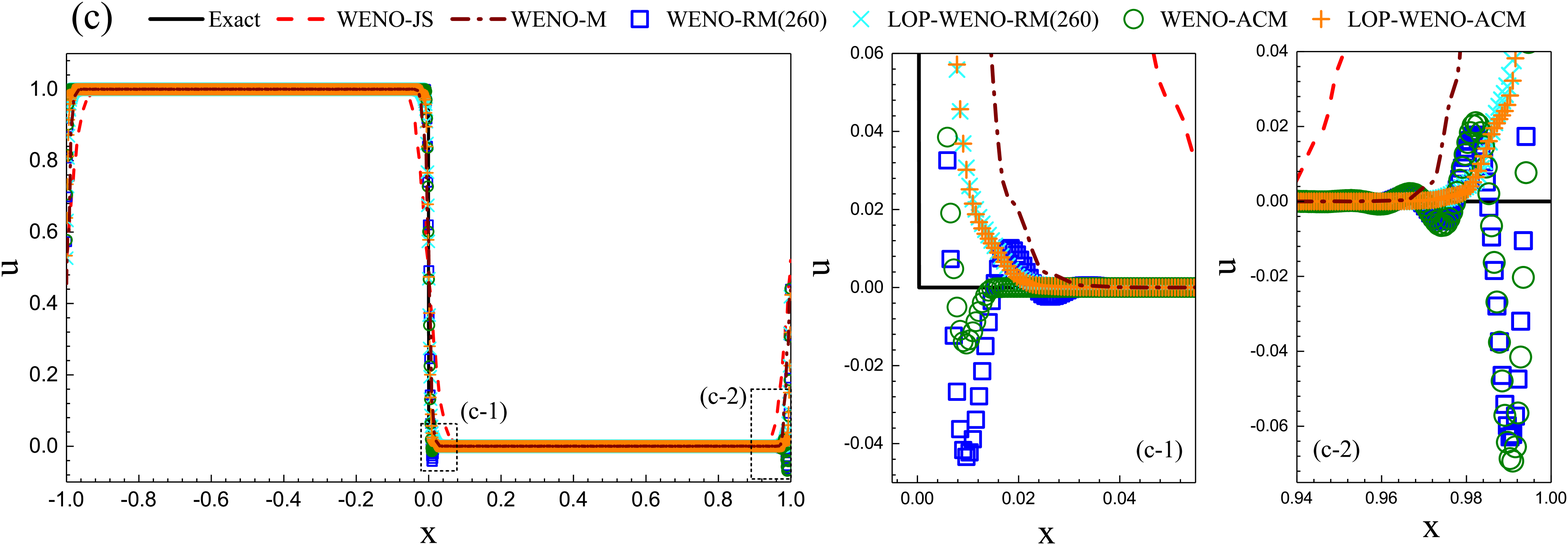}
\caption{Performance of the considered WENO schemes for Example 
\ref{ex:AccuracyTest:Z} at output time $t = 200$ with $N = 3200$.}
\label{fig:Z:N3200}
\end{figure}

The comparison of considered schemes for the case of $t = 200$ and 
$N = 3200$ are shown in Fig. \ref{fig:Z:N3200}. We can find that: (1) as the grid 
number increases, the spurious oscillations produced by the 
WENO-IM(2, 0.1) scheme become more violent and they are easily to be 
observed, however, the LOP-WENO-IM(2, 0.1) scheme can still prevent
the spurious oscillations but provide very high resolutions; (2) all 
the LOP-WENO-X schemes still evidently provide much better 
resolutions than those of the WENO-JS and WENO-M schemes; (3) as the 
grid number increases, the spurious oscillations generated by the 
WENO-X schemes appear to be closer to the discontinuities, and the 
amplitudes of these spurious oscillations become larger; (4) 
furthermore, even though the grid number increases, the LOP-WENO-X 
schemes can still avoid spurious oscillations but obtain the great 
improvement of the resolution.

Taken together, the results above suggest that the \textit{LOP} 
property, satisfied by the posteriori adaptive \textit{OP} mapping 
in the present study, can help to preserve high resolutions and 
meanwhile avoid spurious oscillations in the simulation of problems 
with discontinuities, especially for long output times. And this is 
a very important aspect of this paper.

\section{Numerical experiments for Euler equations}
\label{NumericalExperiments}
In this section, we apply the various considered schemes to solve 
the Euler equations. In all calculations, we choose $\epsilon$ in 
Eq.\eqref{eq:weights:WENO-JS} to be $10^{-40}$ as recommended in 
\cite{WENO-M,WENO-IM}. The local characteristic decomposition 
\cite{WENO-JS} is used when the WENO schemes are applied to solve 
the Euler system.

\subsection{Comparisons of the performances on simulating problems 
with high-frequency smooth waves}\label{subsec:1DEuler}
In this subsection, we compare the numerical results of the 
LOP-WENO-X schemes with those of the MOP-WENO-X, WENO-X and WENO-JS 
schemes on calculating problems with high-frequency smooth waves. 
The considered problems are governed by the following 
one-dimensional Euler systems of compressible gas dynamics
\begin{equation}
\dfrac{\partial \mathbf{U}}{\partial t}    + 
\dfrac{\partial \mathbf{F(U)}}{\partial x} = \mathbf{0},
\label{2DEulerEquations}
\end{equation}
with
\begin{equation*}
\begin{array}{l}
\mathbf{U} = \Big( \rho, \rho u, E \Big)^{\mathrm{T}}, \quad
\mathbf{F(U)} = \Big(\rho u, \rho u^{2} + p, u(E+p) 
\Big)^{\mathrm{T}}, 
\end{array}
\end{equation*}
where $\rho, u, p$ and $E$ are the density, velocity in the $x$ 
coordinate direction, pressure and total energy, respectively. The 
relation of pressure $p$ and total energy for ideal gases is defined 
by
\begin{equation*}
p = (\gamma - 1)\Big( E - \dfrac{1}{2}\rho u^{2}\Big), \quad
\gamma = 1.4.
\end{equation*}

\subsubsection{Shu-Osher problem}
\begin{example}
This problem was presented by Shu and Osher \cite{ENO-Shu1989}. The 
computational domain of $[-5, 5]$ is initialized by
\begin{equation}
\big( \rho, u, p \big)(x, 0) =\left\{
\begin{array}{ll}
(3.857143, 2.629369, 10.333333), & x \in [-5.0, -4.0], \\
(1.0 + 0.2\sin(5x), 0, 1), & x \in [-4.0, 5.0].
\end{array}\right.
 \label{initial:Shu-Osher}
\end{equation}
The transmissive boundary conditions are used at $x = \pm 5$, and 
the output time is set to be $t = 1.8$.
\label{ex:Shu-Osher}
\end{example}

\begin{figure}[!ht]
\centering
  \includegraphics[height=0.37\textwidth]{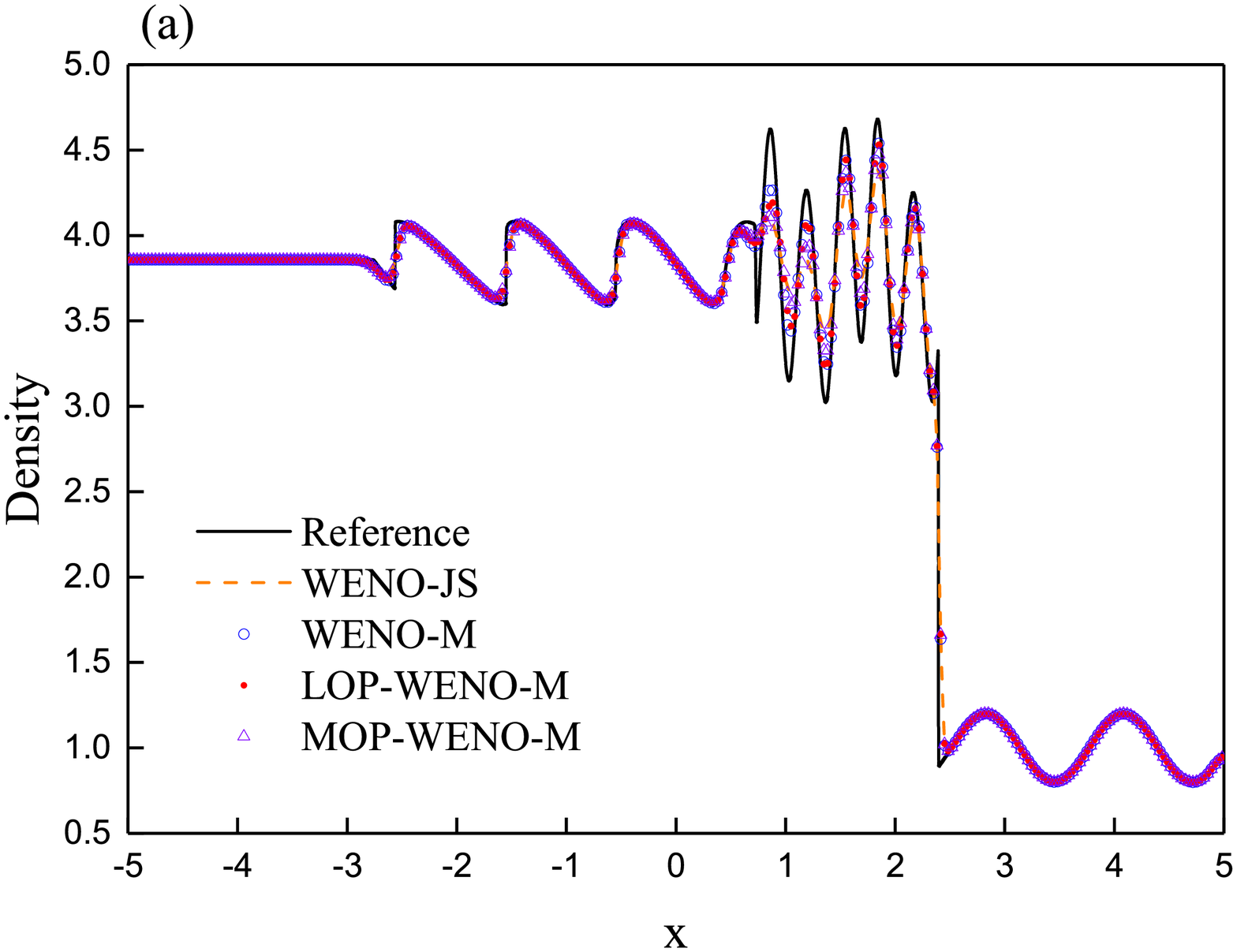}
  \includegraphics[height=0.37\textwidth]{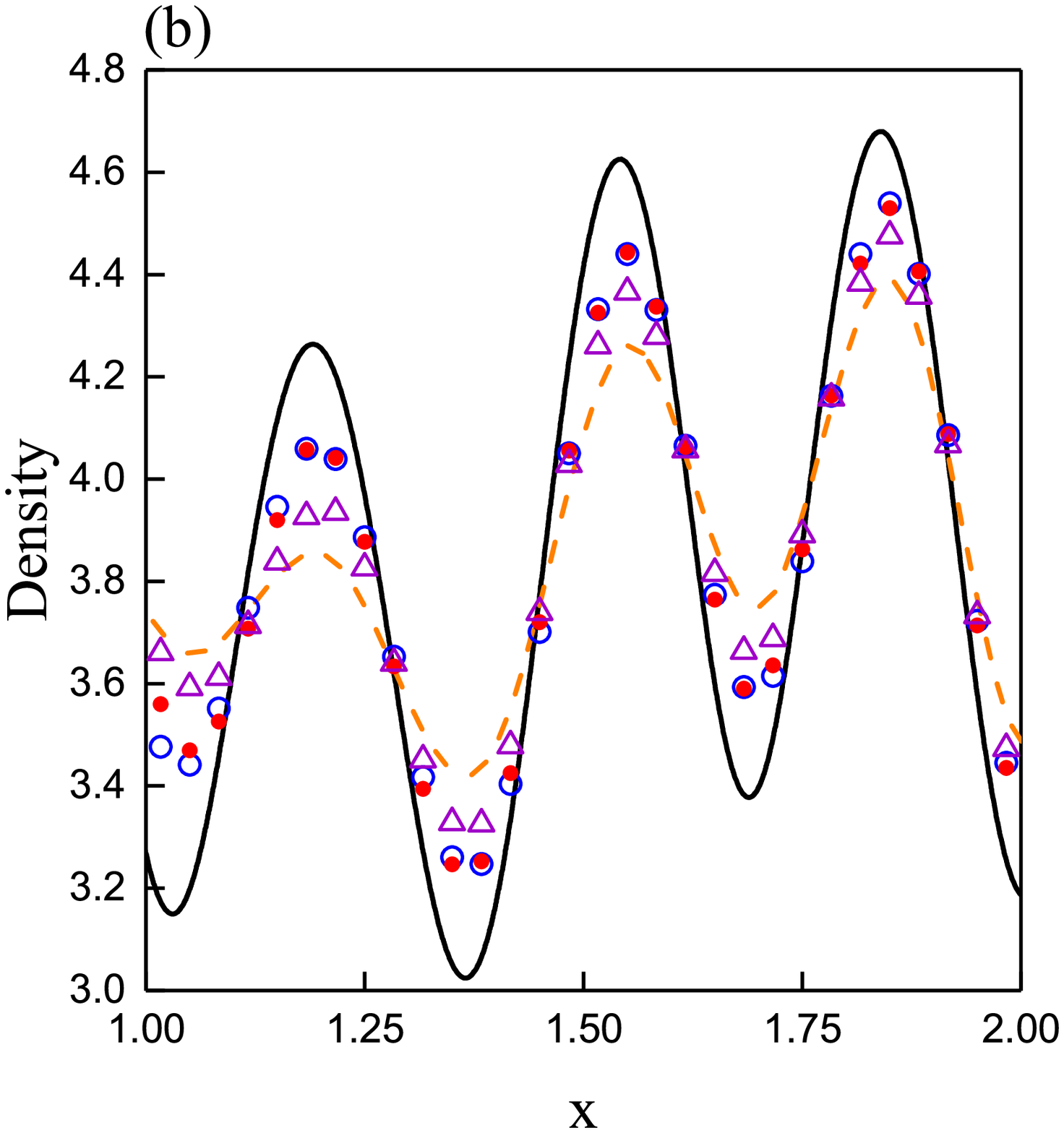} 
     \caption{Results of LOP-/MOP-/WENO-M and WENO-JS on solving the Shu-Osher problem.}
     \label{fig:Shu-Osher:WENO-M}
\end{figure}

\begin{figure}[!ht]
\centering
  \includegraphics[height=0.37\textwidth]{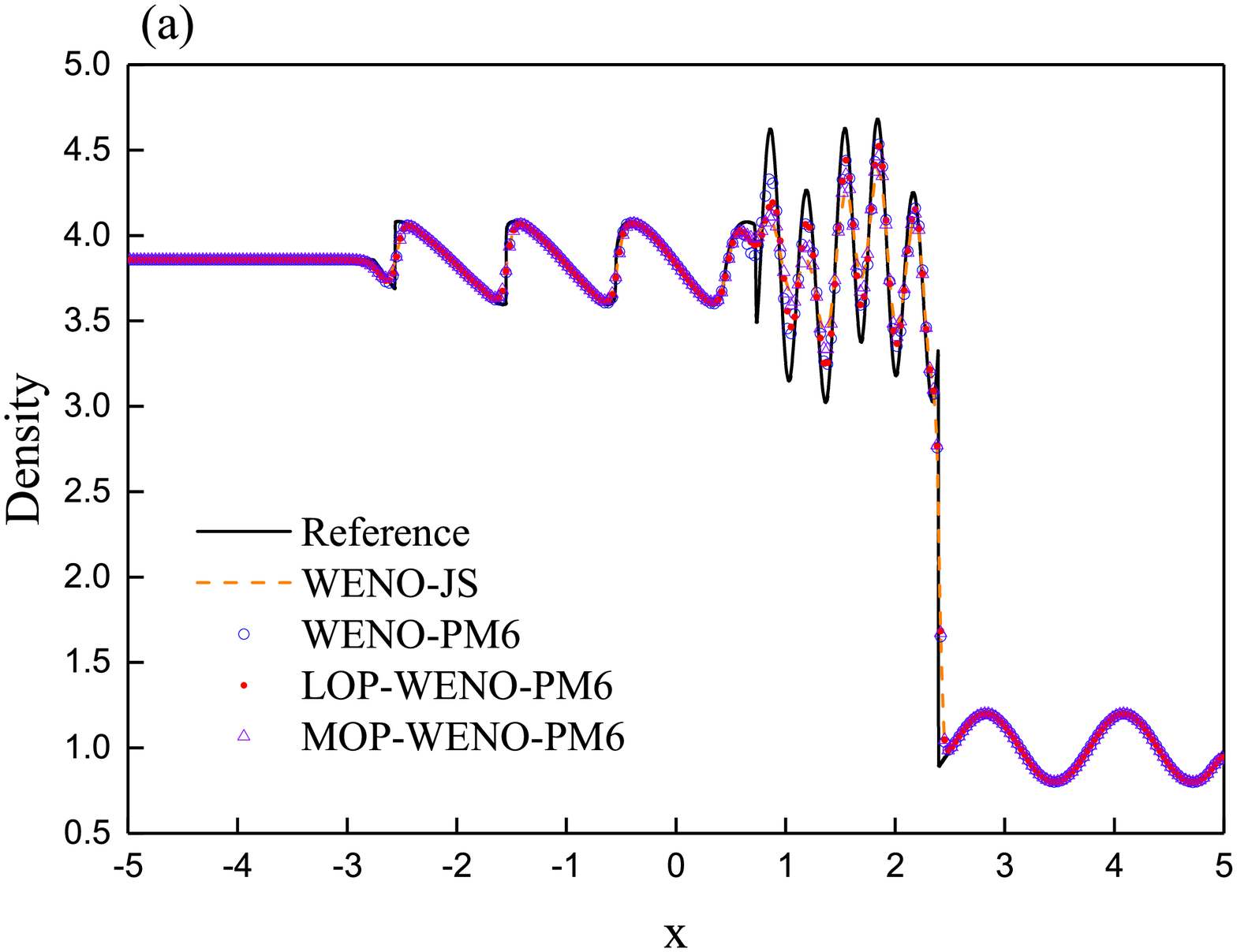}
  \includegraphics[height=0.37\textwidth]{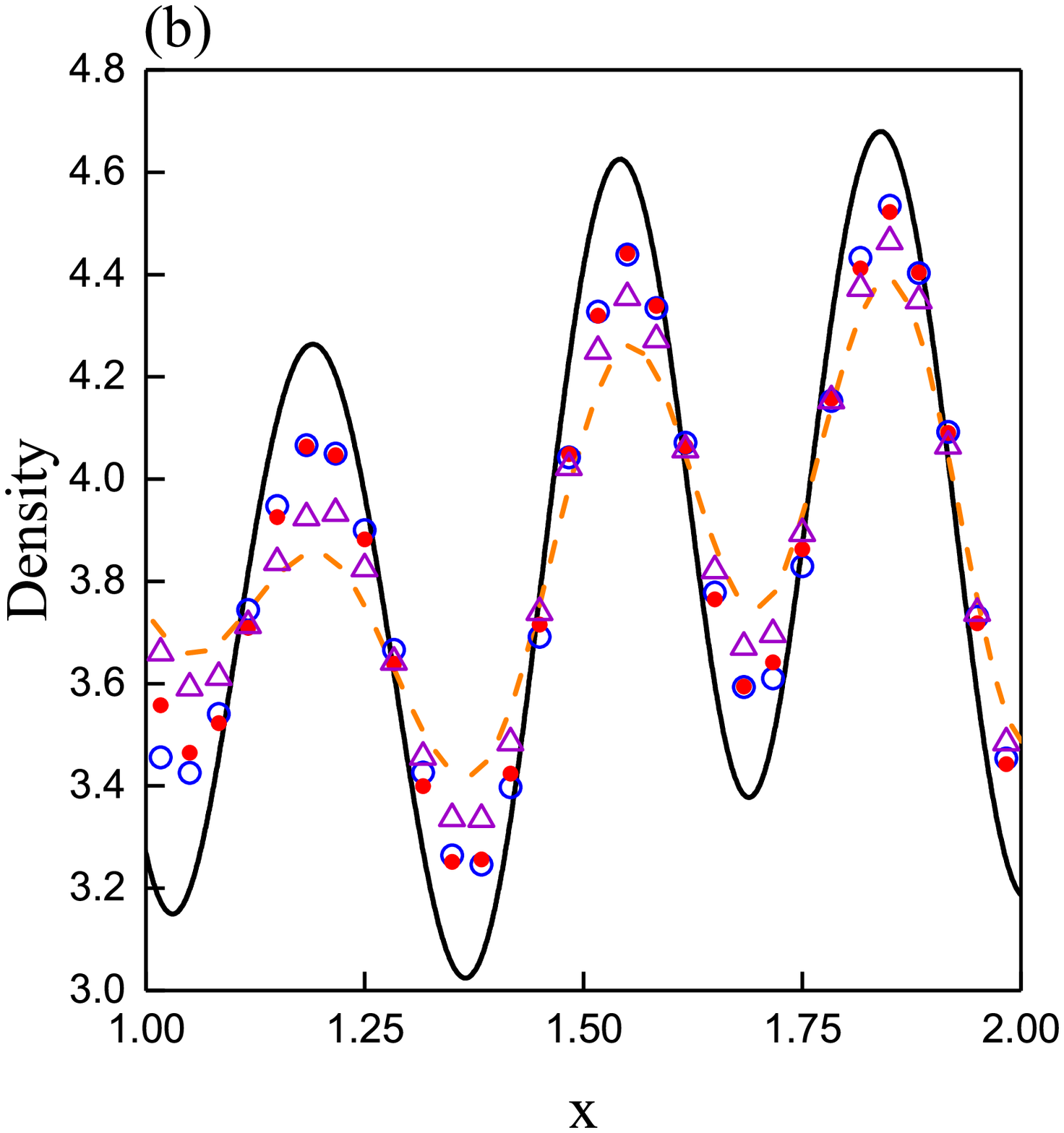} 
     \caption{Results of LOP-/MOP-/WENO-PM6 and WENO-JS on solving the Shu-Osher problem.}
     \label{fig:Shu-Osher:WENO-PM6}
\end{figure}

\begin{figure}[!ht]
\centering
  \includegraphics[height=0.37\textwidth]{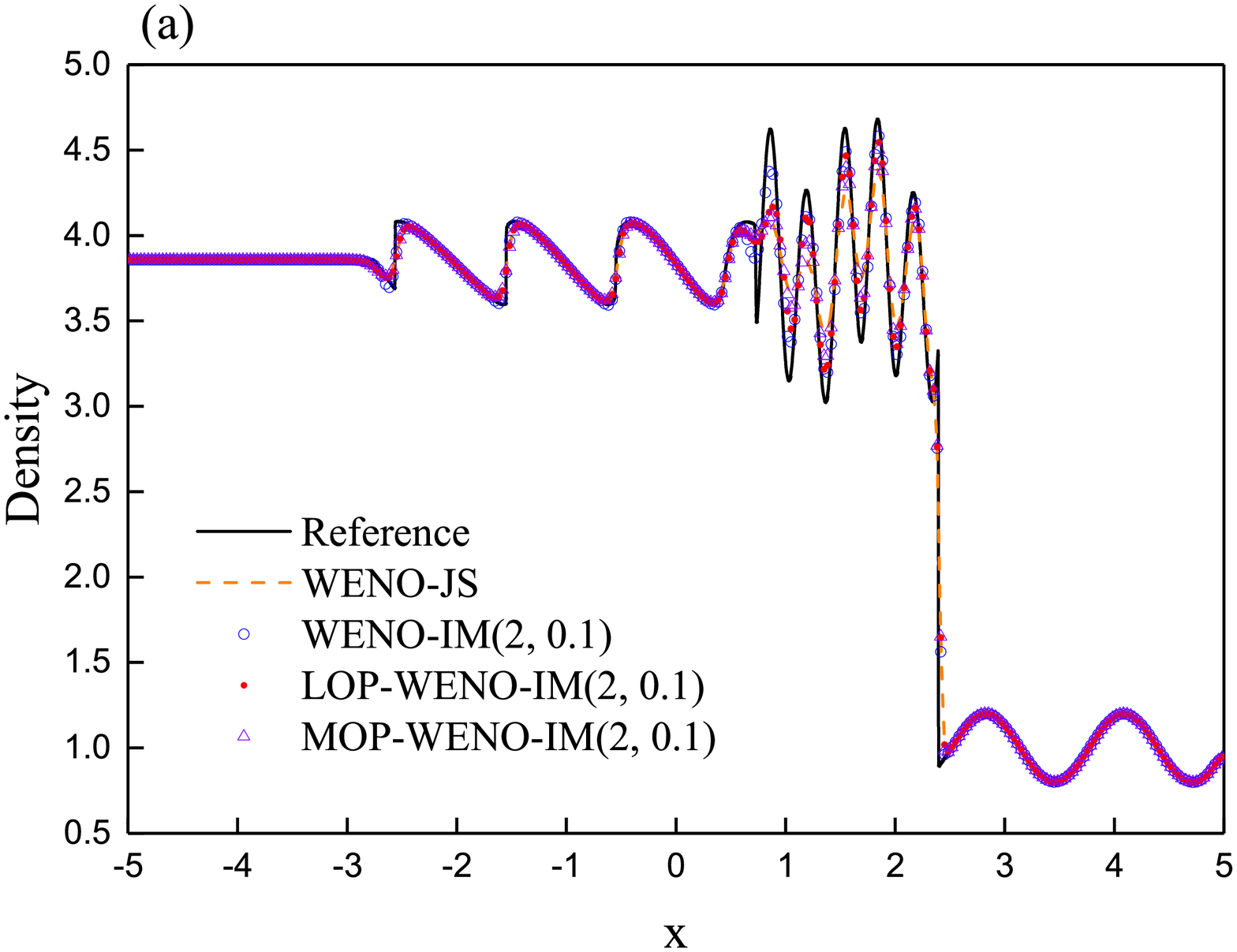}
  \includegraphics[height=0.37\textwidth]{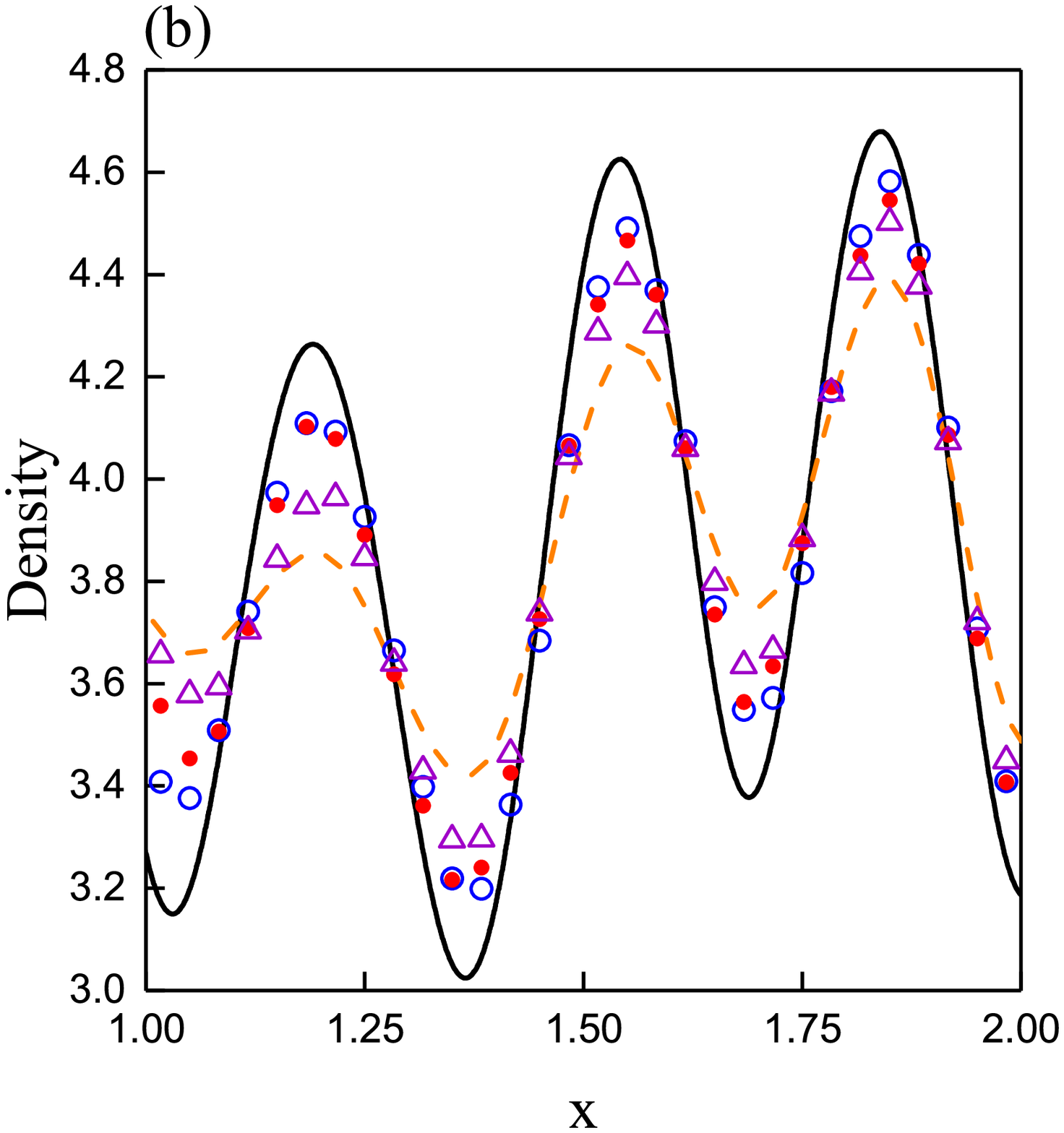}  
     \caption{Results of LOP-/MOP-/WENO-IM(2, 0.1) and WENO-JS on solving the Shu-Osher problem.}
     \label{fig:Shu-Osher:WENO-IM}
\end{figure}

\begin{figure}[!ht]
\centering
  \includegraphics[height=0.37\textwidth]{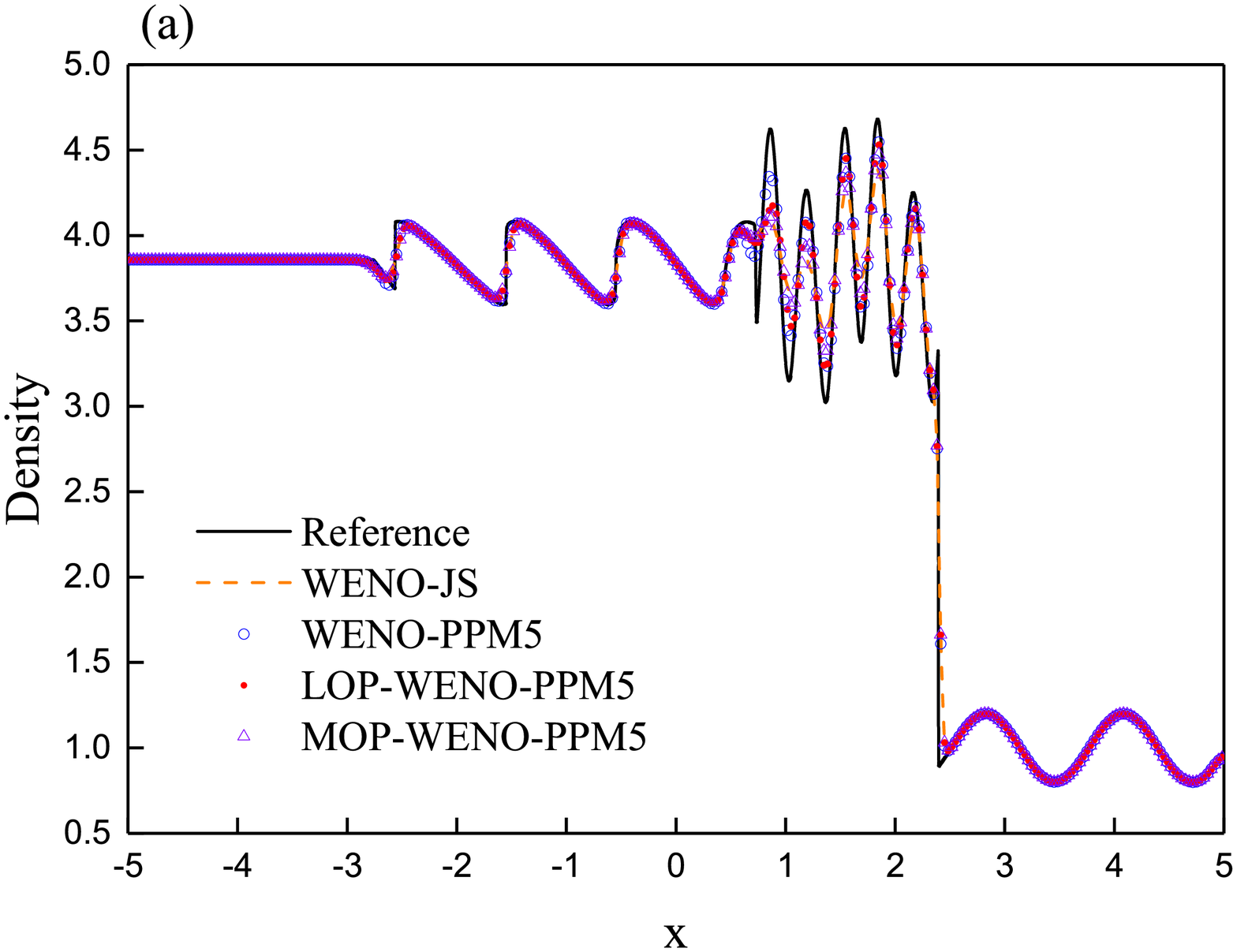}
  \includegraphics[height=0.37\textwidth]{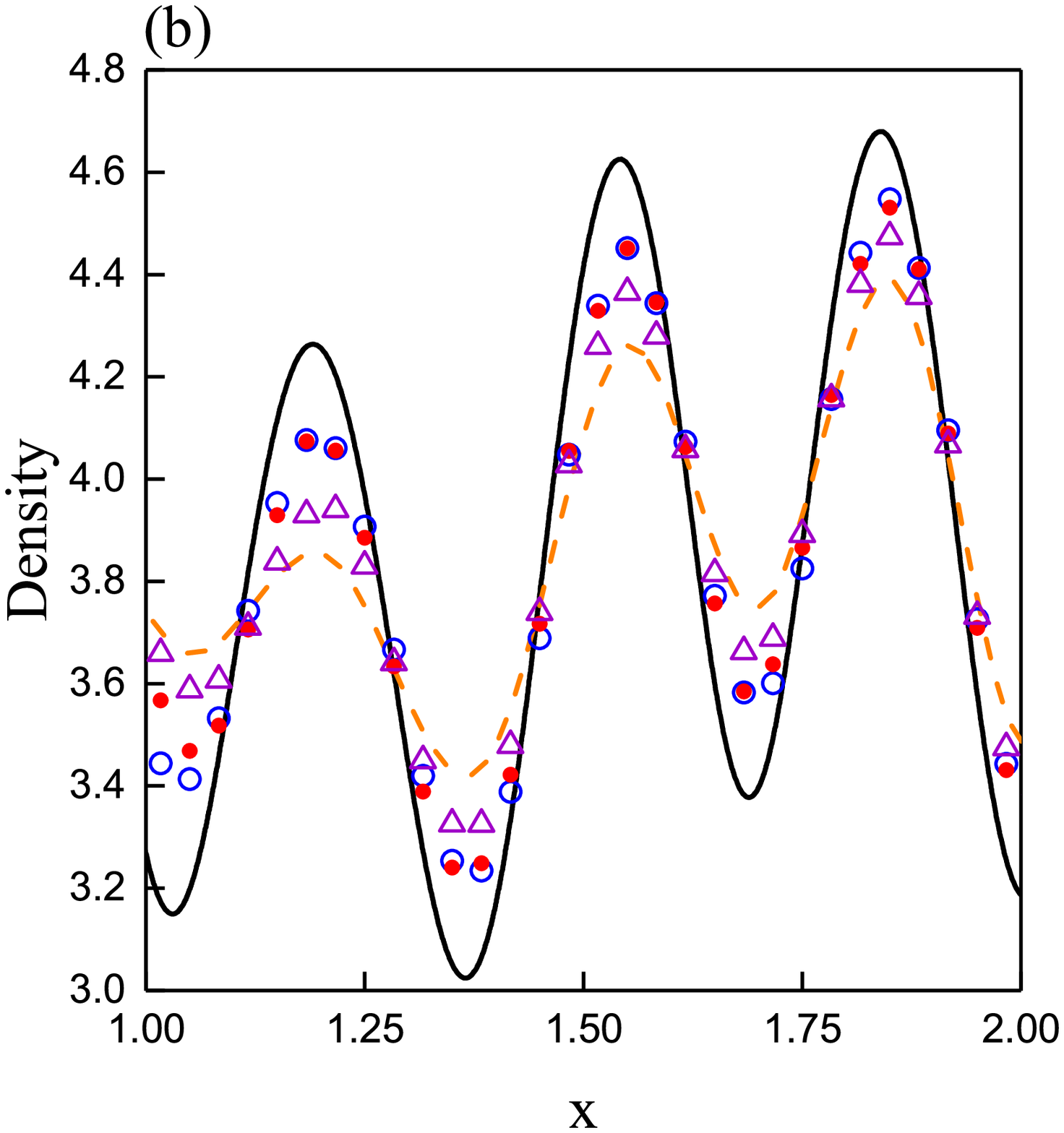}  
     \caption{Results of LOP-/MOP-/WENO-PPM5 and WENO-JS on solving the Shu-Osher problem.}
     \label{fig:Shu-Osher:WENO-PPM5}
\end{figure}

\begin{figure}[!ht]
\centering
  \includegraphics[height=0.37\textwidth]{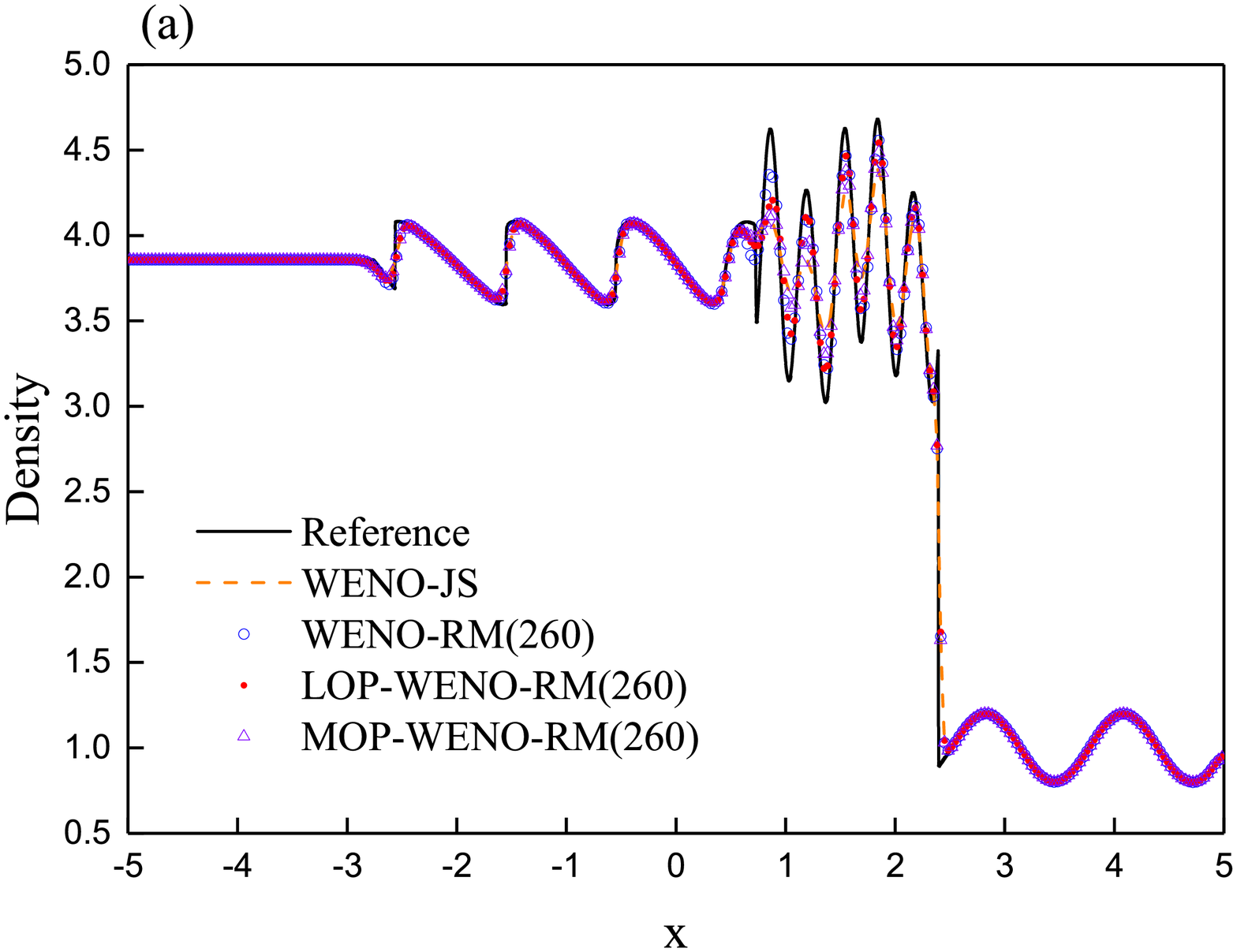}
  \includegraphics[height=0.37\textwidth]{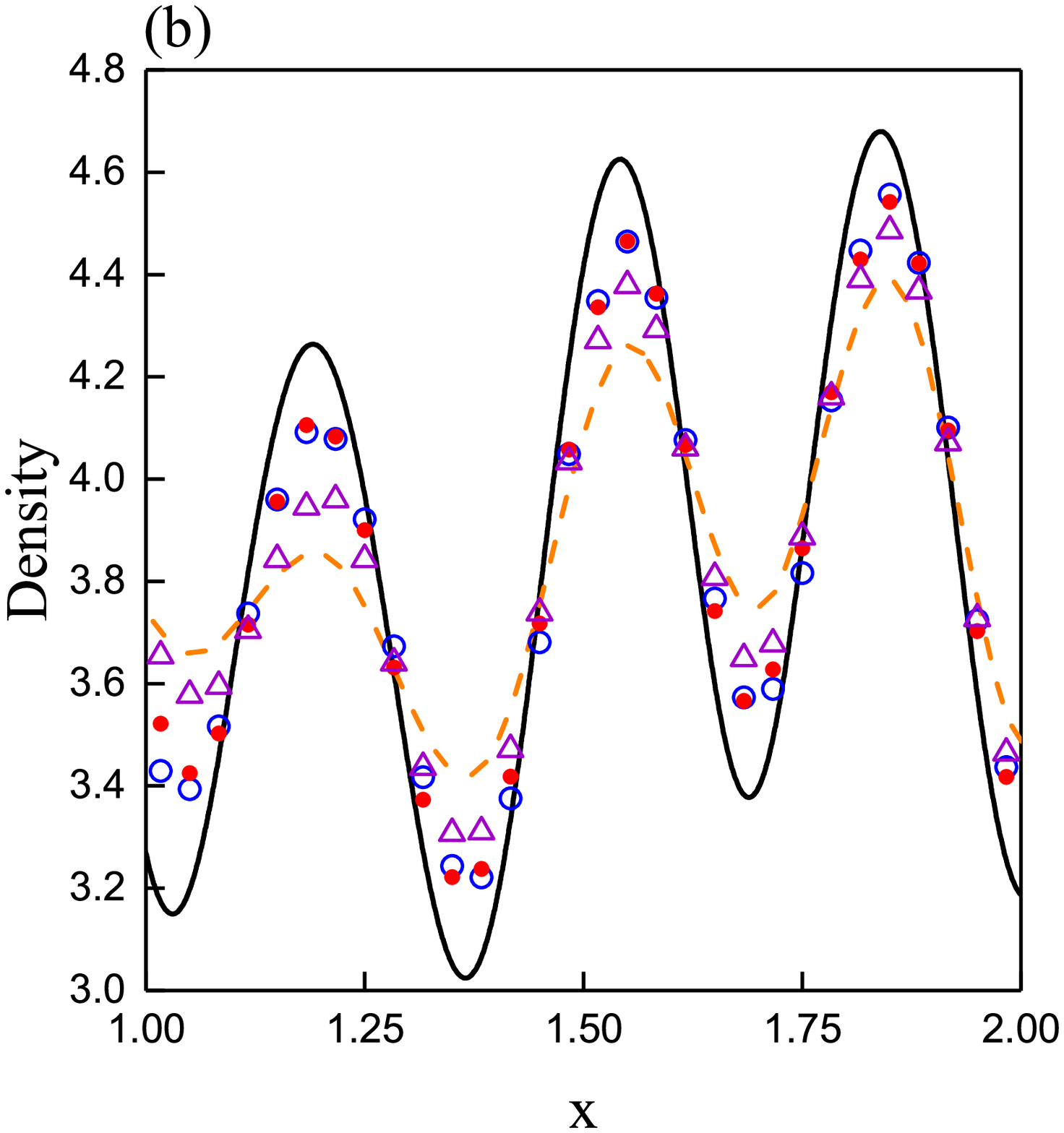} 
     \caption{Results of LOP-/MOP-/WENO-RM(260) and WENO-JS on solving the Shu-Osher problem.}
     \label{fig:Shu-Osher:WENO-RM260}
\end{figure}

\begin{figure}[!ht]
\centering
  \includegraphics[height=0.37\textwidth]{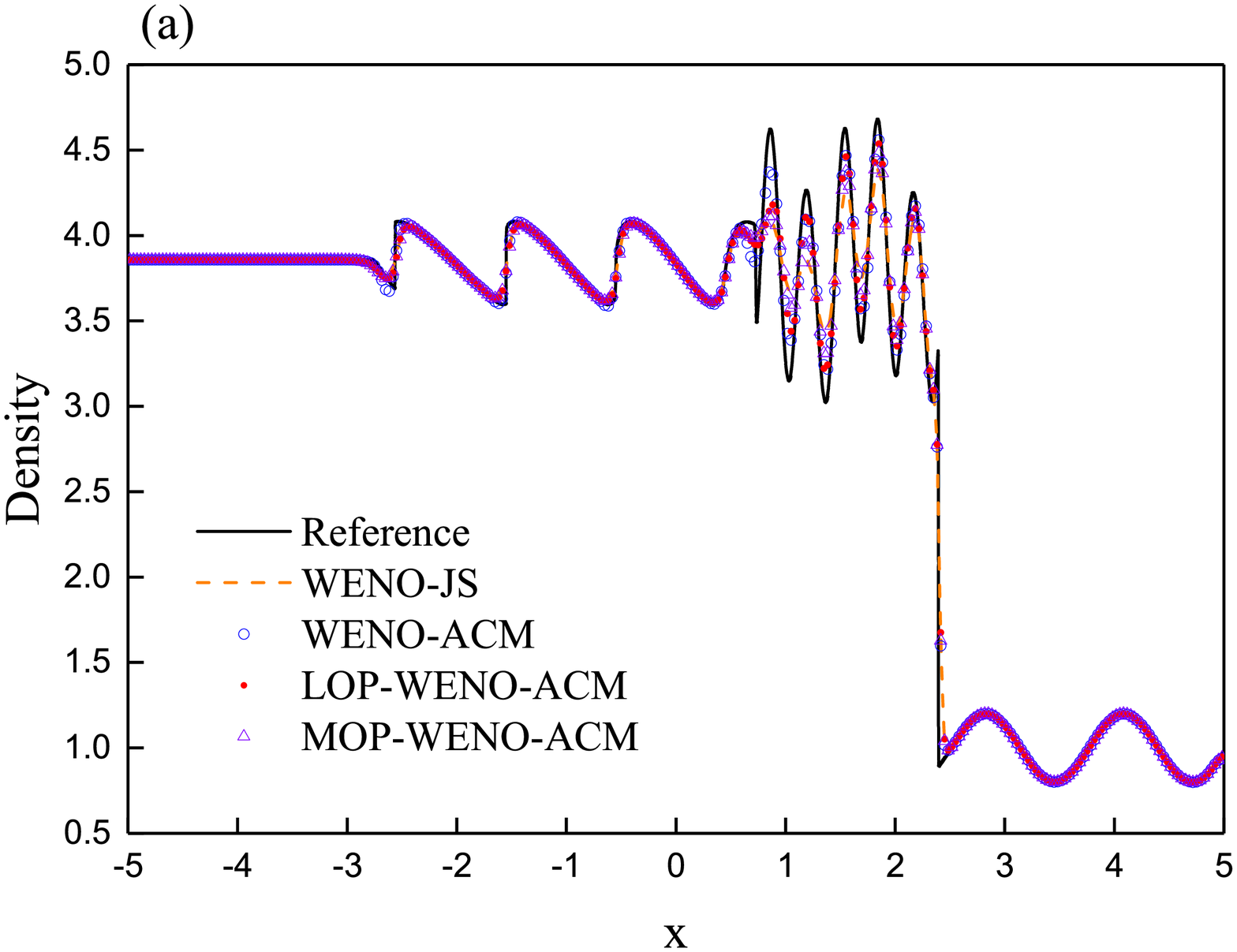}
  \includegraphics[height=0.37\textwidth]{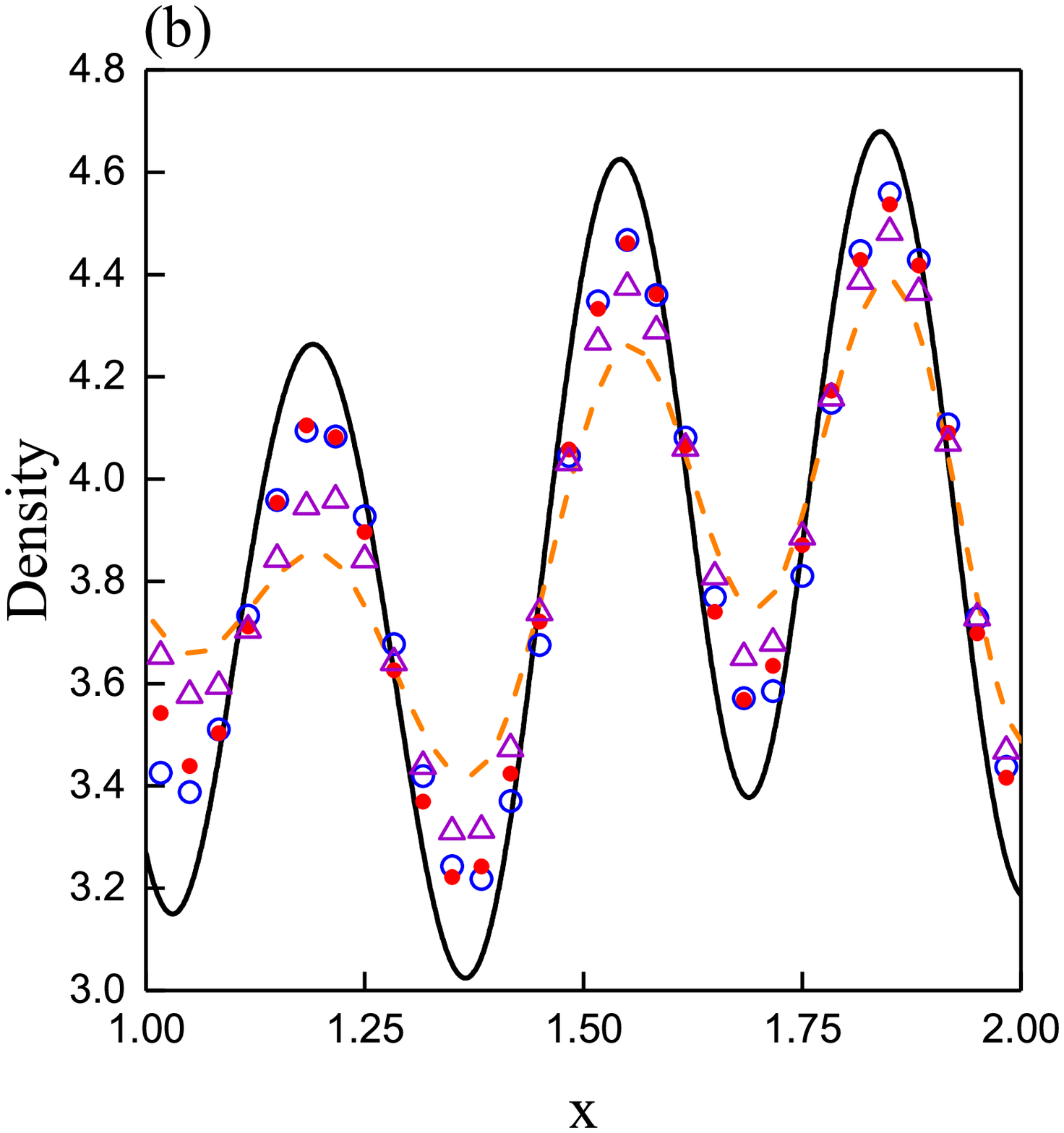} 
     \caption{Results of LOP-/MOP-/WENO-ACM and WENO-JS on solving the Shu-Osher problem.}
     \label{fig:Shu-Osher:WENO-ACM}
\end{figure}

We compute this problem with a uniform cell number of $N = 300$ by 
setting the CFL number to be 0.1. The solutions of density are given 
in Fig. \ref{fig:Shu-Osher:WENO-M} to Fig. 
\ref{fig:Shu-Osher:WENO-ACM} where the reference solution is 
computed by employing WENO-JS with $N = 10000$. For comparison 
purpose, we also present the solutions of the associated MOP-WENO-X 
schemes proposed in our previous work \cite{MOP-WENO-X} and that of 
WENO-JS. Obviously, WENO-JS provides the lowest resolution. 
Unfortunately, the resolutions of the MOP-WENO-X schemes are much 
lower than those of the WENO-X schemes. However, the LOP-WENO-X 
schemes can get comparable resolutions with those of the WENO-X 
schemes.

\subsubsection{Titarev-Toro problem}
\begin{example}
This problem was presented by Titarev and Toro \cite{Titarev-Toro-1,
Titarev-Toro-2,Titarev-Toro-3} and it is a more severe version of 
the Shu-Osher problem. The computational domain of $[-5, 5]$ is 
initialized by
\begin{equation}
\big( \rho, u, p \big)(x, 0) =\left\{
\begin{array}{ll}
(1.515695, 0.5233346, 1.80500), & x \in [-5.0, -4.5], \\
(1.0 + 0.1\sin(20\pi x), 0, 1), & x \in [-4.5, 5.0].
\end{array}\right.
 \label{initial:Titarev-Toro}
\end{equation} 
Also, the transmissive boundary conditions are used at $x = \pm 5$, 
while the output time is set to be $t = 5.0$ in this case.
\label{ex:Titarev-Toro}
\end{example}

\begin{figure}[!ht]
\centering
  \includegraphics[height=0.37\textwidth]{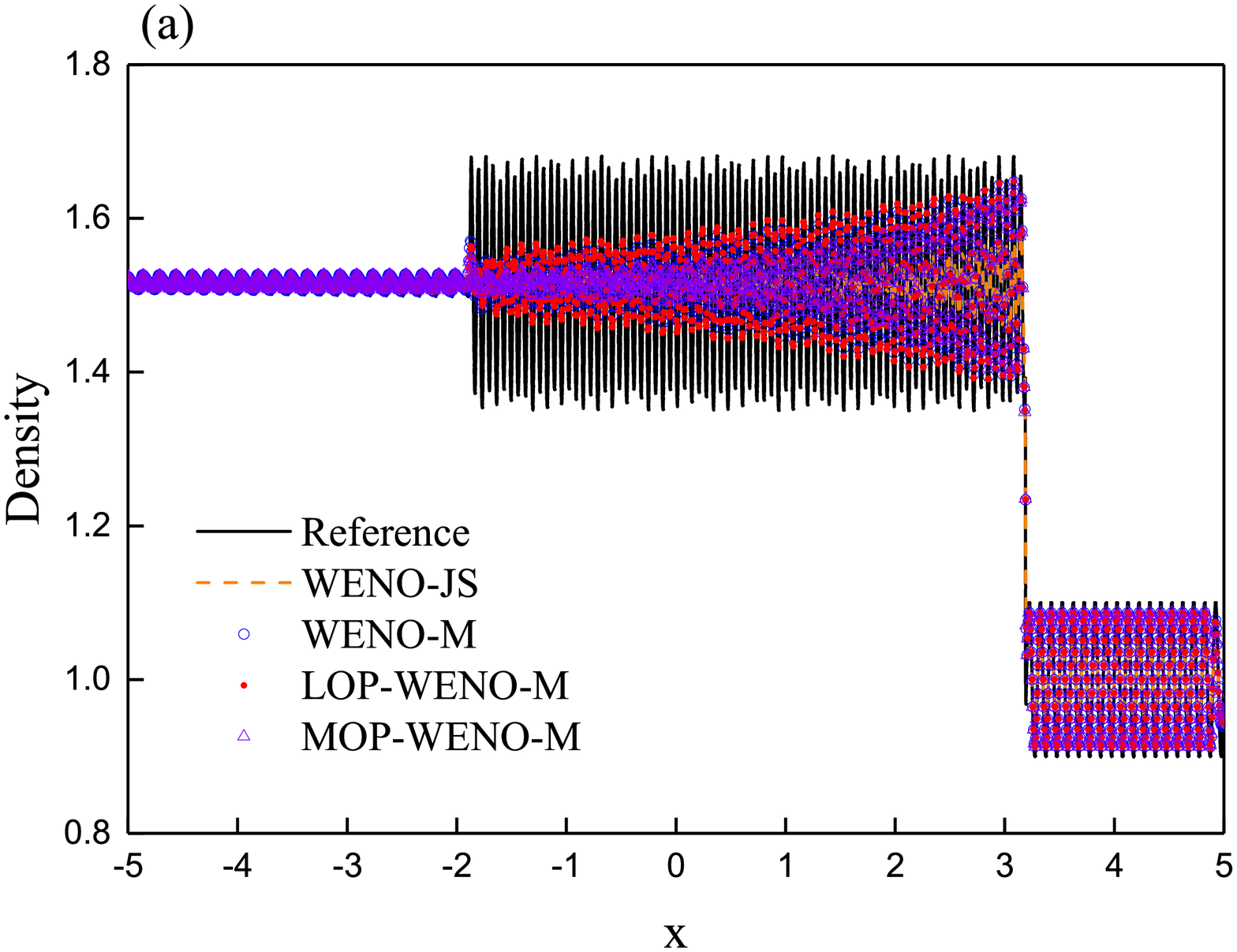}
  \includegraphics[height=0.37\textwidth]{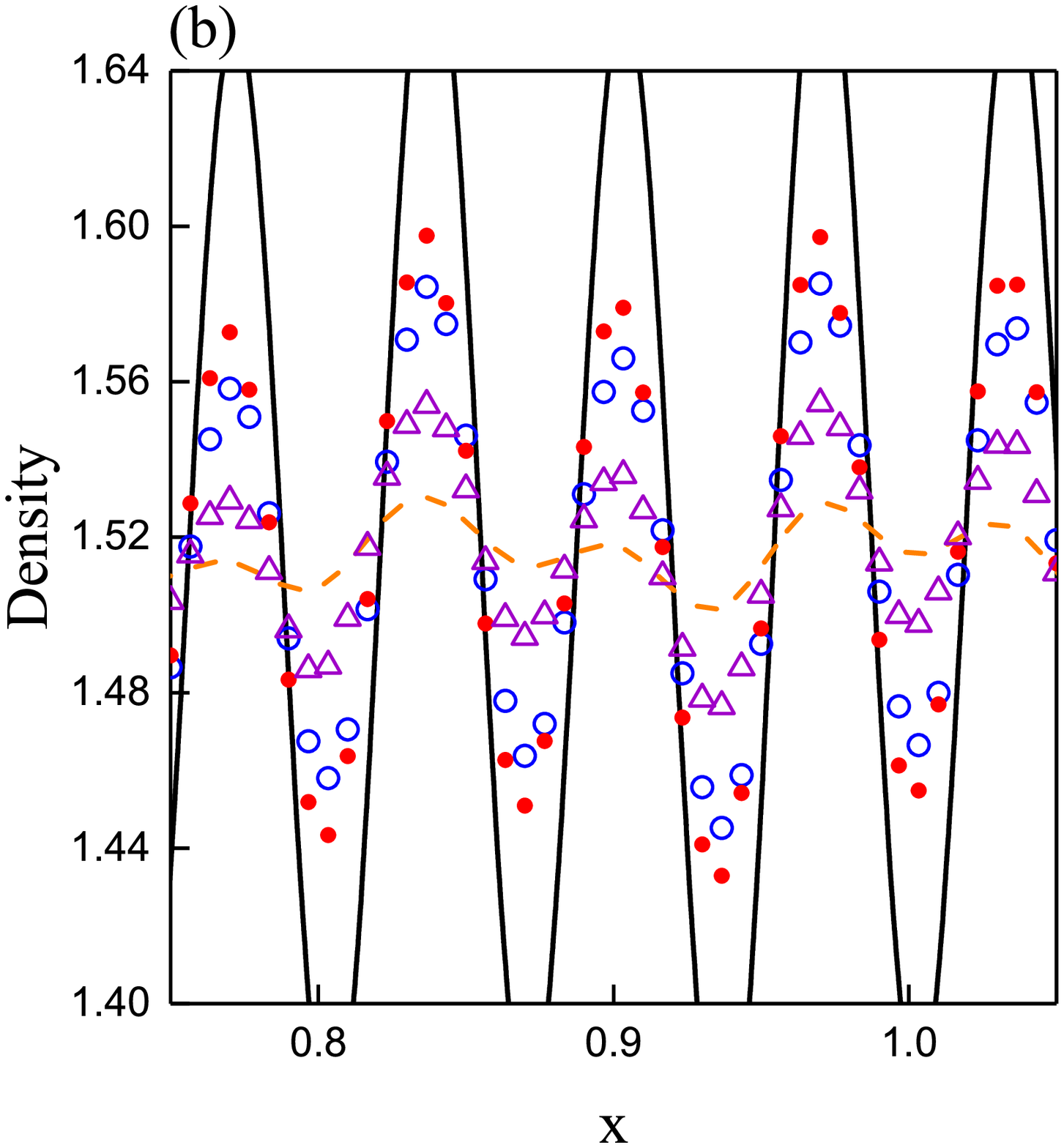} 
     \caption{Results of LOP-/MOP-/WENO-M and WENO-JS on solving the Titarev-Toro problem.}
     \label{fig:Titarev-Toro:WENO-M}
\end{figure}

\begin{figure}[!ht]
\centering
  \includegraphics[height=0.37\textwidth]{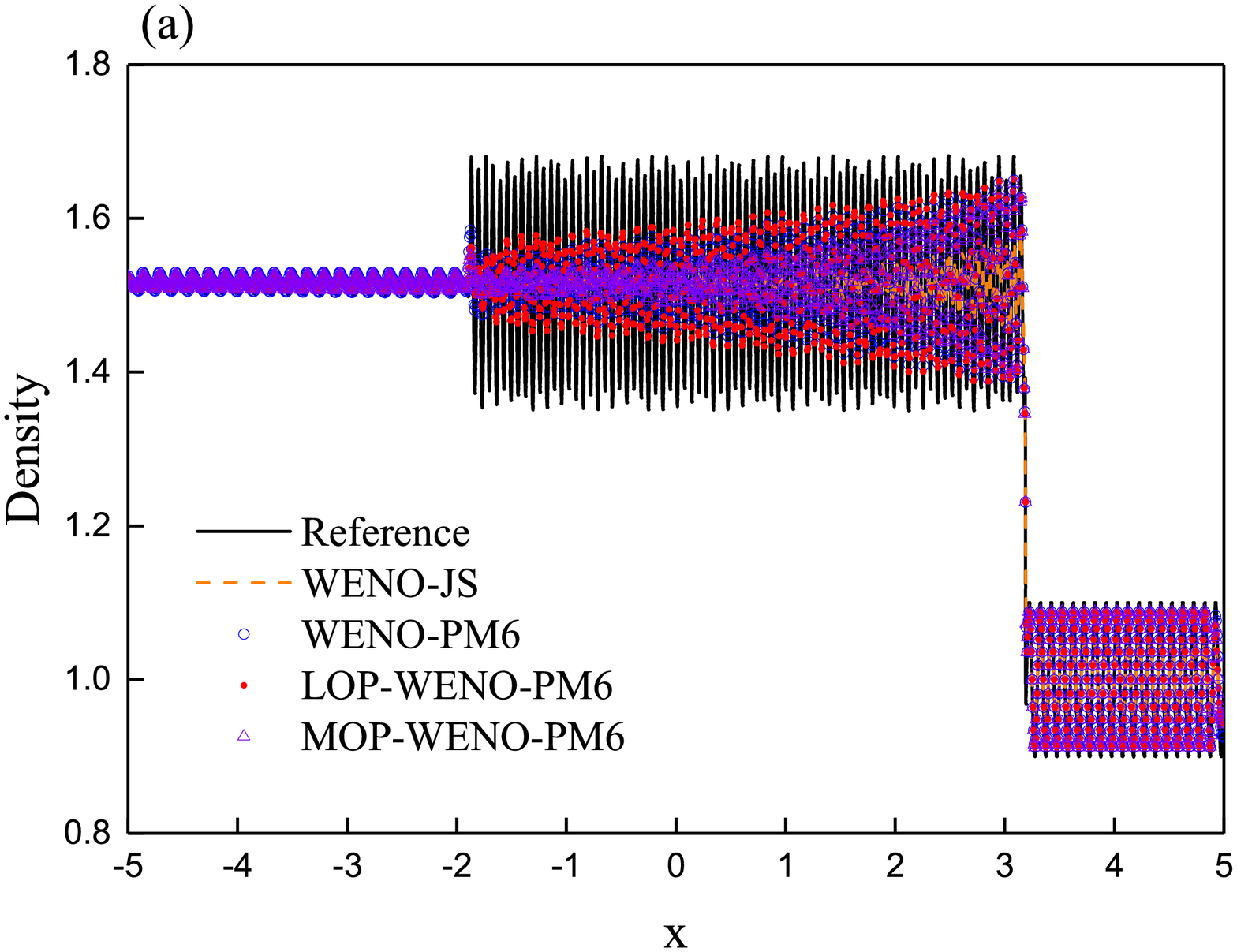}
  \includegraphics[height=0.37\textwidth]{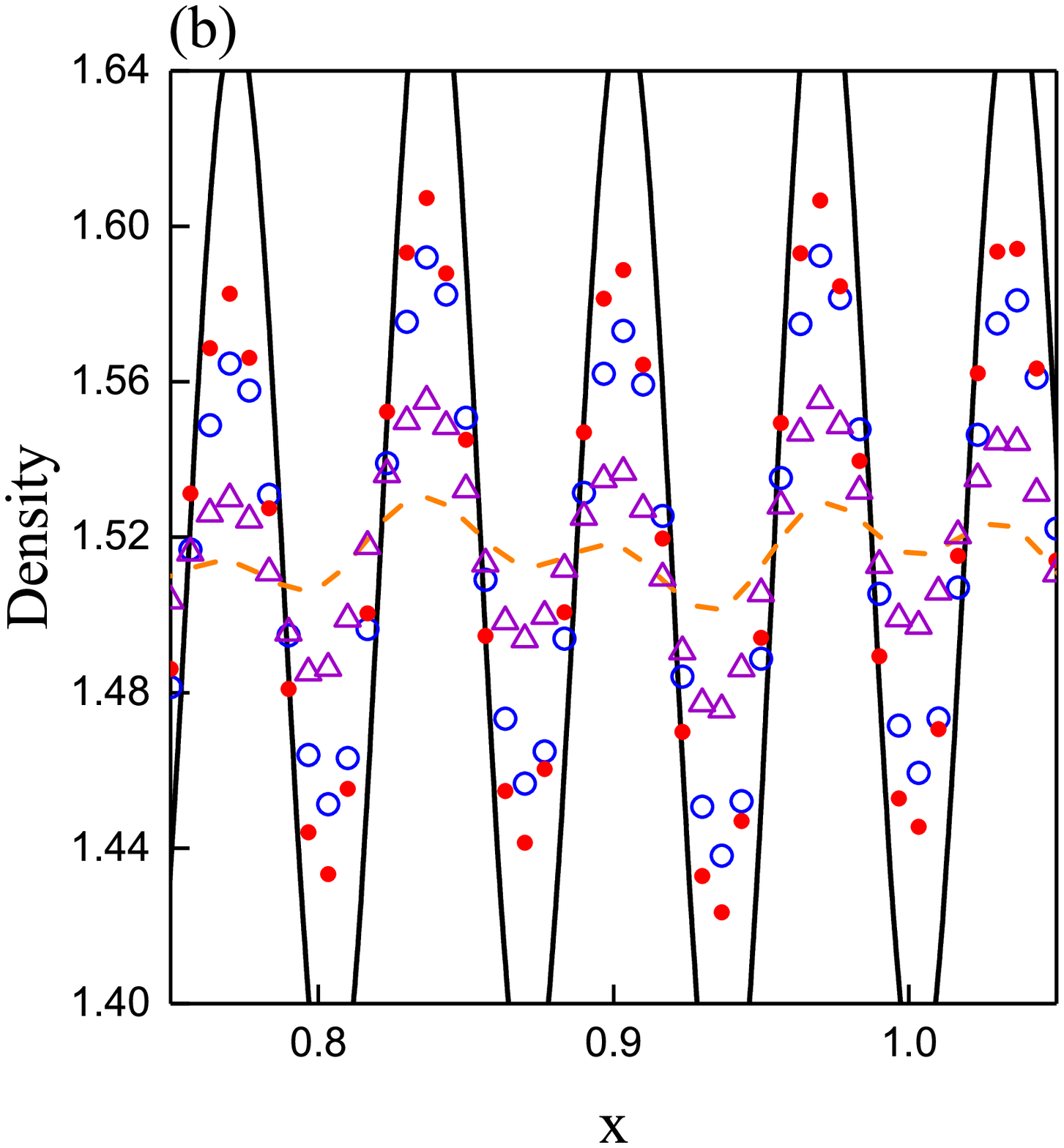} 
     \caption{Results of LOP-/MOP-/WENO-PM6 and WENO-JS on solving the Titarev-Toro problem.}
     \label{fig:Titarev-Toro:WENO-PM6}
\end{figure}

\begin{figure}[!ht]
\centering
  \includegraphics[height=0.37\textwidth]{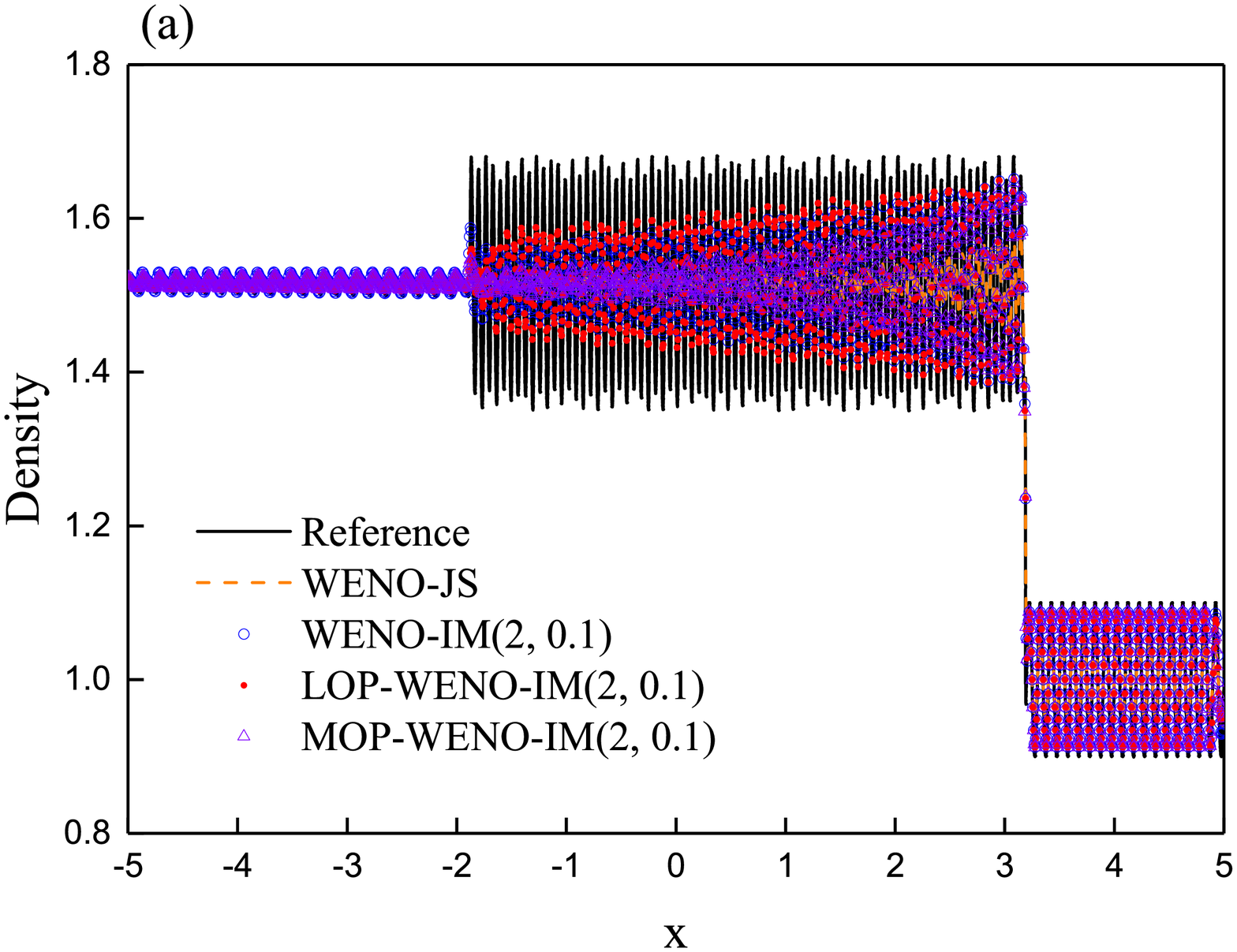}
  \includegraphics[height=0.37\textwidth]{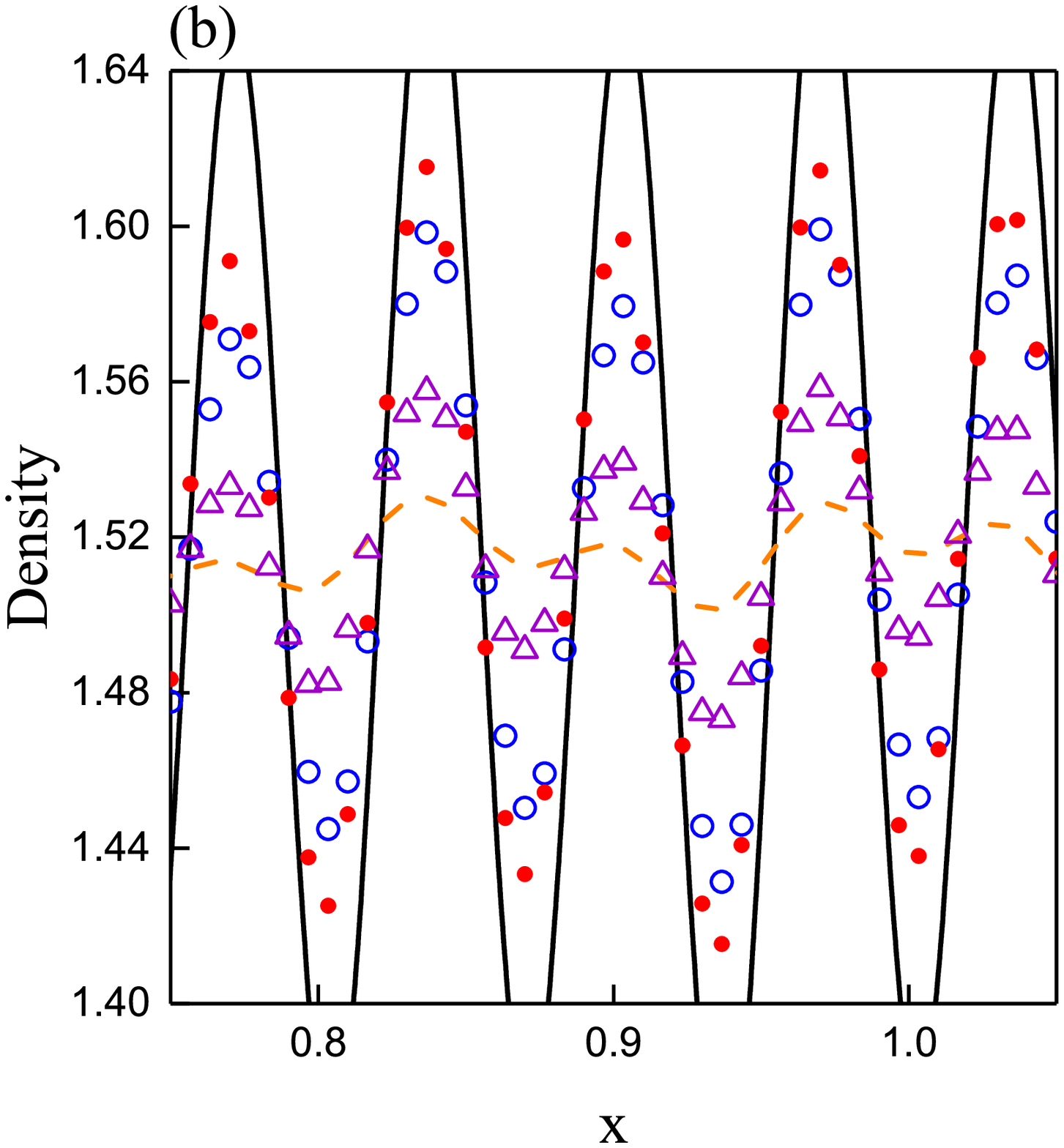}  
     \caption{Results of LOP-/MOP-/WENO-IM(2, 0.1) and WENO-JS on solving the Titarev-Toro problem.}
     \label{fig:Titarev-Toro:WENO-IM}
\end{figure}

\begin{figure}[!ht]
\centering
  \includegraphics[height=0.37\textwidth]{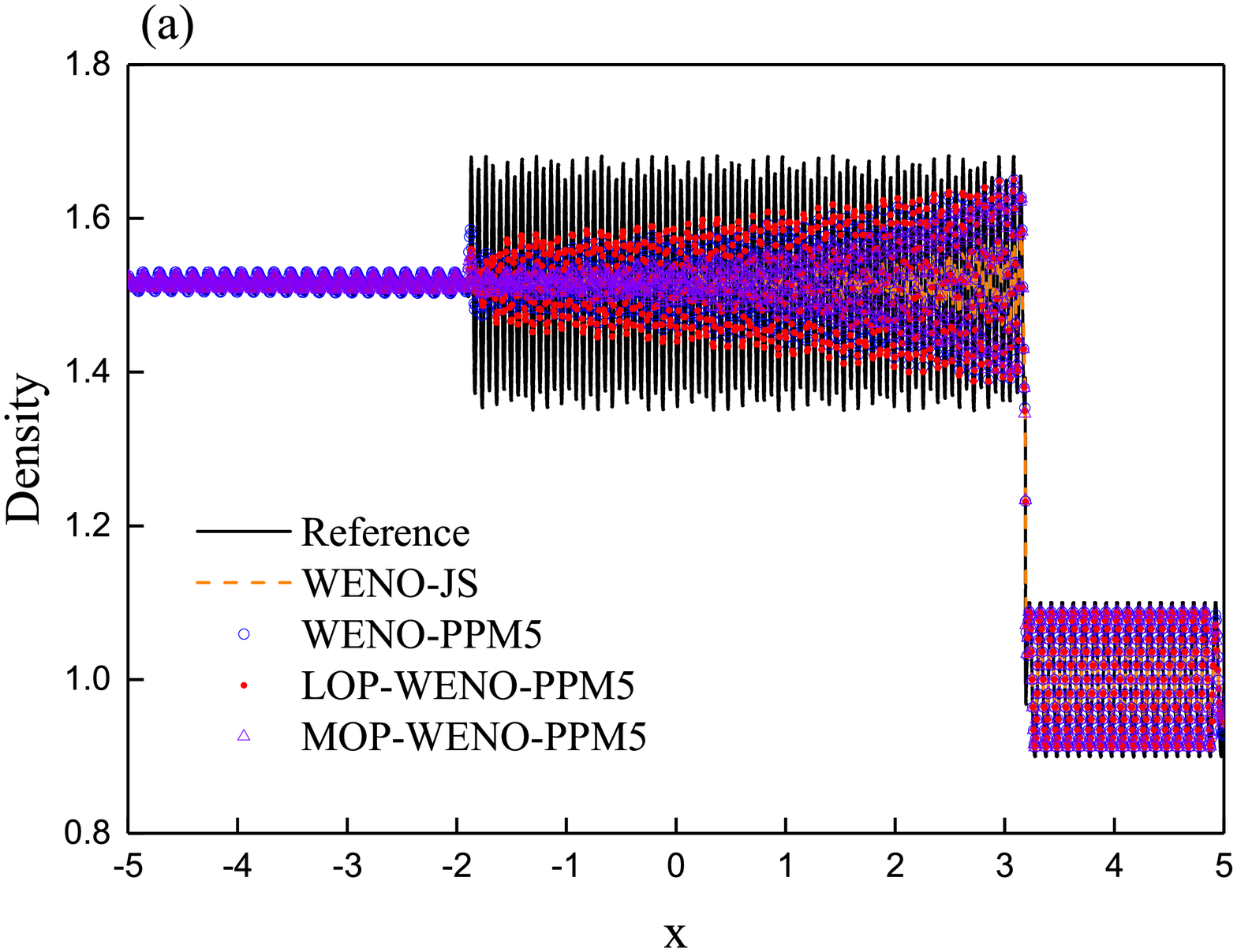}
  \includegraphics[height=0.37\textwidth]{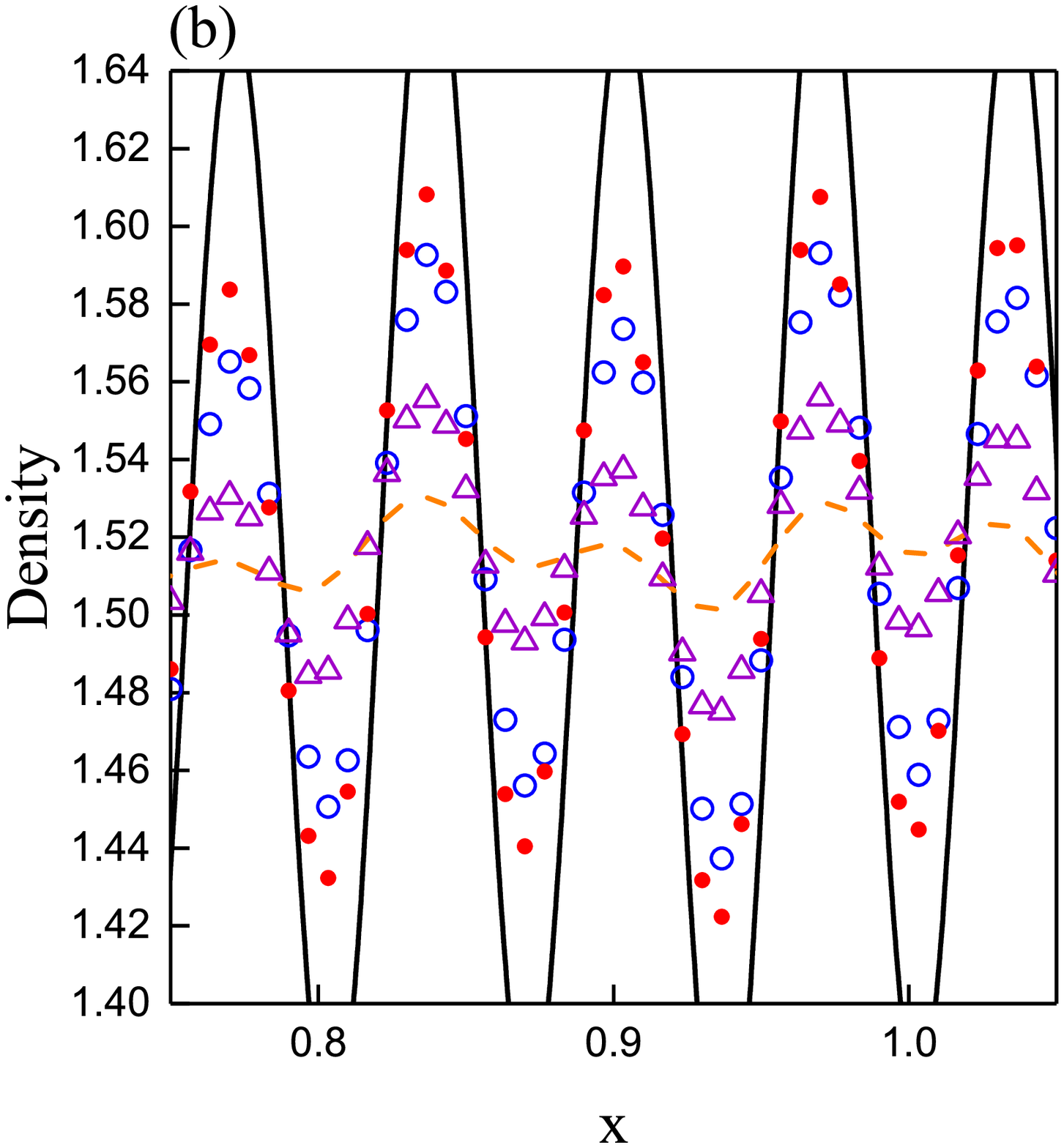}  
     \caption{Results of LOP-/MOP-/WENO-PPM5 and WENO-JS on solving the Titarev-Toro problem.}
     \label{fig:Titarev-Toro:WENO-PPM5}
\end{figure}

\begin{figure}[!ht]
\centering
  \includegraphics[height=0.37\textwidth]{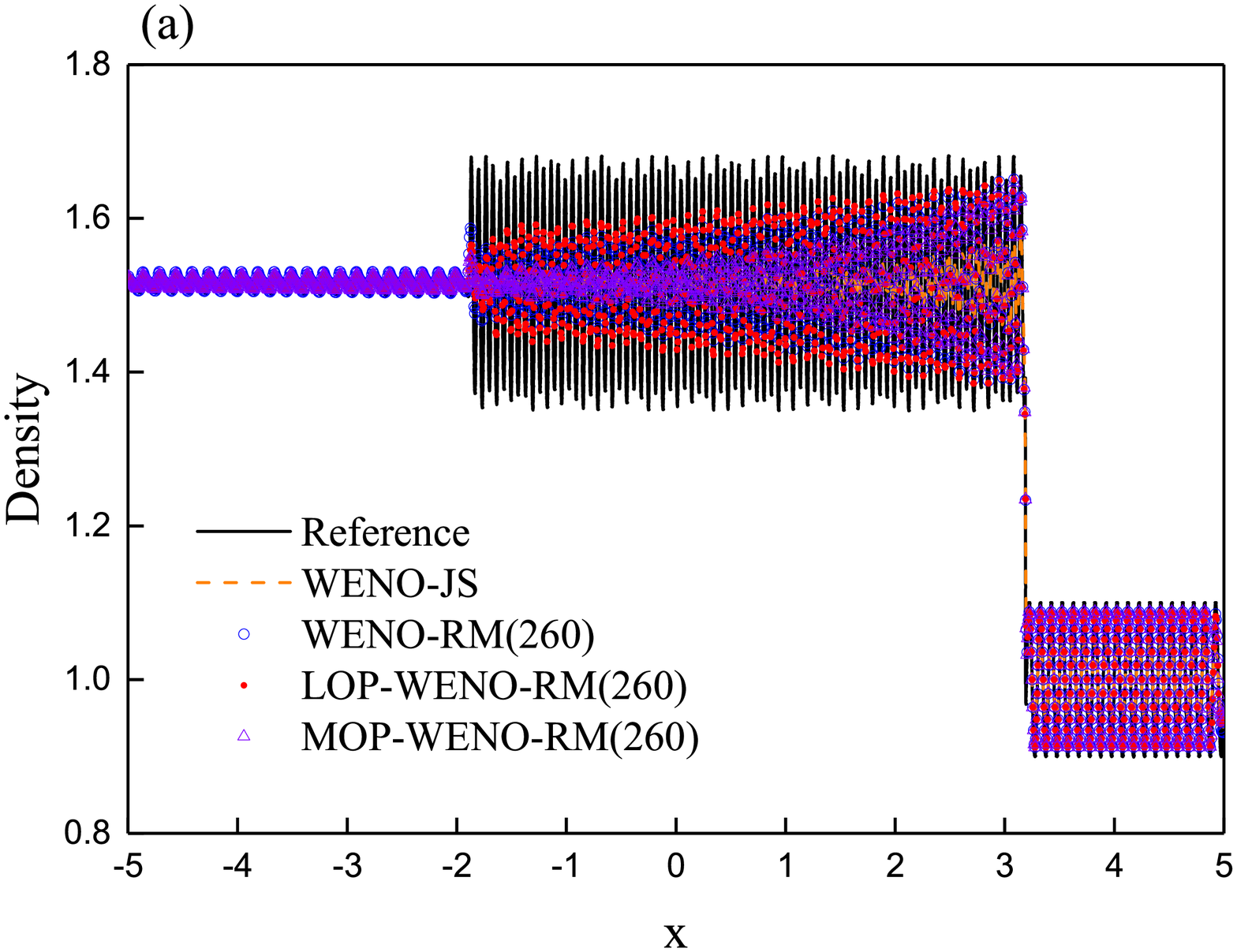}
  \includegraphics[height=0.37\textwidth]{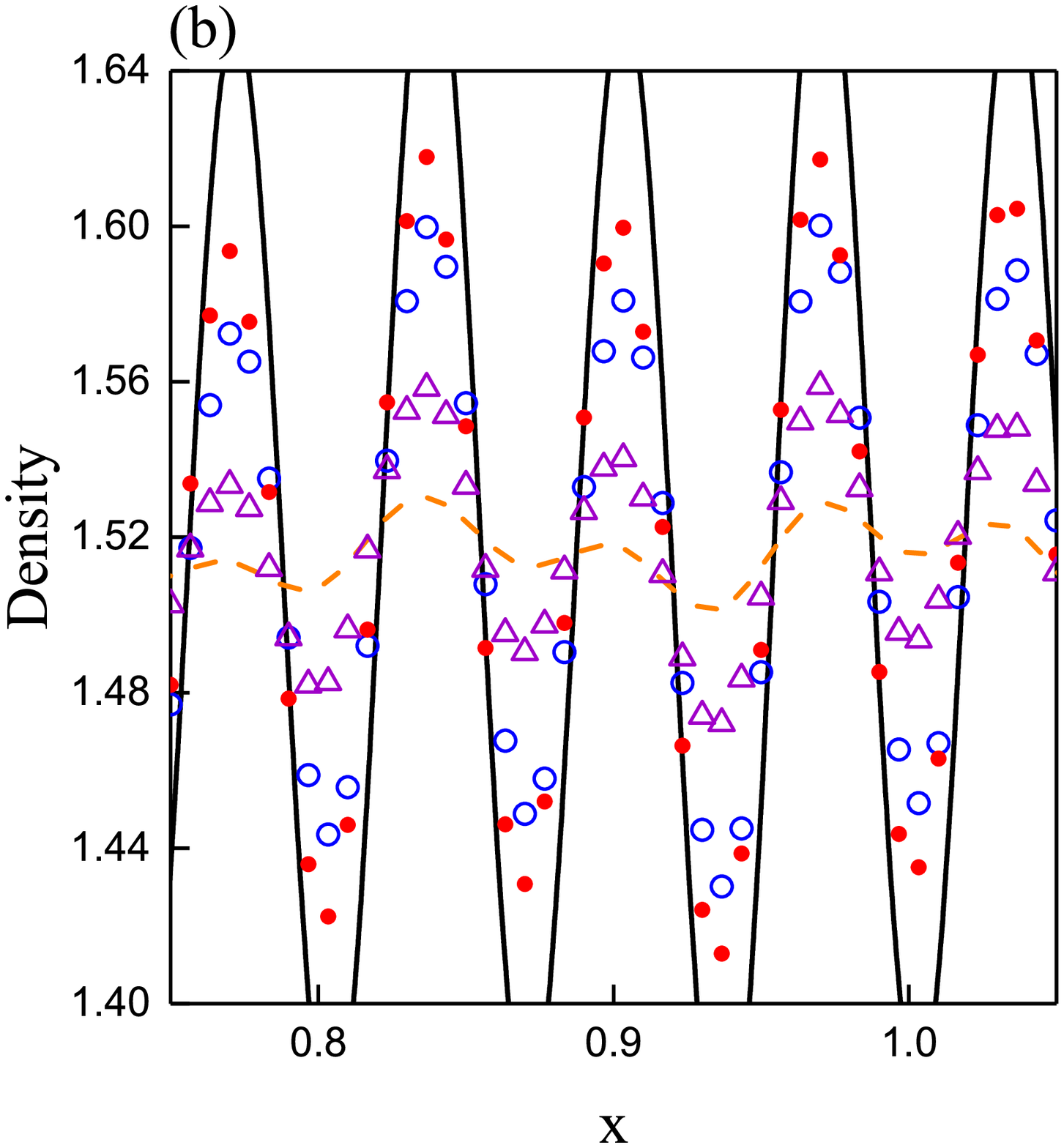} 
     \caption{Results of LOP-/MOP-/WENO-RM(260) and WENO-JS on solving the Titarev-Toro problem.}
     \label{fig:Titarev-Toro:WENO-RM260}
\end{figure}

\begin{figure}[!ht]
\centering
  \includegraphics[height=0.37\textwidth]{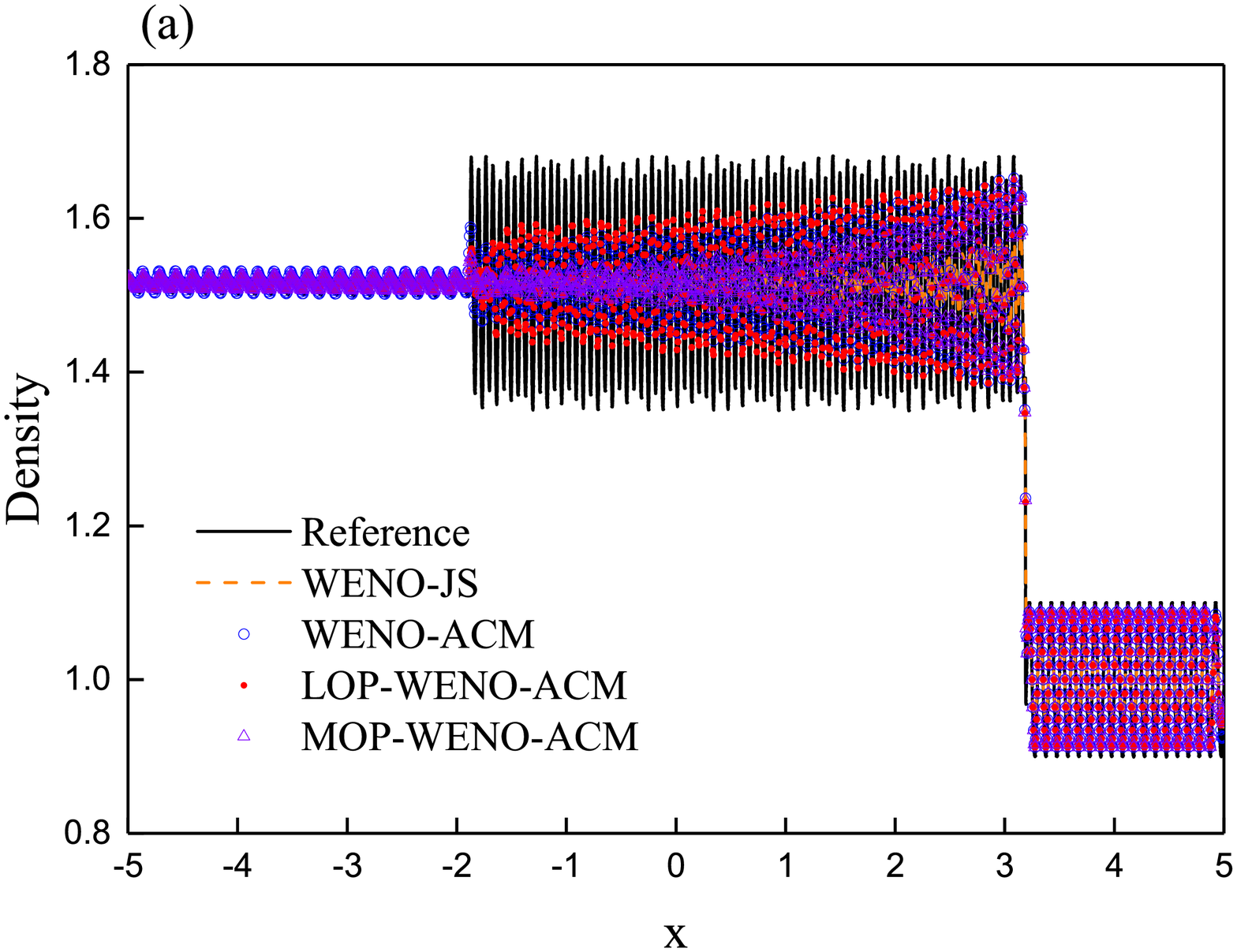}
  \includegraphics[height=0.37\textwidth]{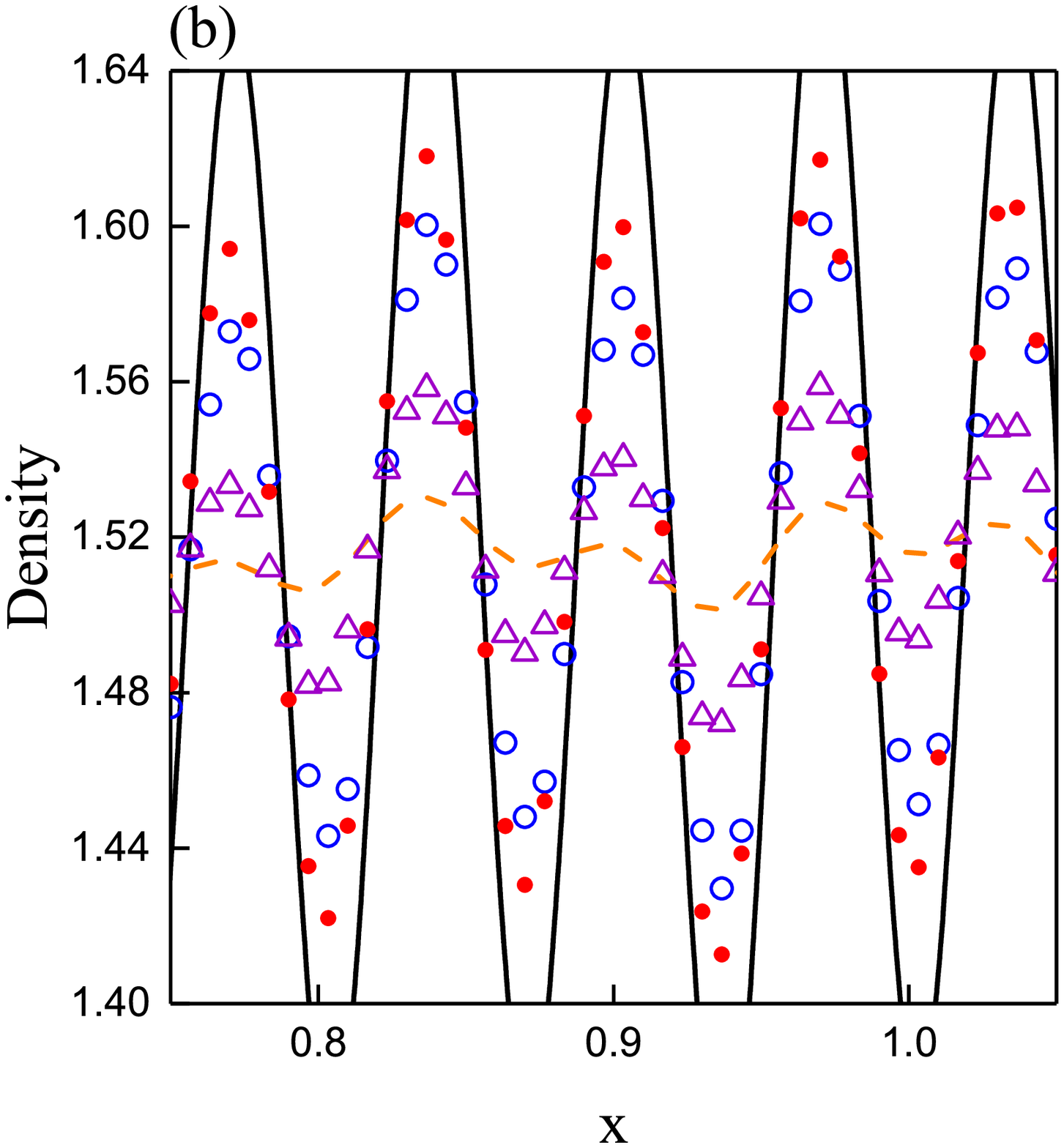} 
     \caption{Results of LOP-/MOP-/WENO-ACM and WENO-JS on solving the Titarev-Toro problem.}
     \label{fig:Titarev-Toro:WENO-ACM}
\end{figure}

We compute this problem with a uniform cell number of $N = 1500$ by 
setting the CFL number to be 0.4. The solutions of density are given 
in Fig. \ref{fig:Titarev-Toro:WENO-M} to Fig. 
\ref{fig:Titarev-Toro:WENO-ACM} where the reference solution is 
computed by employing WENO-JS with $N = 10000$. Again, for 
comparison purpose, we show the solutions of the associated 
MOP-WENO-X schemes and that of WENO-JS. Not surprisingly, WENO-JS 
provides the lowest resolution and the resolutions of the MOP-WENO-X 
schemes are much lower than those of the WENO-X schemes. 
Particularly, the resolutions of the LOP-WENO-X schemes are 
significantly higher than those of the WENO-X schemes. This is a 
remarkable competitive advantage of the LOP-WENO-X schemes compared 
to the MOP-WENO-X schemes.

\subsection{2D Euler equations}
In this subsection, we consider the following two-dimensional Euler 
systems of compressible gas dynamics
\begin{equation}
\dfrac{\partial \mathbf{U}}{\partial t}    + 
\dfrac{\partial \mathbf{F(U)}}{\partial x} +
\dfrac{\partial \mathbf{G(U)}}{\partial y} = \mathbf{0},
\label{2DEulerEquations}
\end{equation}
with
\begin{equation*}
\begin{array}{l}
\mathbf{U} = \Big( \rho, \rho u, \rho v, E \Big)^{\mathrm{T}}, \\
\mathbf{F(U)} = \Big(\rho u, \rho u^{2} + p, \rho uv, u(E+p) 
\Big)^{\mathrm{T}}, \\
\mathbf{G(U)} = \Big(\rho v, \rho vu, \rho v^{2} + p, v(E+p) 
\Big)^{\mathrm{T}},
\end{array}
\end{equation*}
where $\rho, u, v, p$ and $E$ are the density, components of 
velocity in the $x$ and $y$ coordinate directions, pressure and 
total energy, respectively. The relation of pressure $p$ and total 
energy for ideal gases is defined by
\begin{equation*}
p = (\gamma - 1)\Big( E - \dfrac{1}{2}\rho(u^{2} + v^{2})\Big), \quad
\gamma = 1.4.
\end{equation*}
Two commonly used classes of finite volume WENO schemes in 
two-dimensional Cartesian meshes were studied carefully by Zhang et 
al. \cite{FVMaccuracyProofs03}. Here, we implement the one denoted 
as class A for our calculations. Now, we examine the performances of 
the considered WENO shcemes on solving the following three benchmark 
tests.

\subsubsection{Accuracy test 1}
\begin{example}
\rm{We use this density wave propagation problem 
\cite{AccuracyTest-Euler2D} to test the convergence orders of the 
considered WENO schemes. The initial condition on the computational 
domain $[-1.0, 1.0]\times[-1.0, 1.0]$ is given by}
\label{ex:AccuracyTest:case1}
\end{example}
\begin{equation}
(\rho, u, v, p)(x, y, 0) = \Big( 1.0 + 0.2 \sin\big( \pi(x + y) \big)
, 0.7, 0.3, 1.0\Big).
\end{equation}
Here, the $L_{1}$ and $L_{\infty}$ errors are computed by
\begin{equation*}
\displaystyle
\begin{aligned}
&L_{1} = h_{x}h_{y} \cdot \displaystyle\sum\nolimits_{j=1}^{N_{y}}
\sum\nolimits_{i=1}^{N_{x}} \Big\lvert \rho_{i,j}^{\mathrm{exact}} - 
(\rho_{h})_{i,j} \Big\rvert, \quad
&L_{\infty} = \displaystyle\max_{\substack{1\leq i\leq N_{x}\\1\leq j
\leq N_{y}}} \Big\lvert \rho_{i,j}^{\mathrm{exact}} -(\rho_{h})_{i,j}
\Big\rvert,
\end{aligned}
\end{equation*}
where $N_{x}, N_{y}$ is the number of cells in $x-$ and $y-$ 
direction, and $h_{x}, h_{y}$ is the associated uniform spatial step 
size and we set $h = h_{x} = h_{y}$ in all calculations of this 
paper. $(\rho_{h})_{i,j}$ is the numerical solution of the density 
and $\rho_{i,j}^{\mathrm{exact}}$ is its exact solution. We can 
easily check that the exact solution is $\rho(x,y,t) = 1.0 + 
0.2\sin\Big( \pi\big(x + y - (u + v)t \big)\Big)$, $u(x, y ,t) = 0.7$
, $v(x, y, t) = 0.3$, and $p(x, y, t) = 1.0$.

The computational time is advanced until $t = 2.0$. The periodic 
boundary condition is used and the CFL number is taken to be 
$h^{2/3}$ so that the error for the overall scheme is a measure of 
the spatial convergence only.

The numerical errors and corresponding convergence orders of 
accuracy for the density $\rho$ are shown in Table 
\ref{table:AccuracyTest:case1}. Again, for comparison purpose, we 
also present the results computed by the WENO5-ILW scheme. We can 
see that all the considered WENO schemes can achieve the disgned 
order of accuracy, while the error magnitude is larger for the 
WENO-JS scheme than for the other schemes. Moreover, as expected, 
the numerical errors with respect to all grid numbers of the 
LOP-WENO-X schemes are almost the same to those of the WENO-X 
schemes. It should be noted that, in terms of accuracy, the 
LOP-WENO-PM6, LOP-WENO-PPM5, LOP-WENO-RM(260) and LOP-WENO-ACM 
schemes provide the numerical errors equivalent to that of the 
WENO5-ILW scheme, but of course this is not always the case and we 
will show it in the next test.

\begin{table}[ht]
\begin{myFontSize}
\centering
\caption{Numerical errors and convergence orders of accuracy for the 
density $\rho$ on Example \ref{ex:AccuracyTest:case1} at $t = 2.0$.}
\label{table:AccuracyTest:case1}
\begin{tabular*}{\hsize}
{@{}@{\extracolsep{\fill}}cllllllll@{}}
\hline
\space    &\multicolumn{4}{l}{\cellcolor{gray!35}{WENO5-ILW}}  
          &\multicolumn{4}{l}{\cellcolor{gray!35}{WENO-JS}}  \\
\cline{2-5}  \cline{6-9}
$N_{x}\times N_{y}$   & $L_{1}$ error      & $L_{1}$ order 
                      & $L_{\infty}$ error & $L_{\infty}$ order
		  			  & $L_{1}$ error      & $L_{1}$ order 
      				  & $L_{\infty}$ error & $L_{\infty}$ order \\
\Xhline{0.65pt}
40 $\times$ 40        & 2.05111E-05        & -
                      & 8.06075E-06        & -
                      & 1.44379E-04        & -
                      & 6.11870E-05        & -                  \\
60 $\times$ 60        & 2.71152E-06        & 4.9905
                      & 1.06517E-06        & 4.9915
                      & 1.90416E-05        & 4.9963
                      & 8.53414E-06        & 4.8583             \\
80 $\times$ 80        & 6.44325E-07        & 4.9953
                      & 2.53073E-07        & 4.9958
                      & 4.51609E-06        & 5.0020
                      & 2.03906E-06        & 4.9763             \\ 
100 $\times$ 100      & 2.11264E-07        & 4.9972
                      & 8.29736E-08        & 4.9975
                      & 1.47974E-06        & 5.0003
                      & 6.75061E-07        & 4.9539             \\ 
\hline
\space    &\multicolumn{4}{l}{\cellcolor{gray!35}{WENO-M}}
          &\multicolumn{4}{l}{\cellcolor{gray!35}{LOP-WENO-M}}\\
\cline{2-5}  \cline{6-9}
$N_{x}\times N_{y}$   & $L_{1}$ error      & $L_{1}$ order 
                      & $L_{\infty}$ error & $L_{\infty}$ order
		  			  & $L_{1}$ error      & $L_{1}$ order 
      				  & $L_{\infty}$ error & $L_{\infty}$ order \\
\Xhline{0.65pt}
40 $\times$ 40        & 2.05584E-05        & -
                      & 8.07114E-06        & -
                      & 2.05584E-05        & -
                      & 8.07114E-06        & -                  \\
60 $\times$ 60        & 2.71274E-06        & 4.9950
                      & 1.06546E-06        & 4.9940
                      & 2.71274E-06        & 4.9950
                      & 1.06546E-06        & 4.9940             \\
80 $\times$ 80        & 6.44416E-07        & 4.9964
                      & 2.53097E-07        & 4.9965
                      & 6.44416E-07        & 4.9964
                      & 2.53097E-07        & 4.9965             \\ 
100 $\times$ 100      & 2.11276E-07        & 4.9976
                      & 8.29769E-08        & 4.9977
                      & 2.11276E-07        & 4.9976
                      & 8.29769E-08        & 4.9977             \\
\hline
\space    &\multicolumn{4}{l}{\cellcolor{gray!35}{WENO-PM6}} 
          &\multicolumn{4}{l}{\cellcolor{gray!35}{LOP-WENO-PM6}}\\
\cline{2-5}  \cline{6-9}
$N_{x}\times N_{y}$   & $L_{1}$ error      & $L_{1}$ order 
                      & $L_{\infty}$ error & $L_{\infty}$ order
		  			  & $L_{1}$ error      & $L_{1}$ order 
      				  & $L_{\infty}$ error & $L_{\infty}$ order \\
\Xhline{0.65pt}
40 $\times$ 40        & 2.05111E-05        & -
                      & 8.06076E-06        & -
                      & 2.05111E-05        & -
                      & 8.06076E-06        & -                 \\
60 $\times$ 60        & 2.71152E-06        & 4.9905
                      & 1.06517E-06        & 4.9915
                      & 2.71152E-06        & 4.9905
                      & 1.06517E-06        & 4.9915            \\
80 $\times$ 80        & 6.44325E-07        & 4.9953
                      & 2.53073E-07        & 4.9958
                      & 6.44325E-07        & 4.9953
                      & 2.53073E-07        & 4.9958             \\ 
100 $\times$ 100      & 2.11264E-07        & 4.9972
                      & 8.29736E-08        & 4.9975
                      & 2.11264E-07        & 4.9972
                      & 8.29736E-08        & 4.9975             \\
\hline
\space    &\multicolumn{4}{l}{\cellcolor{gray!35}{WENO-IM(2, 0.1)}} 
    &\multicolumn{4}{l}{\cellcolor{gray!35}{LOP-WENO-IM(2, 0.1)}}\\
\cline{2-5}  \cline{6-9}
$N_{x}\times N_{y}$   & $L_{1}$ error      & $L_{1}$ order 
                      & $L_{\infty}$ error & $L_{\infty}$ order
		  			  & $L_{1}$ error      & $L_{1}$ order 
      				  & $L_{\infty}$ error & $L_{\infty}$ order \\
\Xhline{0.65pt}
40 $\times$ 40        & 2.05159E-05        & -
                      & 8.06179E-06        & -
                      & 2.05159E-05        & -
                      & 8.06179E-06        & -                  \\
60 $\times$ 60        & 2.71164E-06        & 4.9909
                      & 1.06520E-06        & 4.9917
                      & 2.71164E-06        & 4.9909
                      & 1.06520E-06        & 4.9917             \\
80 $\times$ 80        & 6.44334E-07        & 4.9954
                      & 2.53076E-07        & 4.9959
                      & 6.44334E-07        & 4.9954
                      & 2.53076E-07        & 4.9959             \\ 
100 $\times$ 100      & 2.11265E-07        & 4.9972
                      & 8.29739E-08        & 4.9975
                      & 2.11265E-07        & 4.9972
                      & 8.29739E-08        & 4.9975             \\
\hline
\space    &\multicolumn{4}{l}{\cellcolor{gray!35}{WENO-PPM5}} 
          &\multicolumn{4}{l}{\cellcolor{gray!35}{LOP-WENO-PPM5}}\\
\cline{2-5}  \cline{6-9}
$N_{x}\times N_{y}$   & $L_{1}$ error      & $L_{1}$ order 
                      & $L_{\infty}$ error & $L_{\infty}$ order
		  			  & $L_{1}$ error      & $L_{1}$ order 
      				  & $L_{\infty}$ error & $L_{\infty}$ order \\
\Xhline{0.65pt}
40 $\times$ 40        & 2.05111E-05        & -
                      & 8.06083E-06        & -
                      & 2.05111E-05        & -
                      & 8.06083E-06        & -                  \\
60 $\times$ 60        & 2.71152E-06        & 4.9905
                      & 1.06517E-06        & 4.9915
                      & 2.71152E-06        & 4.9905
                      & 1.06517E-06        & 4.9915             \\
80 $\times$ 80        & 6.44325E-07        & 4.9953
                      & 2.53073E-07        & 4.9958
                      & 6.44325E-07        & 4.9953
                      & 2.53073E-07        & 4.9958             \\ 
100 $\times$ 100      & 2.11264E-07        & 4.9972
                      & 8.29736E-08        & 4.9975
                      & 2.11264E-07        & 4.9972
                      & 8.29736E-08        & 4.9975             \\
\hline
\space    &\multicolumn{4}{l}{\cellcolor{gray!35}{WENO-RM(260)}} 
      &\multicolumn{4}{l}{\cellcolor{gray!35}{LOP-WENO-RM(260)}}\\
\cline{2-5}  \cline{6-9}
$N_{x}\times N_{y}$   & $L_{1}$ error      & $L_{1}$ order 
                      & $L_{\infty}$ error & $L_{\infty}$ order
		  			  & $L_{1}$ error      & $L_{1}$ order 
      				  & $L_{\infty}$ error & $L_{\infty}$ order \\
\Xhline{0.65pt}
40 $\times$ 40        & 2.05111E-05        & -
                      & 8.06075E-06        & -
                      & 2.05111E-05        & -
                      & 8.06075E-06        & -                  \\
60 $\times$ 60        & 2.71152E-06        & 4.9905
                      & 1.06517E-06        & 4.9915
                      & 2.71152E-06        & 4.9905
                      & 1.06517E-06        & 4.9915             \\
80 $\times$ 80        & 6.44325E-07        & 4.9953
                      & 2.53073E-07        & 4.9958
                      & 6.44325E-07        & 4.9953
                      & 2.53073E-07        & 4.9958             \\ 
100 $\times$ 100      & 2.11264E-07        & 4.9972
                      & 8.29736E-08        & 4.9975
                      & 2.11264E-07        & 4.9972
                      & 8.29736E-08        & 4.9975             \\
\hline
\space    &\multicolumn{4}{l}{\cellcolor{gray!35}{WENO-ACM}} 
          &\multicolumn{4}{l}{\cellcolor{gray!35}{LOP-WENO-ACM}}\\
\cline{2-5}  \cline{6-9}
$N_{x}\times N_{y}$   & $L_{1}$ error      & $L_{1}$ order 
                      & $L_{\infty}$ error & $L_{\infty}$ order
		  			  & $L_{1}$ error      & $L_{1}$ order 
      				  & $L_{\infty}$ error & $L_{\infty}$ order \\
\Xhline{0.65pt}
40 $\times$ 40        & 2.05111E-05        & -
                      & 8.06075E-06        & -
                      & 2.05111E-05        & -
                      & 8.06075E-06        & -                  \\
60 $\times$ 60        & 2.71152E-06        & 4.9905
                      & 1.06517E-06        & 4.9915
                      & 2.71152E-06        & 4.9905
                      & 1.06517E-06        & 4.9915             \\
80 $\times$ 80        & 6.44325E-07        & 4.9953
                      & 2.53073E-07        & 4.9958
                      & 6.44325E-07        & 4.9953
                      & 2.53073E-07        & 4.9958             \\ 
100 $\times$ 100      & 2.11264E-07        & 4.9972
                      & 8.29736E-08        & 4.9975
                      & 2.11264E-07        & 4.9972
                      & 8.29736E-08        & 4.9975             \\
\hline
\end{tabular*}
\end{myFontSize}
\end{table}

\subsubsection{Accuracy test 2}
\begin{example}
\rm{Now we use a modified version of the density wave propagation 
problem \cite{AccuracyTest-Euler2D} to test the convergence orders 
of the considered WENO schemes. Here, the initial condition on the 
computational domain $[-1.0, 1.0]\times[-1.0, 1.0]$ is given by}
\label{ex:AccuracyTest:case2}
\end{example}
\begin{equation}
(\rho, u, v, p)(x, y, 0) = \Bigg( 1.0 + 0.2 \sin\bigg( \pi(x + y) - 
\dfrac{\sin(\pi(x + y))}{\pi} \bigg)
, 0.7, 0.3, 1.0\Bigg).
\end{equation}
Again, the computational time is advanced until $t = 2.0$. And also 
the periodic boundary condition is used and the CFL number is taken 
to be $h^{2/3}$. Trivially, the exact solution is $\rho(x,y,t) = 1.0 
+ 0.2\sin\Bigg( \pi\bigg(x + y - (u + v)t - \frac{\sin\Big(\pi \big(x
+ y - (u + v) t\big)\Big)}{\pi}\bigg)\Bigg)$, $u(x,y,t) = 0.7$, $v(x
,y,t) = 0.3$, and $p(x,y,t) = 1.0$.

The numerical errors and corresponding convergence orders of 
accuracy for the density $\rho$ are shown in Table 
\ref{table:AccuracyTest:case2}. It is noted that the 
$L_{\infty}$ convergence order of the WENO-JS scheme drops by 
nearly 2 orders that leads to an overall accuracy loss shown with 
the $L_{1}$ convergence order. However, it is evident that 
the other schemes can retain the optimal convergence orders even in 
the presence of critical points. Unsurprisingly, in terms of 
accuracy, the LOP-WENO-X schemes give equally accurate numerical 
solutions like those of their associated WENO-X schemes. We point 
out that, for this test, only the LOP-WENO-ACM/WENO-ACM scheme 
provides the numerical errors equivalent to that of the WENO5-ILW 
scheme.

\begin{table}[ht]
\begin{myFontSize}
\centering
\caption{Numerical errors and convergence orders of accuracy for the 
density $\rho$ on Example \ref{ex:AccuracyTest:case2} at $t = 2.0$.}
\label{table:AccuracyTest:case2}
\begin{tabular*}{\hsize}
{@{}@{\extracolsep{\fill}}cllllllll@{}}
\hline
\space    &\multicolumn{4}{l}{\cellcolor{gray!35}{WENO5-ILW}}  
          &\multicolumn{4}{l}{\cellcolor{gray!35}{WENO-JS}}  \\
\cline{2-5}  \cline{6-9}
$N_{x}\times N_{y}$   & $L_{1}$ error      & $L_{1}$ order 
                      & $L_{\infty}$ error & $L_{\infty}$ order
		  			  & $L_{1}$ error      & $L_{1}$ order 
      				  & $L_{\infty}$ error & $L_{\infty}$ order \\
\Xhline{0.65pt}
40 $\times$ 40        & 2.31214E-04        & -
                      & 1.58230E-04        & -
                      & 8.15797E-04        & -
                      & 5.48728E-04        & -                  \\
60 $\times$ 60        & 3.13106E-05        & 4.9311
                      & 2.18798E-05        & 4.8795
                      & 1.61432E-04        & 3.9956
                      & 1.32554E-04        & 3.5037             \\
80 $\times$ 80        & 7.48937E-06        & 4.9724
                      & 5.24972E-06        & 4.9617
                      & 4.67993E-05        & 4.3041
                      & 4.84021E-05        & 3.5019             \\ 
100 $\times$ 100      & 2.46221E-06        & 4.9852
                      & 1.72697E-06        & 4.9825
                      & 1.76222E-05        & 4.3770
                      & 2.23573E-05        & 3.4614             \\
\hline
\space    &\multicolumn{4}{l}{\cellcolor{gray!35}{WENO-M}}
          &\multicolumn{4}{l}{\cellcolor{gray!35}{LOP-WENO-M}}\\
\cline{2-5}  \cline{6-9}
$N_{x}\times N_{y}$   & $L_{1}$ error      & $L_{1}$ order 
                      & $L_{\infty}$ error & $L_{\infty}$ order
		  			  & $L_{1}$ error      & $L_{1}$ order 
      				  & $L_{\infty}$ error & $L_{\infty}$ order \\
\Xhline{0.65pt}
40 $\times$ 40        & 2.21884E-04        & -
                      & 1.57466E-04        & -
                      & 2.21884E-04        & -
                      & 1.57466E-04        & -                  \\
60 $\times$ 60        & 3.06949E-05        & 4.8785
                      & 2.20005E-05        & 4.8540
                      & 3.06949E-05        & 4.8785
                      & 2.20005E-05        & 4.8540             \\
80 $\times$ 80        & 7.40640E-06        & 4.9421
                      & 5.25840E-06        & 4.9751
                      & 7.40640E-06        & 4.9421
                      & 5.25840E-06        & 4.9751             \\ 
100 $\times$ 100      & 2.44462E-06        & 4.9675
                      & 1.72528E-06        & 4.9943
                      & 2.44462E-06        & 4.9675
                      & 1.72528E-06        & 4.9943             \\
\hline
\space    &\multicolumn{4}{l}{\cellcolor{gray!35}{WENO-PM6}} 
          &\multicolumn{4}{l}{\cellcolor{gray!35}{LOP-WENO-PM6}}\\
\cline{2-5}  \cline{6-9}
$N_{x}\times N_{y}$   & $L_{1}$ error      & $L_{1}$ order 
                      & $L_{\infty}$ error & $L_{\infty}$ order
		  			  & $L_{1}$ error      & $L_{1}$ order 
      				  & $L_{\infty}$ error & $L_{\infty}$ order \\
\Xhline{0.65pt}
40 $\times$ 40        & 2.35238E-04        & -
                      & 1.57970E-04        & -
                      & 2.35238E-04        & -
                      & 1.57970E-04        & -                  \\
60 $\times$ 60        & 3.14340E-05        & 4.9639
                      & 2.18667E-05        & 4.8770
                      & 3.14340E-05        & 4.9639
                      & 2.18667E-05        & 4.8770             \\
80 $\times$ 80        & 7.49935E-06        & 4.9814
                      & 5.25053E-06        & 4.9591
                      & 7.49935E-06        & 4.9814
                      & 5.25053E-06        & 4.9591             \\ 
100 $\times$ 100      & 2.46354E-06        & 4.9888
                      & 1.72711E-06        & 4.9828
                      & 2.46354E-06        & 4.9888
                      & 1.72711E-06        & 4.9828             \\
\hline
\space    &\multicolumn{4}{l}{\cellcolor{gray!35}{WENO-IM(2, 0.1)}} 
    &\multicolumn{4}{l}{\cellcolor{gray!35}{LOP-WENO-IM(2, 0.1)}}\\
\cline{2-5}  \cline{6-9}
$N_{x}\times N_{y}$   & $L_{1}$ error      & $L_{1}$ order 
                      & $L_{\infty}$ error & $L_{\infty}$ order
		  			  & $L_{1}$ error      & $L_{1}$ order 
      				  & $L_{\infty}$ error & $L_{\infty}$ order \\
\Xhline{0.65pt}
40 $\times$ 40        & 2.30237E-04        & -
                      & 1.57910E-04        & -
                      & 2.30237E-04        & -
                      & 1.57910E-04        & -                  \\
60 $\times$ 60        & 3.12478E-05        & 4.9256
                      & 2.18921E-05        & 4.8732
                      & 3.12478E-05        & 4.9256
                      & 2.18921E-05        & 4.8732             \\
80 $\times$ 80        & 7.48097E-06        & 4.9693
                      & 5.25058E-06        & 4.9631
                      & 7.48097E-06        & 4.9693
                      & 5.25058E-06        & 4.9631             \\ 
100 $\times$ 100      & 2.46044E-06        & 4.9834
                      & 1.72680E-06        & 4.9836
                      & 2.46044E-06        & 4.9834
                      & 1.72680E-06        & 4.9836             \\
\hline
\space    &\multicolumn{4}{l}{\cellcolor{gray!35}{WENO-PPM5}} 
          &\multicolumn{4}{l}{\cellcolor{gray!35}{LOP-WENO-PPM5}}\\
\cline{2-5}  \cline{6-9}
$N_{x}\times N_{y}$   & $L_{1}$ error      & $L_{1}$ order 
                      & $L_{\infty}$ error & $L_{\infty}$ order
		  			  & $L_{1}$ error      & $L_{1}$ order 
      				  & $L_{\infty}$ error & $L_{\infty}$ order \\
\Xhline{0.65pt}
40 $\times$ 40        & 2.35717E-04        & -
                      & 1.57956E-04        & -
                      & 2.35717E-04        & -
                      & 1.57956E-04        & -                  \\
60 $\times$ 60        & 3.15372E-05        & 4.9609
                      & 2.18541E-05        & 4.8782
                      & 3.15372E-05        & 4.9609
                      & 2.18541E-05        & 4.8782             \\
80 $\times$ 80        & 7.51693E-06        & 4.9847
                      & 5.25176E-06        & 4.9563
                      & 7.51693E-06        & 4.9847
                      & 5.25176E-06        & 4.9563             \\ 
100 $\times$ 100      & 2.46739E-06        & 4.9923
                      & 1.72749E-06        & 4.9829
                      & 2.46739E-06        & 4.9923
                      & 1.72749E-06        & 4.9829             \\
\hline
\space    &\multicolumn{4}{l}{\cellcolor{gray!35}{WENO-RM(260)}} 
      &\multicolumn{4}{l}{\cellcolor{gray!35}{LOP-WENO-RM(260)}}\\
\cline{2-5}  \cline{6-9}
$N_{x}\times N_{y}$   & $L_{1}$ error      & $L_{1}$ order 
                      & $L_{\infty}$ error & $L_{\infty}$ order
		  			  & $L_{1}$ error      & $L_{1}$ order 
      				  & $L_{\infty}$ error & $L_{\infty}$ order \\
\Xhline{0.65pt}
40 $\times$ 40        & 2.31192E-04        & -
                      & 1.58226E-04        & -
                      & 2.31192E-04        & -
                      & 1.58226E-04        & -                  \\
60 $\times$ 60        & 3.13103E-05        & 4.9309
                      & 2.18799E-05        & 4.8795
                      & 3.13103E-05        & 4.9309
                      & 2.18799E-05        & 4.8795             \\
80 $\times$ 80        & 7.48936E-06        & 4.9724
                      & 5.24972E-06        & 4.9617
                      & 7.48936E-06        & 4.9724
                      & 5.24972E-06        & 4.9617             \\ 
100 $\times$ 100      & 2.46221E-06        & 4.9852
                      & 1.72697E-06        & 4.9825
                      & 2.46221E-06        & 4.9852
                      & 1.72697E-06        & 4.9825             \\
\hline
\space    &\multicolumn{4}{l}{\cellcolor{gray!35}{WENO-ACM}} 
          &\multicolumn{4}{l}{\cellcolor{gray!35}{LOP-WENO-ACM}}\\
\cline{2-5}  \cline{6-9}
$N_{x}\times N_{y}$   & $L_{1}$ error      & $L_{1}$ order 
                      & $L_{\infty}$ error & $L_{\infty}$ order
		  			  & $L_{1}$ error      & $L_{1}$ order 
      				  & $L_{\infty}$ error & $L_{\infty}$ order \\
\Xhline{0.65pt}
40 $\times$ 40        & 2.31214E-04        & -
                      & 1.58230E-04        & -
                      & 2.31214E-04        & -
                      & 1.58230E-04        & -                  \\
60 $\times$ 60        & 3.13106E-05        & 4.9311
                      & 2.18798E-05        & 4.8795
                      & 3.13106E-05        & 4.9311
                      & 2.18798E-05        & 4.8795             \\
80 $\times$ 80        & 7.48937E-06        & 4.9724
                      & 5.24972E-06        & 4.9617
                      & 7.48937E-06        & 4.9724
                      & 5.24972E-06        & 4.9617             \\ 
100 $\times$ 100      & 2.46221E-06        & 4.9852
                      & 1.72697E-06        & 4.9825
                      & 2.46221E-06        & 4.9852
                      & 1.72697E-06        & 4.9825             \\
\hline
\end{tabular*}
\end{myFontSize}
\end{table}

\subsubsection{Shock-vortex interaction problem}
\begin{example}
\rm{The shock-vortex interaction problem is a very favorable 2D test 
case for high-resolution schemes \cite{Shock-vortex_interaction-1,
Shock-vortex_interaction-2,Shock-vortex_interaction-3}. The initial 
condition is given by}
\label{ex:shock-vortex}
\end{example}
\begin{equation*}
\big( \rho, u, v, p \big)(x, y, 0) = \left\{
\begin{aligned}
\begin{array}{ll}
(1, \sqrt{\gamma}, 0, 1), & x < 0.5, \\
\Bigg(\rho_{\mathrm{L}}\bigg( 
\dfrac{\gamma - 1 + (\gamma + 1)p_{\mathrm{R}}}{\gamma + 1 + (\gamma
- 1)p_{\mathrm{R}}} \bigg), u_{\mathrm{L}}\bigg( \dfrac{1 - 
p_{\mathrm{R}}}{\sqrt{\gamma-1 + p_{\mathrm{R}}(\gamma + 1)}}\bigg), 
0, 1.3 \Bigg), & x \geq 0.5. \\
\end{array}
\end{aligned}
\right.
\label{eq:initial_Euer2D:shock-vortex-interaction}
\end{equation*}
The following perturbations is superimposed onto the left state,
\begin{equation*}
\delta \rho = \dfrac{\rho_{\mathrm{L}}^{2}}{(\gamma - 1)
p_{\mathrm{L}}}\delta T, 
\delta u = \epsilon \dfrac{y - y_{\mathrm{c}}}{r_\mathrm{c}}
\mathrm{e}^{\alpha(1-r^{2})}, 
\delta v = - \epsilon \dfrac{x - x_{\mathrm{c}}}{r_\mathrm{c}}
\mathrm{e}^{\alpha(1-r^{2})}, 
\delta p = \dfrac{\gamma \rho_{\mathrm{L}}^{2}}{(\gamma - 1)
\rho_{\mathrm{L}}}\delta T,
\label{eq:Euler2D:shock-vortex-interactions:Perturbations}
\end{equation*}
where $\epsilon = 0.3, r_{\mathrm{c}} = 0.05, \alpha = 0.204, 
x_{\mathrm{c}} = 0.25, y_{\mathrm{c}} = 0.5,
r = \sqrt{((x - x_{\mathrm{c}})^{2} + (y - y_{\mathrm{c}})^{2})/r_{
\mathrm{c}}^{2}}, \delta T = - (\gamma - 1)\epsilon^{2}\mathrm{e}^{2
\alpha (1 - r^{2})}/(4\alpha \gamma)$.
The transmissive boundary condition is used. We compute the solution 
up to two different output times $t = 0.35, 0.60$ by using all 
considered schemes with a uniform mesh size of $800 \times 800$. 
Here, we set the CFL number to be $0.5$.

\begin{figure}[ht]
\centering
  \includegraphics[height=0.31\textwidth]
  {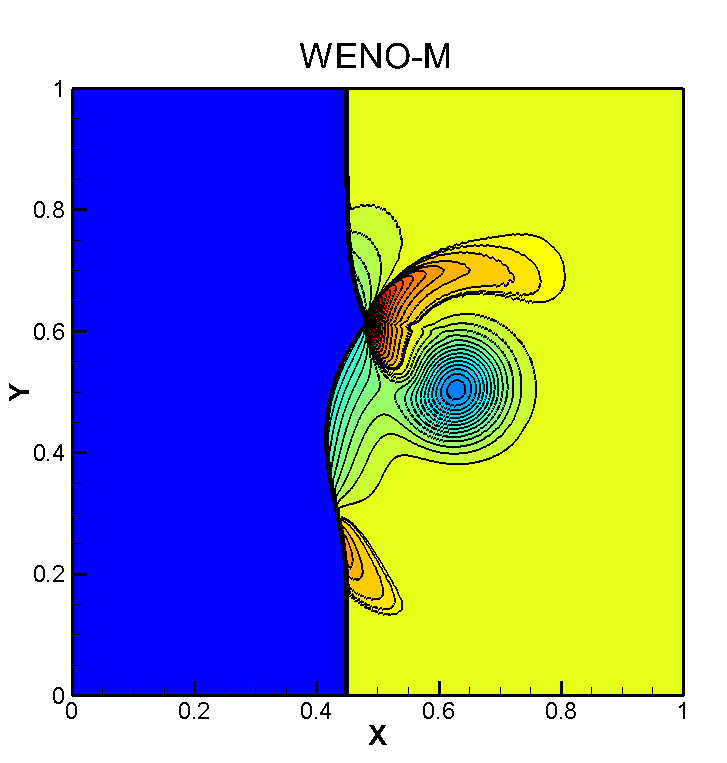}
  \includegraphics[height=0.31\textwidth]
  {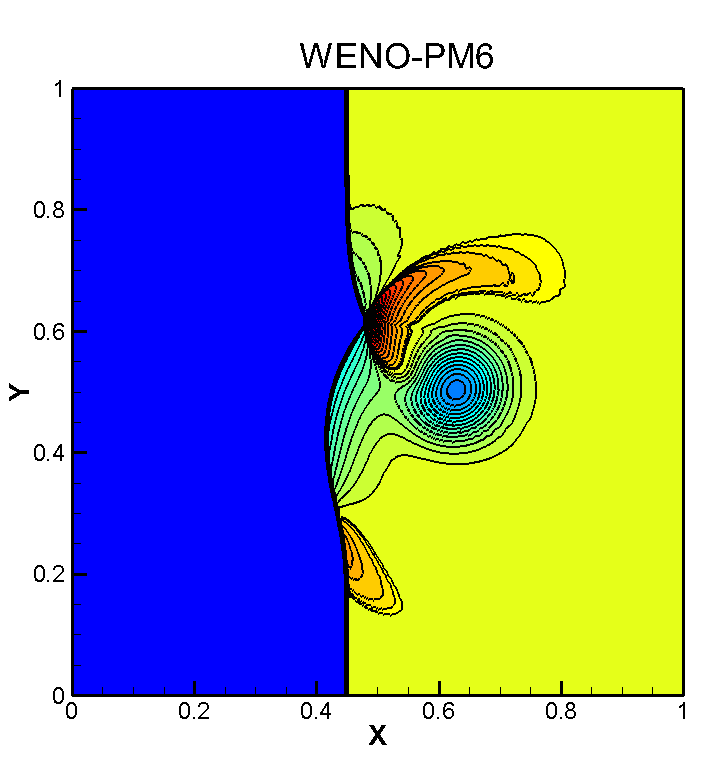}
  \includegraphics[height=0.31\textwidth]
  {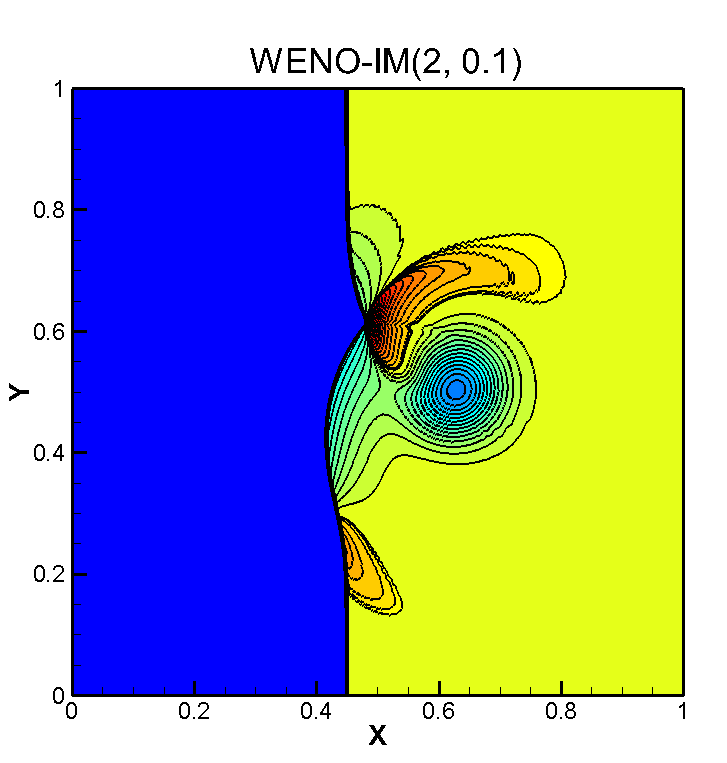}\\
  \includegraphics[height=0.31\textwidth]
  {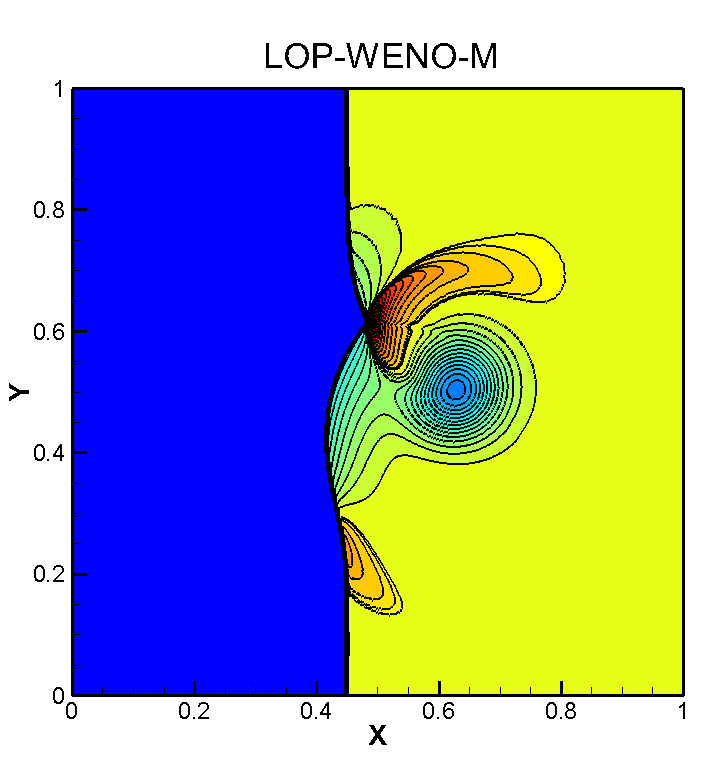}
  \includegraphics[height=0.31\textwidth]
  {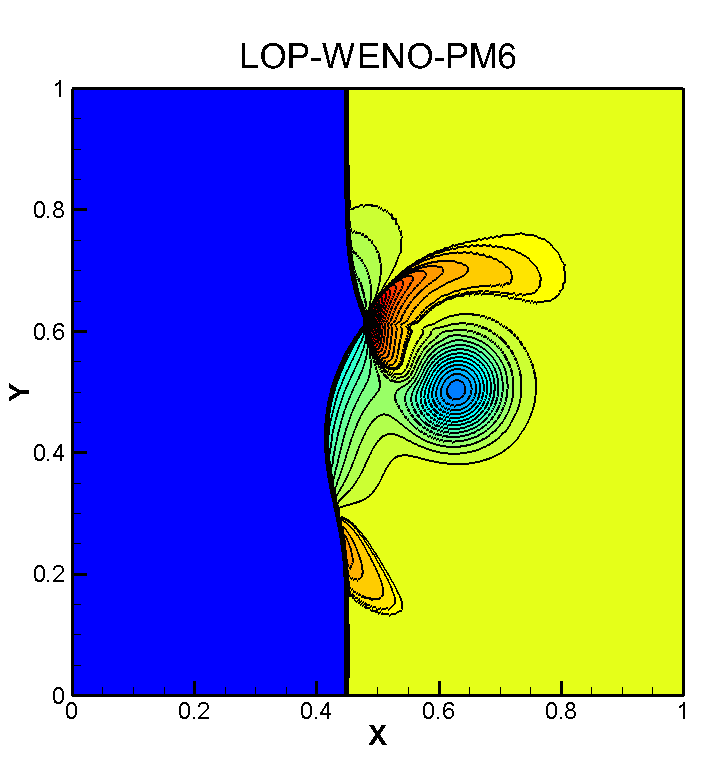}  
  \includegraphics[height=0.31\textwidth]
  {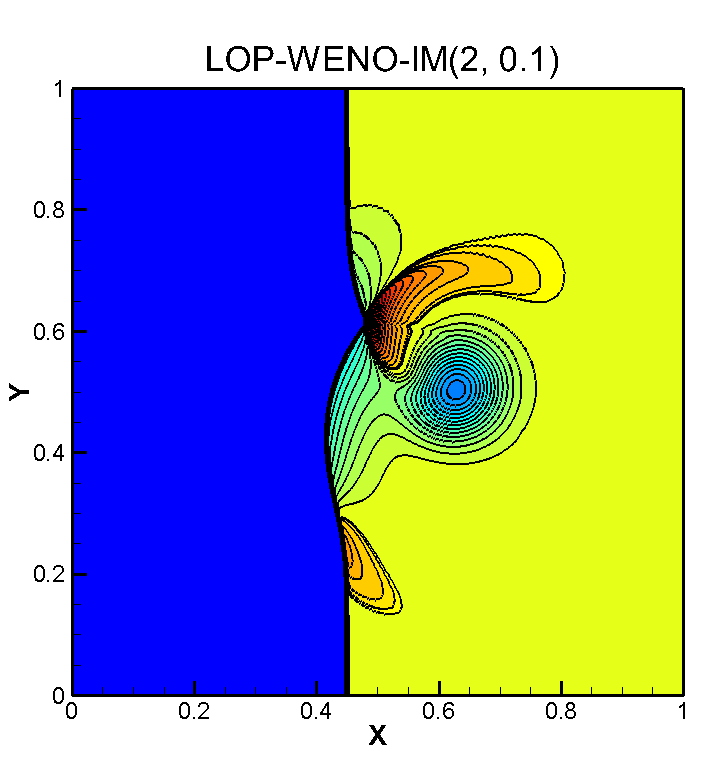}
\caption{Density plots for the Shock-vortex interaction, $t = 0.35$.}
\label{fig:ex:SVI:1}
\end{figure}

\begin{figure}[ht]
\centering
  \includegraphics[height=0.26\textwidth]
  {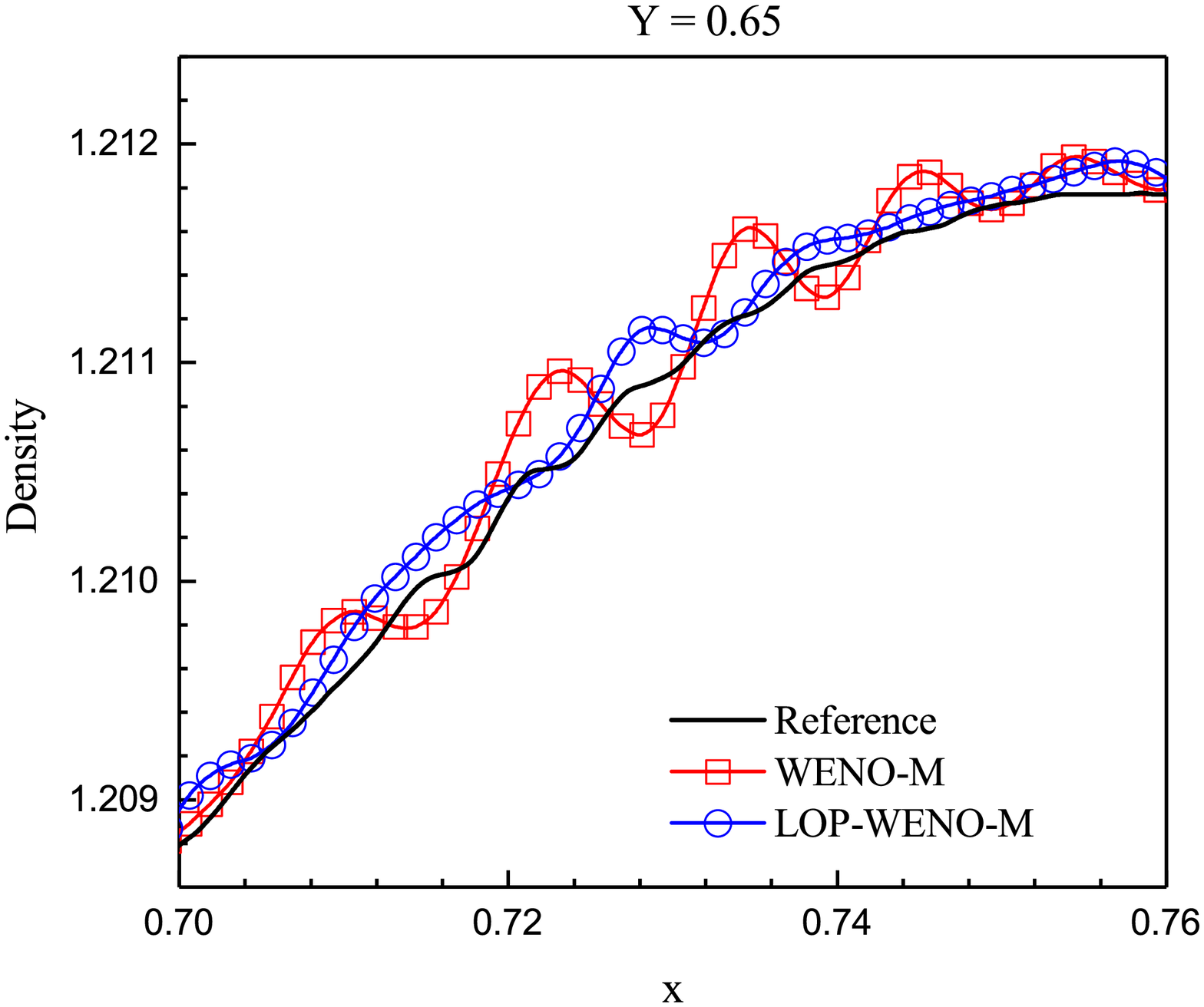}
  \includegraphics[height=0.26\textwidth]
  {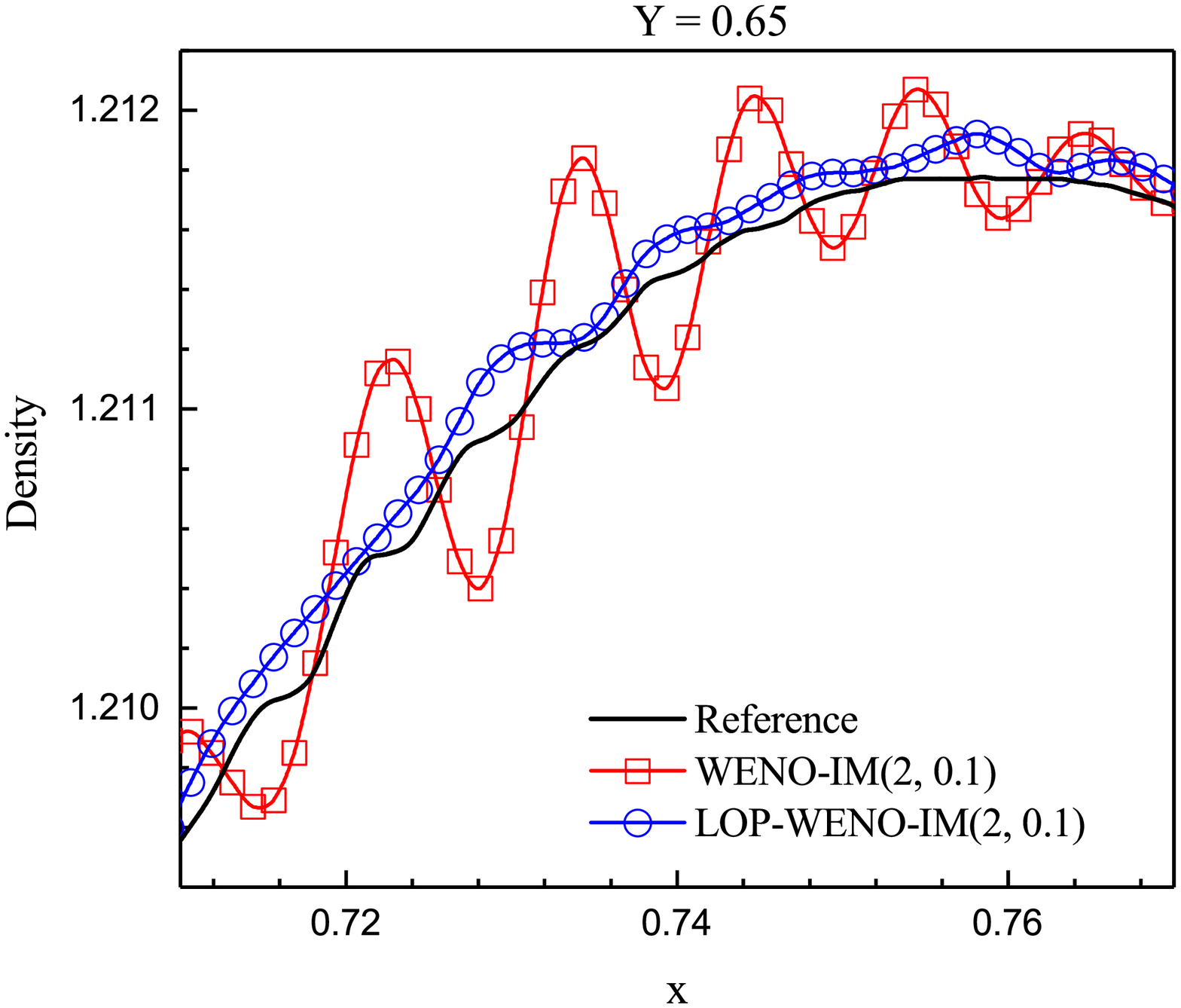}
  \includegraphics[height=0.26\textwidth]
  {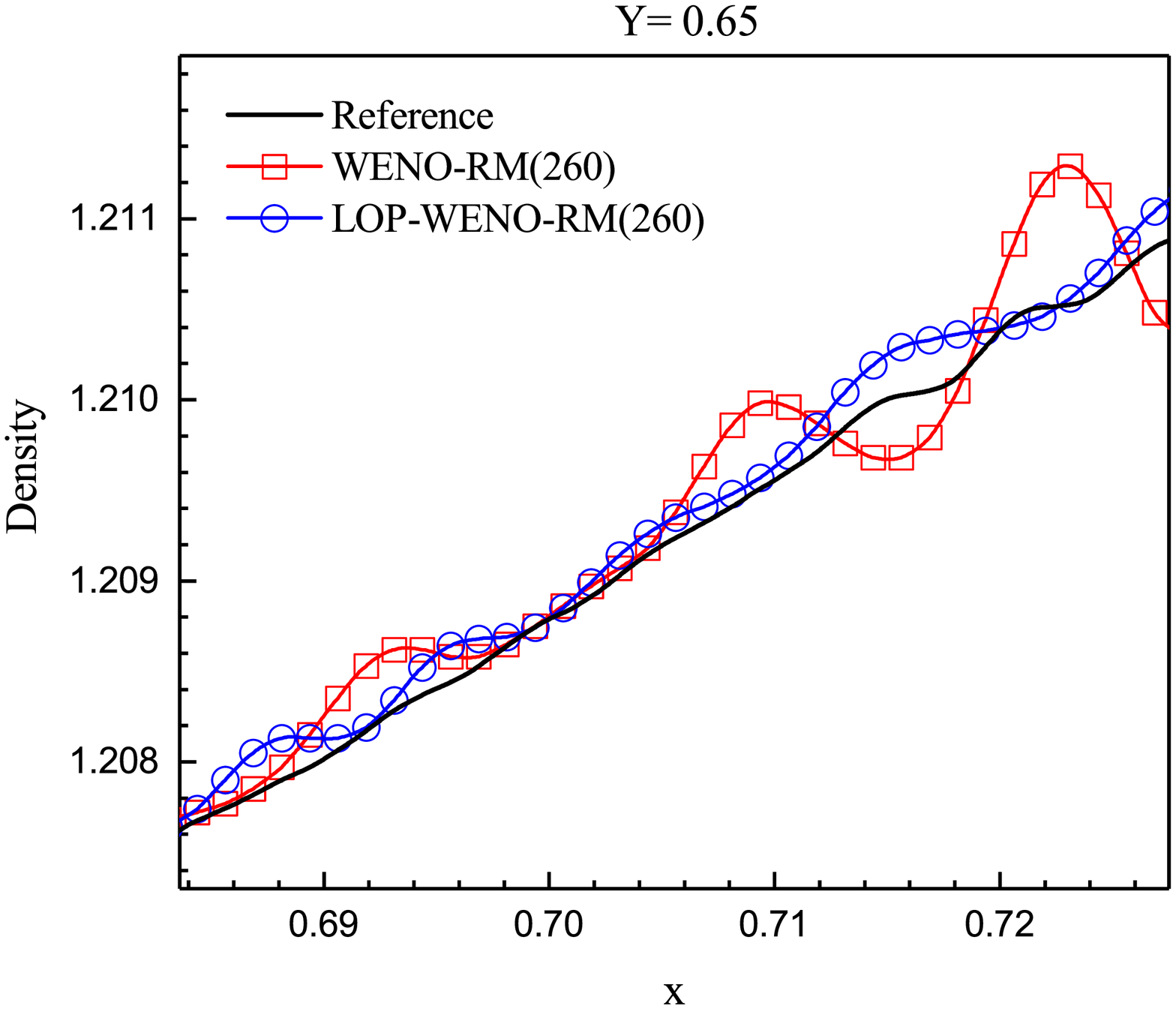}\\
  \includegraphics[height=0.26\textwidth]
  {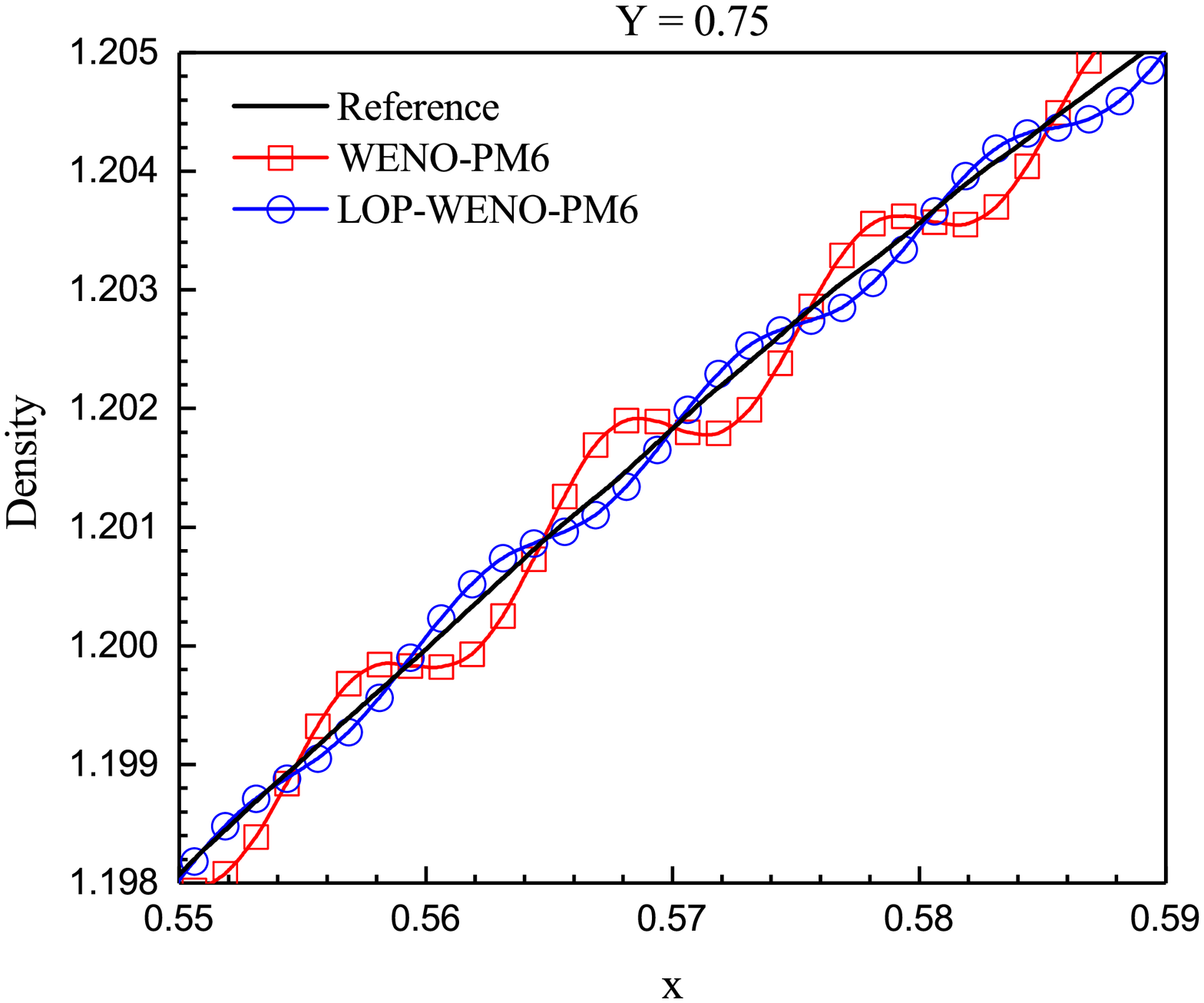}
  \includegraphics[height=0.26\textwidth]
  {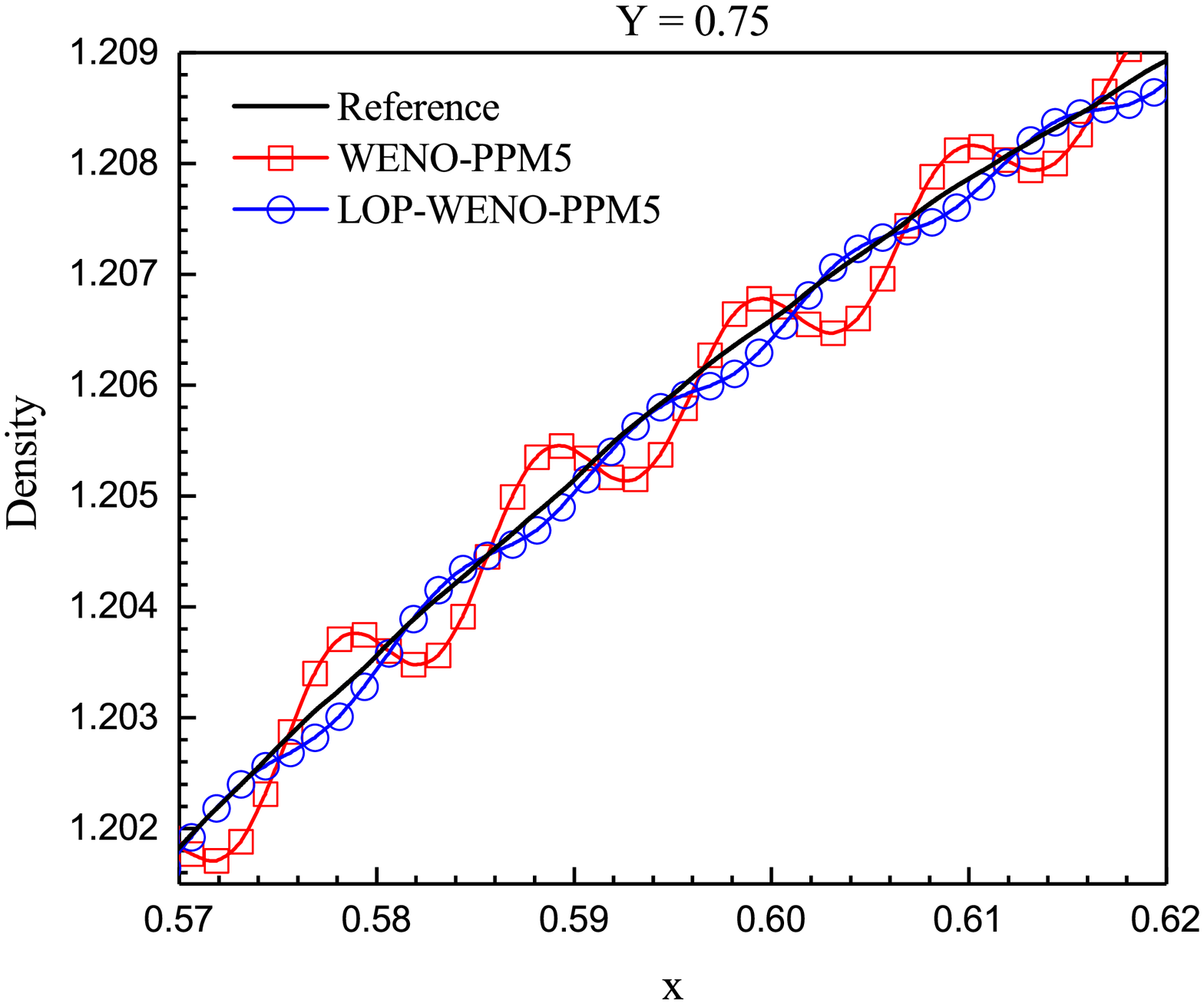}
  \includegraphics[height=0.26\textwidth]
  {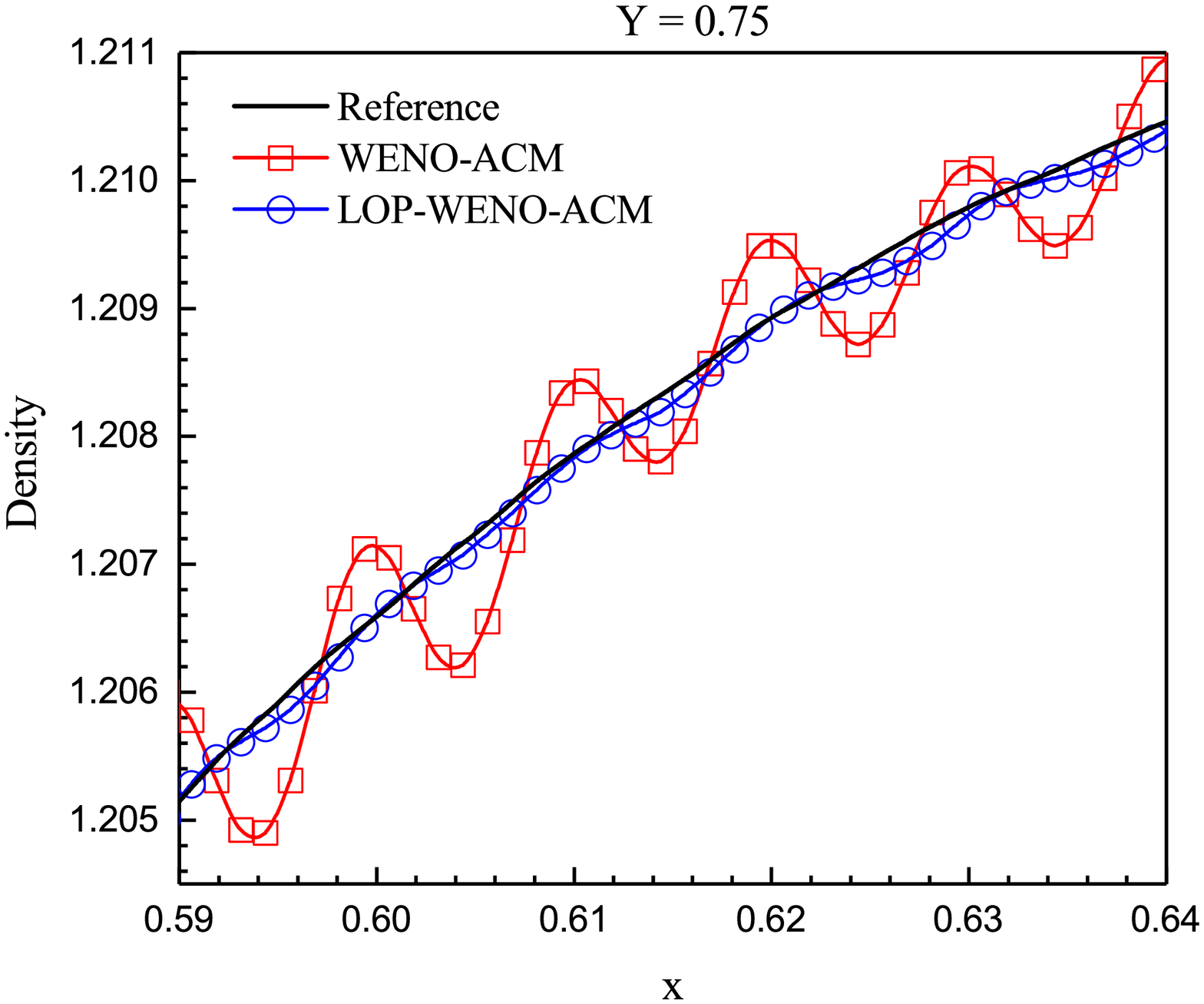}  
\caption{Density cross-sectional slices plotted along the plane
$y = 0.65, 0.75$ with $t = 0.35$.}
\label{fig:ex:SVI:3}
\end{figure}

\begin{figure}[ht]
\centering
  \includegraphics[height=0.31\textwidth]
  {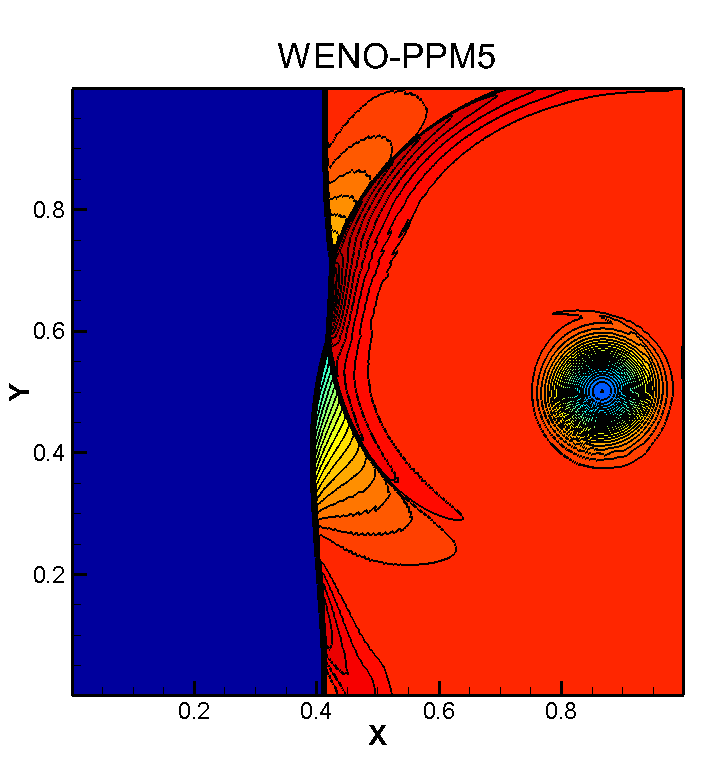}   
  \includegraphics[height=0.31\textwidth]
  {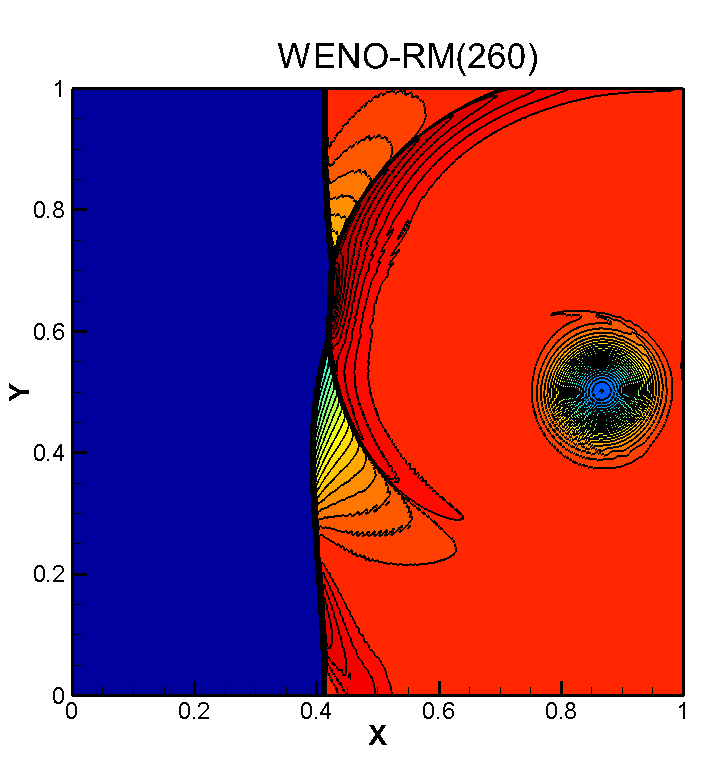} 
  \includegraphics[height=0.31\textwidth]
  {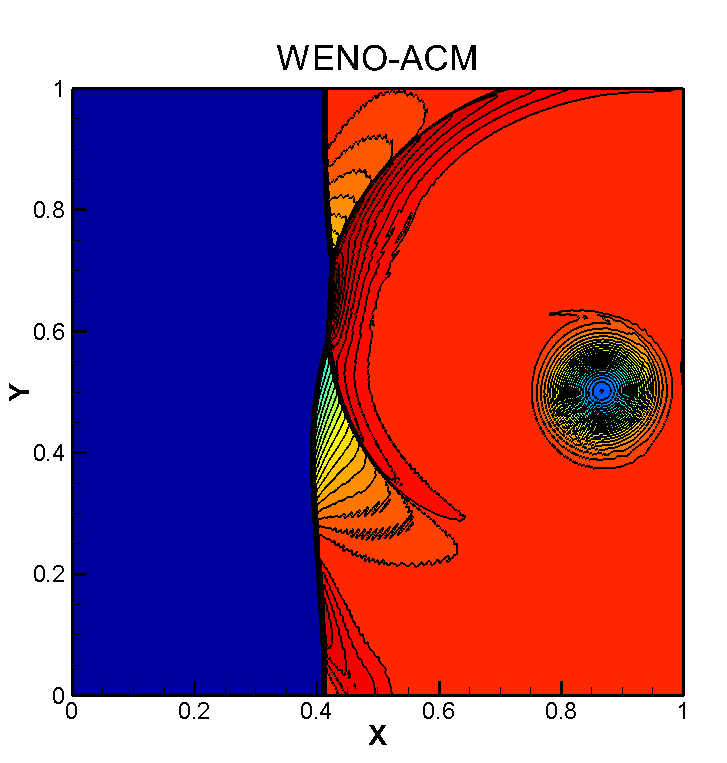} \\
  \includegraphics[height=0.31\textwidth]
  {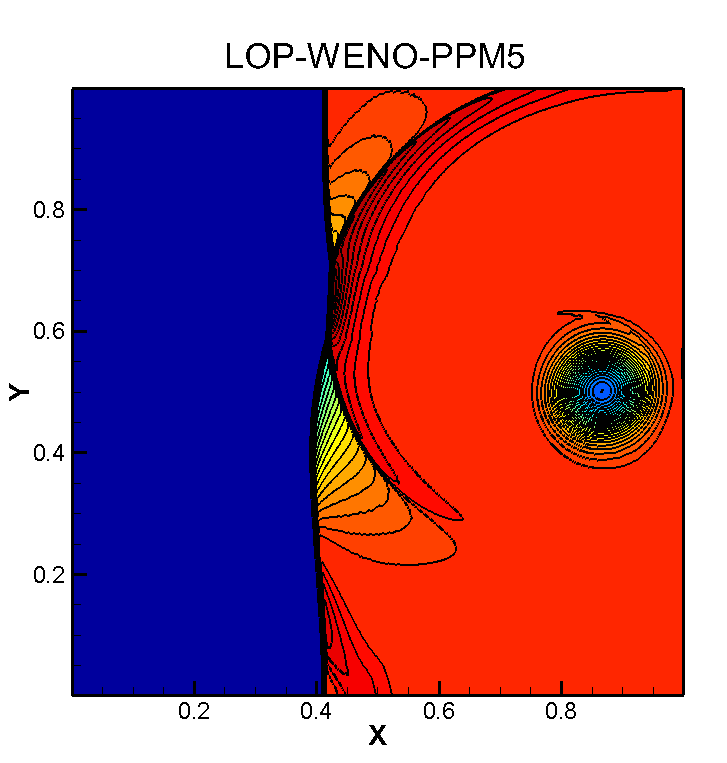}
  \includegraphics[height=0.31\textwidth]
  {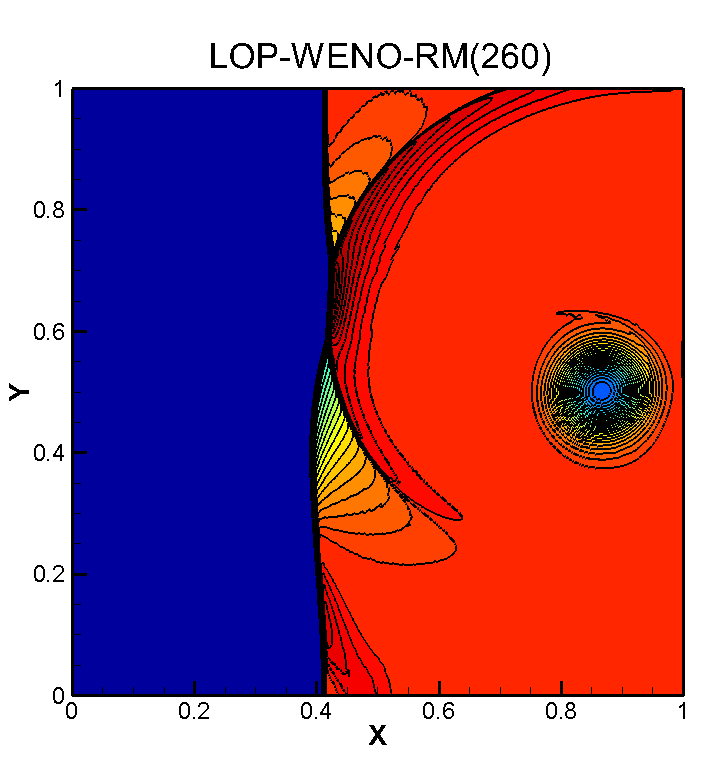}
  \includegraphics[height=0.31\textwidth]
  {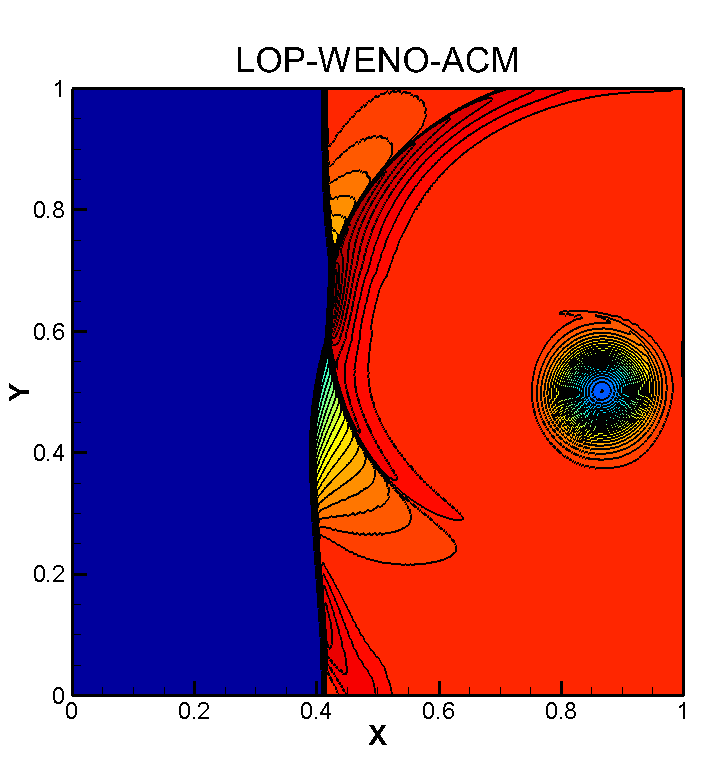}\\   
\caption{Density plots for the Shock-vortex interaction, $t = 0.6$.}
\label{fig:ex:SVI:2}
\end{figure}

\begin{figure}[ht]
\centering
  \includegraphics[height=0.26\textwidth]
  {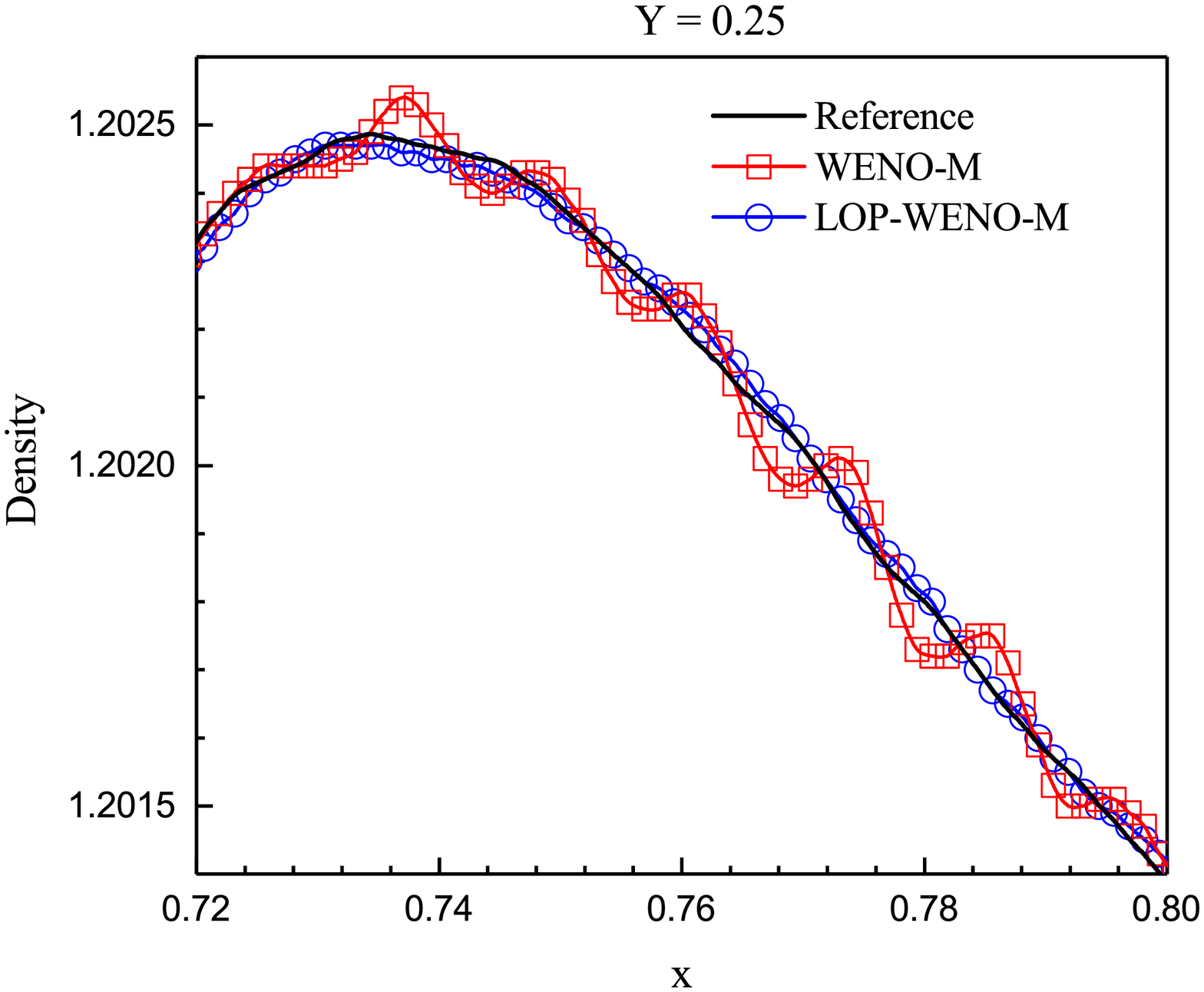}
  \includegraphics[height=0.26\textwidth]
  {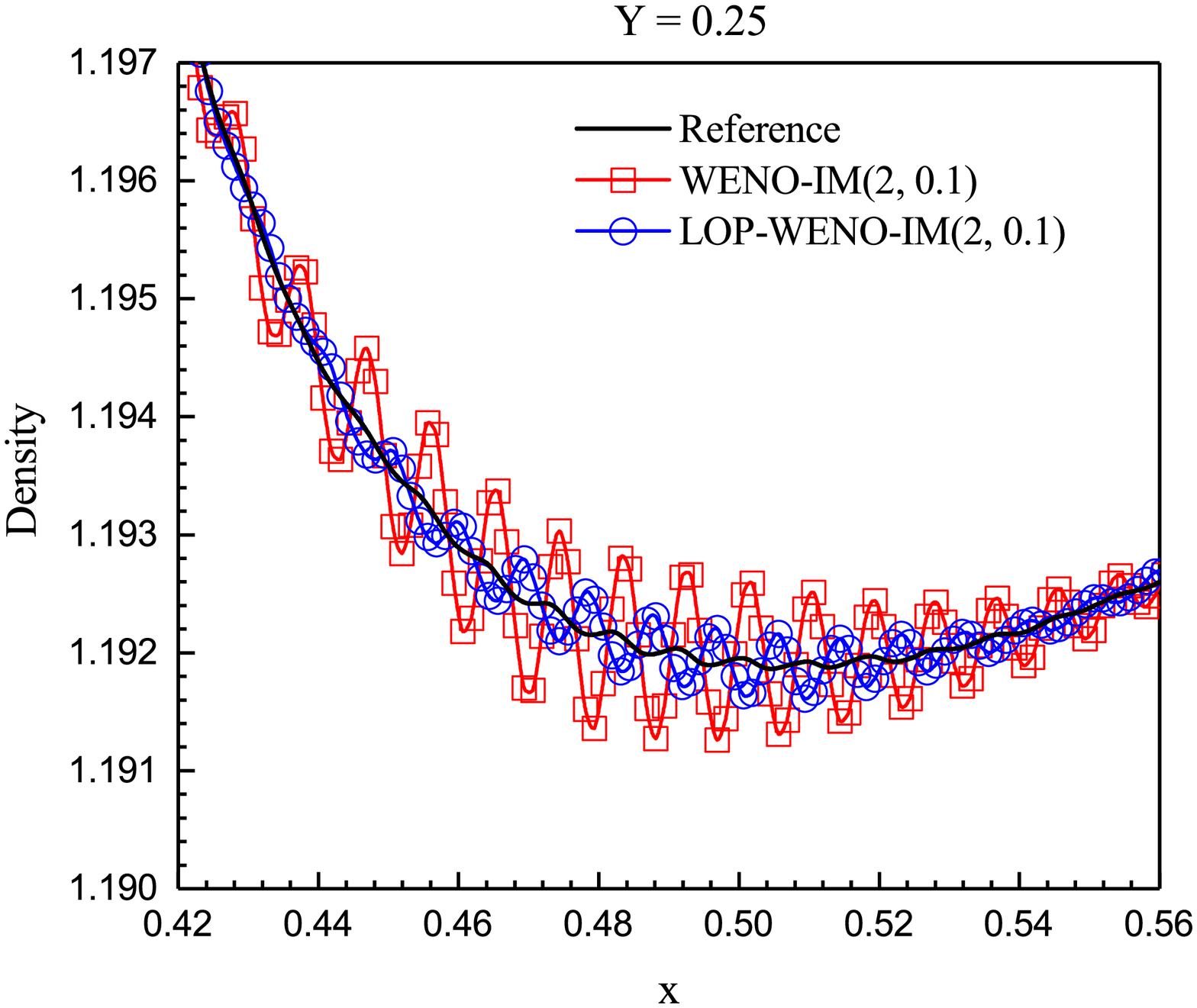}
  \includegraphics[height=0.26\textwidth]
  {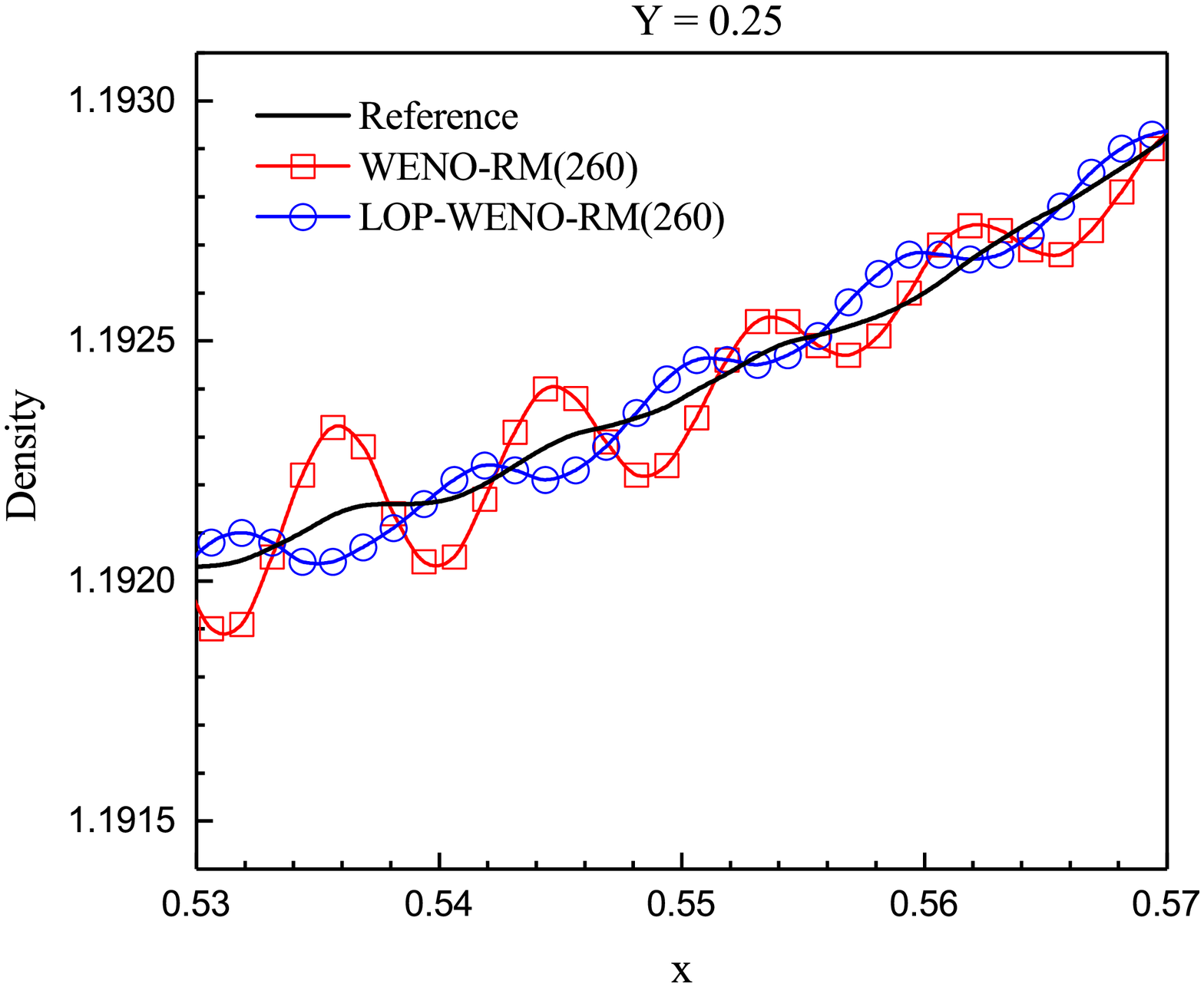}\\
  \includegraphics[height=0.26\textwidth]
  {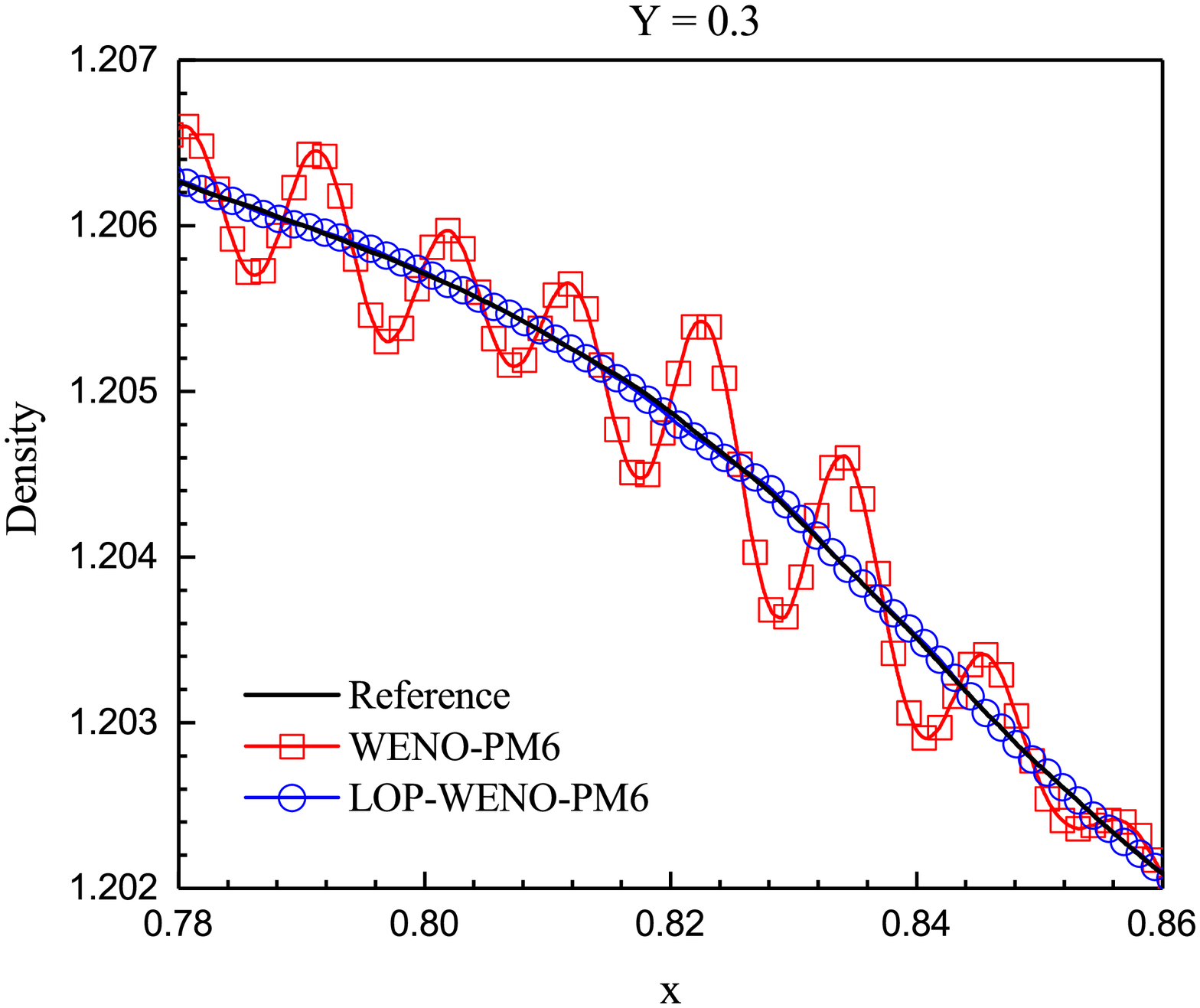}
  \includegraphics[height=0.26\textwidth]
  {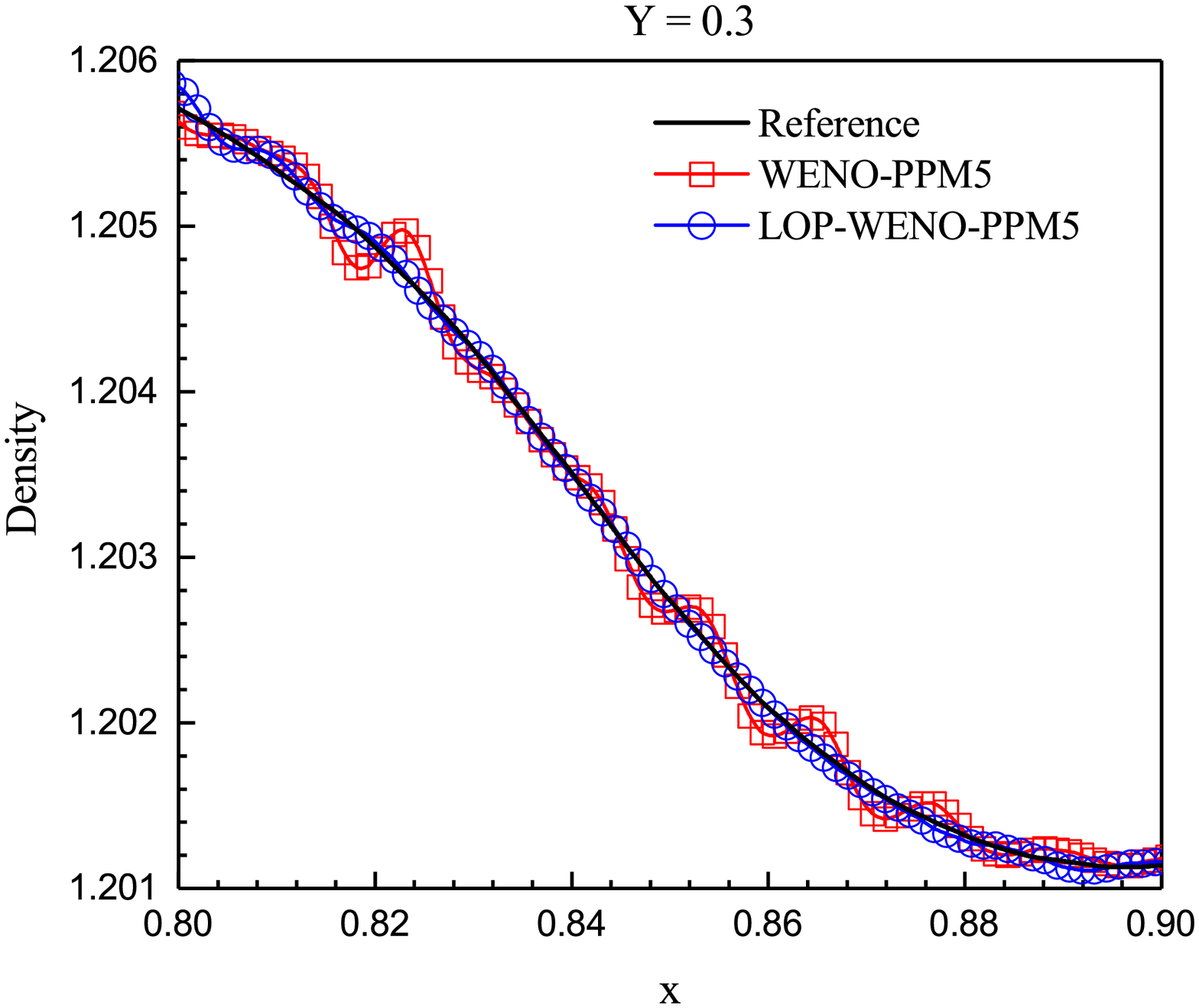}   
  \includegraphics[height=0.26\textwidth]
  {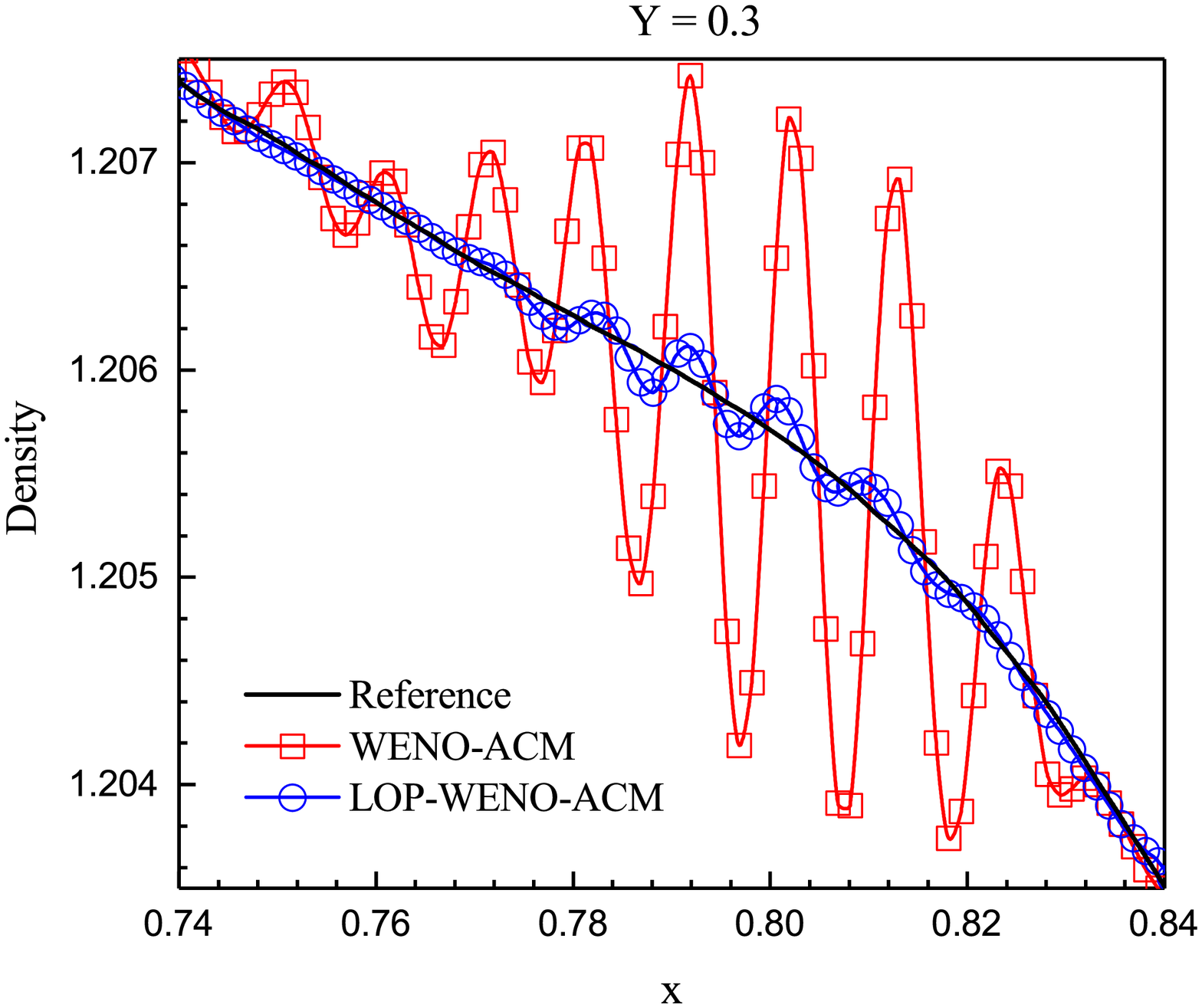}  
\caption{Density cross-sectional slices plotted along the plane
$y = 0.25, 0.30$ with $t = 0.60$.}  
\label{fig:ex:SVI:4}
\end{figure}

Just for the sake of simplicity in presentation, we only show the 
density profiles of the WENO-M, WENO-PM6, WENO-IM(2, 0.1) schemes 
for $t = 0.35$ (see Fig. \ref{fig:ex:SVI:1}), and the density 
profiles of the WENO-PPM5, WENO-RM(260), WENO-ACM schemes for 
$t = 0.6$ (see Fig. \ref{fig:ex:SVI:2}). To unveil the advantage of 
the LOP-WENO-X schemes more precisely, we present the 
cross-sectional slices of density plots along the plane 
$y = 0.65, 0.75$ (see Fig. \ref{fig:ex:SVI:3}) and  $y = 0.25, 0.3$ 
(see Fig. \ref{fig:ex:SVI:4}) of all considered schemes for 
$t = 0.35$ and $t = 0.6$, respectively. It can be seen that: (1) the 
main structure of the shock and vortex after the interaction were 
captured properly by all the considered schemes; (2) in the 
solutions of the WENO-X schemes, clear post-shock oscillations can 
be observed, whereas the post-shock oscillations are considerably 
reduced in the solutions of the associated LOP-WENO-X schemes; (3) 
it is easy to find that the amplitudes of the post-shock 
oscillations produced by the WENO-X schemes are much greater than 
those of their associated LOP-WENO-X schemes. In a word, the 
LOP-WENO-X schemes only produce some highly tolerable post-shock 
oscillations. This should be another merit of the mapped WENO 
schemes with \textit{LOP} mappings.

\section{Conclusions}
\label{secConclusions} 

We aim to develop a method to address the drawback that all the 
MOP-WENO-X schemes proposed in our previous work fail to achieve the 
same resolutions as their associated WENO-X schemes in the region 
with high-frequency smooth waves. To do this, we develop the 
\textit{locally order-preserving (LOP)} mapping in this paper. By 
providing a posteriori adaptive technique, we apply the LOP mapping 
to many previously published mapped WENO schemes. We firstly find 
the global stencil in which the existing mapping is 
\textit{non-order-preserving (non-OP)} through manipulating its 
mapped nonlinear weights of the associated substencils. Then, in 
order to recover the \textit{LOP} property, we abandon these non-OP 
mapped weights and replace them with the weights of the classic 
WENO-JS scheme. We conduct numerical experiments to show that the 
LOP-WENO-X schemes provide similar or even higher resolutions than 
those of their associated WENO-X schemes in the region with 
high-frequency smooth waves. This is the major improvement of the 
LOP-WENO-X schemes. In addition, they can not only preserve high 
resolutions but also prevent spurious oscillations on solving 
problems with high-order critical points or discontinuities, 
especially for long-run simulations. We also find that there should 
be another competitive advancement that, when solving the 2D 
problems with shock waves, the LOP-WENO-X schemes can properly 
capture the main structures of the complicated flows and perform 
admirably in reducing the post-shock oscillations.



\bibliographystyle{model1b-shortjournal-num-names}
\bibliography{../refs}

\begin{thebibliography}{29}
\expandafter\ifx\csname natexlab\endcsname\relax\def\natexlab#1{#1}\fi
\providecommand{\bibinfo}[2]{#2}
\ifx\xfnm\relax \def\xfnm[#1]{\unskip,\space#1}\fi
\bibitem[{Borges et~al.(2008)Borges, Carmona, Costa and Don}]{WENO-Z}
\bibinfo{author}{R.~Borges}, \bibinfo{author}{M.~Carmona},
  \bibinfo{author}{B.~Costa}, \bibinfo{author}{W.S. Don}, \bibinfo{title}{An
  improved weighted essentially non-oscillatory scheme for hyperbolic
  conservation laws}, \bibinfo{shortjournal}{J. Comput. Phys.}
  \bibinfo{volume}{227} (\bibinfo{year}{2008}) \bibinfo{pages}{3191--3211}.
\bibitem[{Chatterjee(1999)}]{Shock-vortex_interaction-1}
\bibinfo{author}{A.~Chatterjee}, \bibinfo{title}{Shock wave deformation in
  shock-vortex interactions}, \bibinfo{shortjournal}{Shock Waves}
  \bibinfo{volume}{9} (\bibinfo{year}{1999}) \bibinfo{pages}{95--105}.
\bibitem[{Feng et~al.(2012)Feng, Hu and Wang}]{WENO-PM}
\bibinfo{author}{H.~Feng}, \bibinfo{author}{F.~Hu}, \bibinfo{author}{R.~Wang},
  \bibinfo{title}{A new mapped weighted essentially non-oscillatory scheme},
  \bibinfo{shortjournal}{J. Sci. Comput.} \bibinfo{volume}{51}
  (\bibinfo{year}{2012}) \bibinfo{pages}{449--473}.
\bibitem[{Feng et~al.(2014)Feng, Huang and Wang}]{WENO-IM}
\bibinfo{author}{H.~Feng}, \bibinfo{author}{C.~Huang},
  \bibinfo{author}{R.~Wang}, \bibinfo{title}{An improved mapped weighted
  essentially non-oscillatory scheme}, \bibinfo{shortjournal}{Appl. Math.
  Comput.} \bibinfo{volume}{232} (\bibinfo{year}{2014})
  \bibinfo{pages}{453--468}.
\bibitem[{Gottlieb and Shu(1998)}]{SSPRK1998}
\bibinfo{author}{S.~Gottlieb}, \bibinfo{author}{C.W. Shu},
  \bibinfo{title}{Total variation diminishing {R}unge-{K}utta schemes},
  \bibinfo{shortjournal}{Math. Comput.} \bibinfo{volume}{67}
  (\bibinfo{year}{1998}) \bibinfo{pages}{73--85}.
\bibitem[{Gottlieb et~al.(2001)Gottlieb, Shu and Tadmor}]{SSPRK2001}
\bibinfo{author}{S.~Gottlieb}, \bibinfo{author}{C.W. Shu},
  \bibinfo{author}{E.~Tadmor}, \bibinfo{title}{Strong stability-preserving
  high-order time discretization methods}, \bibinfo{shortjournal}{SIAM Rev.}
  \bibinfo{volume}{43} (\bibinfo{year}{2001}) \bibinfo{pages}{89--112}.
\bibitem[{Harten(1989)}]{ENO1987JCP83}
\bibinfo{author}{A.~Harten}, \bibinfo{title}{{ENO} schemes with subcell
  resolution}, \bibinfo{shortjournal}{J. Comput. Phys.} \bibinfo{volume}{83}
  (\bibinfo{year}{1989}) \bibinfo{pages}{148--184}.
\bibitem[{Harten et~al.(1987)Harten, Engquist, Osher and
  Chakravarthy}]{ENO1987JCP71}
\bibinfo{author}{A.~Harten}, \bibinfo{author}{B.~Engquist},
  \bibinfo{author}{S.~Osher}, \bibinfo{author}{S.R. Chakravarthy},
  \bibinfo{title}{Uniformly high order accurate essentially non-oscillatory
  schemes {III}}, \bibinfo{shortjournal}{J. Comput. Phys.} \bibinfo{volume}{71}
  (\bibinfo{year}{1987}) \bibinfo{pages}{231--303}.
\bibitem[{Harten and Osher(1987)}]{ENO1987V24}
\bibinfo{author}{A.~Harten}, \bibinfo{author}{S.~Osher},
  \bibinfo{title}{Uniformly high order accurate essentially non-oscillatory
  schemes {I}}, \bibinfo{shortjournal}{SIAM J. Numer. Anal.}
  \bibinfo{volume}{24} (\bibinfo{year}{1987}) \bibinfo{pages}{279--309}.
\bibitem[{Harten et~al.(1986)Harten, Osher, Engquist and
  Chakravarthy}]{ENO1986}
\bibinfo{author}{A.~Harten}, \bibinfo{author}{S.~Osher},
  \bibinfo{author}{B.~Engquist}, \bibinfo{author}{S.R. Chakravarthy},
  \bibinfo{title}{Some results on uniformly high order accurate essentially
  non-oscillatory schemes}, \bibinfo{shortjournal}{Appl. Numer. Math.}
  \bibinfo{volume}{2} (\bibinfo{year}{1986}) \bibinfo{pages}{347--377}.
\bibitem[{Henrick et~al.(2005)Henrick, Aslam and Powers}]{WENO-M}
\bibinfo{author}{A.K. Henrick}, \bibinfo{author}{T.D. Aslam},
  \bibinfo{author}{J.M. Powers}, \bibinfo{title}{Mapped weighted essentially
  non-oscillatory schemes: Achieving optimal order near critical points},
  \bibinfo{shortjournal}{J. Comput. Phys.} \bibinfo{volume}{207}
  (\bibinfo{year}{2005}) \bibinfo{pages}{542--567}.
\bibitem[{Jiang and Shu(1996)}]{WENO-JS}
\bibinfo{author}{G.S. Jiang}, \bibinfo{author}{C.W. Shu},
  \bibinfo{title}{Efficient implementation of weighted {ENO} schemes},
  \bibinfo{shortjournal}{J. Comput. Phys.} \bibinfo{volume}{126}
  (\bibinfo{year}{1996}) \bibinfo{pages}{202--228}.
\bibitem[{Jiang et~al.(2013)Jiang, Shu and Zhang}]{AccuracyTest-Euler2D}
\bibinfo{author}{Y.~Jiang}, \bibinfo{author}{C.W. Shu},
  \bibinfo{author}{M.~Zhang}, \bibinfo{title}{An alternative formulation of
  finite difference weighted {ENO} schemes with {L}ax-{W}endroff time
  discretization for conservation laws}, \bibinfo{shortjournal}{SIAM J. Sci.
  Comput.} \bibinfo{volume}{35} (\bibinfo{year}{2013})
  \bibinfo{pages}{A1137--A1160}.
\bibitem[{Li et~al.(2015)Li, Liu and Zhang}]{WENO-PPM5}
\bibinfo{author}{Q.~Li}, \bibinfo{author}{P.~Liu}, \bibinfo{author}{H.~Zhang},
  \bibinfo{title}{Piecewise {Polynomial} {Mapping} {Method} and {Corresponding}
  {WENO} {Scheme} with {Improved} {Resolution}}, \bibinfo{shortjournal}{Commun.
  Comput. Phys.} \bibinfo{volume}{18} (\bibinfo{year}{2015})
  \bibinfo{pages}{1417--1444}.
\bibitem[{Li and Zhong(2021{\natexlab{a}})}]{WENO-ACM}
\bibinfo{author}{R.~Li}, \bibinfo{author}{W.~Zhong}, \bibinfo{title}{An
  efficient mapped {WENO} scheme using approximate constant mapping},
  \bibinfo{shortjournal}{Numer. Math. Theor. Meth. Appl.}
  (\bibinfo{year}{2021}{\natexlab{a}}) \bibinfo{pages}{Accepted for
  publication}.
\bibitem[{Li and Zhong(2021{\natexlab{b}})}]{WENO-MAIMi}
\bibinfo{author}{R.~Li}, \bibinfo{author}{W.~Zhong}, \bibinfo{title}{A modified
  adaptive improved mapped {WENO} method}, \bibinfo{shortjournal}{Commun.
  Comput. Phys.}  (\bibinfo{year}{2021}{\natexlab{b}}) \bibinfo{pages}{Accepted
  for publication}.
\bibitem[{Li and Zhong(2021{\natexlab{c}})}]{MOP-WENO-ACMk}
\bibinfo{author}{R.~Li}, \bibinfo{author}{W.~Zhong}, \bibinfo{title}{A new
  mapped {WENO} scheme using order-preserving mapping},
  \bibinfo{shortjournal}{Commun. Comput. Phys.}
  (\bibinfo{year}{2021}{\natexlab{c}}) \bibinfo{pages}{Accepted for
  publication}.
\bibitem[{Li and Zhong(2021{\natexlab{d}})}]{MOP-WENO-X}
\bibinfo{author}{R.~Li}, \bibinfo{author}{W.~Zhong}, \bibinfo{title}{Towards
  building the {OP-Mapped} {WENO} schemes: {A} general methodology},
  \bibinfo{shortjournal}{Math. Comput. Appl.} \bibinfo{volume}{26}
  (\bibinfo{year}{2021}{\natexlab{d}}) \bibinfo{pages}{67}.
\bibitem[{Liu et~al.(1994)Liu, Osher and Chan}]{WENO-LiuXD}
\bibinfo{author}{X.D. Liu}, \bibinfo{author}{S.~Osher},
  \bibinfo{author}{T.~Chan}, \bibinfo{title}{Weighted essentially
  non-oscillatory schemes}, \bibinfo{shortjournal}{J. Comput. Phys.}
  \bibinfo{volume}{115} (\bibinfo{year}{1994}) \bibinfo{pages}{200--212}.
\bibitem[{Pao and Salas(1981)}]{Shock-vortex_interaction-2}
\bibinfo{author}{S.P. Pao}, \bibinfo{author}{M.D. Salas}, \bibinfo{title}{A
  numerical study of two-dimensional shock-vortex interaction}, in:
  \bibinfo{booktitle}{AIAA 14th Fluid and Plasma Dynamics Conference},
  \bibinfo{address}{California, Palo Alto, 1981}.
\bibitem[{Ren et~al.(2003)Ren, Liu and Zhang}]{Shock-vortex_interaction-3}
\bibinfo{author}{Y.X. Ren}, \bibinfo{author}{M.~Liu},
  \bibinfo{author}{H.~Zhang}, \bibinfo{title}{A characteristic-wise hybrid
  compact-{WENO} scheme for solving hyperbolic conservation laws},
  \bibinfo{shortjournal}{J. Comput. Phys.} \bibinfo{volume}{192}
  (\bibinfo{year}{2003}) \bibinfo{pages}{365--386}.
\bibitem[{Shu(1998)}]{WENOoverview}
\bibinfo{author}{C.W. Shu}, \bibinfo{title}{Essentially non-oscillatory and
  weighted essentially non-oscillatory schemes for hyperbolic conservation
  laws}, in: \bibinfo{booktitle}{{Advanced Numerical Approximation of Nonlinear
  Hyperbolic Equations. Lecture Notes in Mathematics}}, volume
  \bibinfo{volume}{1697}, \bibinfo{publisher}{Springer},
  \bibinfo{address}{Berlin}, \bibinfo{year}{1998}, pp.
  \bibinfo{pages}{325--432}.
\bibitem[{Shu and Osher(1988)}]{ENO-Shu1988}
\bibinfo{author}{C.W. Shu}, \bibinfo{author}{S.~Osher},
  \bibinfo{title}{Efficient implementation of essentially non-oscillatory
  shock-capturing schemes}, \bibinfo{shortjournal}{J. Comput. Phys.}
  \bibinfo{volume}{77} (\bibinfo{year}{1988}) \bibinfo{pages}{439--471}.
\bibitem[{Shu and Osher(1989)}]{ENO-Shu1989}
\bibinfo{author}{C.W. Shu}, \bibinfo{author}{S.~Osher},
  \bibinfo{title}{Efficient implementation of essentially non-oscillatory
  shock-capturing schemes {II}}, \bibinfo{shortjournal}{J. Comput. Phys.}
  \bibinfo{volume}{83} (\bibinfo{year}{1989}) \bibinfo{pages}{32--78}.
\bibitem[{Titarev and Toro(2004)}]{Titarev-Toro-1}
\bibinfo{author}{V.~Titarev}, \bibinfo{author}{E.~Toro},
  \bibinfo{title}{Finite-volume {WENO} schemes for three-dimensional
  conservation laws}, \bibinfo{shortjournal}{J. Comput. Phys.}
  \bibinfo{volume}{201} (\bibinfo{year}{2004}) \bibinfo{pages}{238--260}.
\bibitem[{Titarev and Toro(2005)}]{Titarev-Toro-3}
\bibinfo{author}{V.~Titarev}, \bibinfo{author}{E.~Toro}, \bibinfo{title}{{WENO}
  schemes based on upwind and centred {TVD} fluxes},
  \bibinfo{shortjournal}{Comput. Fluids} \bibinfo{volume}{34}
  (\bibinfo{year}{2005}) \bibinfo{pages}{705--720}.
\bibitem[{Toro and Titarev(2005)}]{Titarev-Toro-2}
\bibinfo{author}{E.~Toro}, \bibinfo{author}{V.~Titarev}, \bibinfo{title}{{TVD}
  {Fluxes} for the {High}-{Order} {ADER} {Schemes}}, \bibinfo{shortjournal}{J.
  Sci. Comput.} \bibinfo{volume}{24} (\bibinfo{year}{2005})
  \bibinfo{pages}{285--309}.
\bibitem[{Wang et~al.(2016)Wang, Feng and Huang}]{WENO-RM260}
\bibinfo{author}{R.~Wang}, \bibinfo{author}{H.~Feng},
  \bibinfo{author}{C.~Huang}, \bibinfo{title}{A {New} {Mapped} {Weighted}
  {Essentially} {Non-oscillatory} {Method} {Using} {Rational} {Function}},
  \bibinfo{shortjournal}{J. Sci. Comput.} \bibinfo{volume}{67}
  (\bibinfo{year}{2016}) \bibinfo{pages}{540--580}.
\bibitem[{Zhang et~al.(2011)Zhang, Zhang and Shu}]{FVMaccuracyProofs03}
\bibinfo{author}{R.~Zhang}, \bibinfo{author}{M.~Zhang}, \bibinfo{author}{C.W.
  Shu}, \bibinfo{title}{On the order of accuracy and numerical performance of
  two classes of finite volume {WENO} schemes}, \bibinfo{shortjournal}{Commun.
  Comput. Phys.} \bibinfo{volume}{9} (\bibinfo{year}{2011})
  \bibinfo{pages}{807--827}.

\end{thebibliography}

\end{document}